\documentclass[11pt]{article}
\usepackage[a4paper, total={7.2in, 10.35in}]{geometry}

\usepackage{amssymb,amsmath,amsfonts,mathrsfs,dsfont,stmaryrd}
\usepackage{etoolbox,verbatim,color}
\usepackage{url}
\usepackage{algorithmic}
\usepackage[ruled,vlined]{algorithm2e}

\DeclareMathOperator*{\argmin}{argmin}

\newcommand{\norm}[1]{\left\lVert {#1} \right\rVert}

\usepackage[table]{xcolor}
\usepackage{caption}     % for sub figures
\usepackage{subcaption}  % for sub figures
\usepackage{mathtools}   									   % for table
\newcounter{tableeqn}[table]    							   % for table
\newcounter{tablesubeqn}[tableeqn]							   % for table
\newcommand{\kron}{\boxtimes}  % produit de Kronecker
\makeatletter
\newcommand{\bigotimesa}{\mathop{{\bigotimes}_a}}

\newcommand{\tens}[1]{\mathcal{#1}}
\newcommand{\tG}{\tens{G}}
\newcommand{\tI}{\tens{I}}
\newcommand{\tT}{\tens{T}}
\newcommand{\tX}{\tens{X}}

\definecolor{gris}{gray}{0.90}
\definecolor{gris25}{gray}{0.90}
\definecolor{americanrose}{rgb}{1.0, 0.01, 0.24}
\definecolor{bostonuniversityred}{rgb}{0.8, 0.0, 0.0}
\definecolor{shamrockgreen}{rgb}{0.0, 0.62, 0.38}
\definecolor{selectiveyellow}{rgb}{1.0, 0.73, 0.0}
\definecolor{royalblue}{rgb}{0.25, 0.41, 0.88}
\definecolor{ashgrey}{rgb}{0.7, 0.75, 0.71}
\definecolor{burgundy}{RGB}{159,29,53}
\definecolor{darkgreen}{RGB}{18,53,26}
\definecolor{lightblue}{RGB}{102,217,255}
\definecolor{fakeorange}{RGB}{255,140,102}
\definecolor{arylideyellow}{rgb}{0.91, 0.84, 0.42}
\definecolor{bananayellow}{rgb}{1.0, 0.88, 0.21}
\definecolor{gris_f}{gray}{0.35}
\definecolor{bordure}{rgb}{0.09,0.17,0.68}
\definecolor{aquamarine}{rgb}{0.5, 1.0, 0.83}
\definecolor{apricot}{rgb}{0.98, 0.81, 0.69}
\definecolor{babyblue}{rgb}{0.54, 0.81, 0.94}
\definecolor{uipoppy}{RGB}{225, 64, 5}
\definecolor{uipaleblue}{RGB}{96,123,139}
\definecolor{uiblack}{RGB}{0, 0, 0}
\definecolor{decoda}{RGB}{0,153, 0}
\definecolor{lightgreen}{rgb}{0.56, 0.93, 0.56}
\definecolor{blue_f}{rgb}{0.2, 0.2, 0.6}
\definecolor{cinnamon}{rgb}{0.82, 0.41, 0.12}
\definecolor{darkpastelgreen}{rgb}{0.2, 0.75, 0.24}
\definecolor{drab}{rgb}{0.59, 0.44, 0.09}
\definecolor{anrblue}{rgb}{0.1, 0.1, 0.5}

\usepackage{xcolor}

\definecolor{brightpink}{rgb}{1.0, 0.0, 0.5}

\begin{document}
\title{Accelerating Block Coordinate Descent for \\ Nonnegative Tensor Factorization}
\author{
Andersen Man Shun Ang$^{\dagger}$ 
\quad Jeremy E. Cohen$^{\ddagger}$ 
\quad Nicolas Gillis$^{\dagger}$ 
\quad Le Thi Khanh Hien$^{\dagger}$\thanks{This work was supported by the Fonds de la Recherche Scientifique - FNRS and the Fonds Wetenschappelijk Onderzoek - Vlanderen (FWO) under EOS Project no O005318F-RG47, and by the European Research Council (ERC starting grant no 679515). 
E-mails:\{manshun.ang, nicolas.gillis, thikhanhhien\}@umons.ac.be, jeremy.cohen@irisa.fr
} \\~\\ 
$\dagger$~Department of Mathematics and Operational Research \\
Facult\'e Polytechnique, Universit\'e de Mons \\
Rue de Houdain 9, 7000 Mons, Belgium \\~\\
$\ddagger$~Universit\'{e} de Rennes, INRIA, CNRS, IRISA, France
}
\date{}

\maketitle
\begin{abstract}
%In this paper, we consider the approximate nonnegative tensor canonical polyadic decomposition problem (NTF),
This paper is concerned with improving the empirical convergence speed of block-coordinate descent algorithms for approximate nonnegative tensor factorization (NTF).
We propose an extrapolation strategy in-between block updates, referred to as heuristic extrapolation with restarts (HER). HER significantly accelerates the empirical convergence speed of most existing  block-coordinate algorithms for dense NTF,
in particular for challenging computational scenarios,
while requiring a negligible additional computational budget.
\end{abstract}


\section{Introduction}
\label{sec:intro}
In this paper, we consider the approximate nonnegative tensor canonical polyadic decomposition (CPD) problem, which we refer to as nonnegative tensor factorization (NTF).
A $N$-way array or $N$-th order tensor $\tT$ is a multidimensional array in the product $\mathbb{R}^{I_1\times\ldots\times I_N}$ of the vector spaces $\mathbb{R}^{I_i}$ for $i=1,2,\dots,N$. 
A vector $x\in\mathbb{R}^{I_1}$ is a first-order tensor, and a matrix $M\in\mathbb{R}^{I_1\times I_2}$ is a second-order tensor.
The goal of NTF is to approximate a tensor $\tT$ by a structured tensor $\tX$. Using the squared Frobenius norm as a distance metric, defined as 
$\norm{\mathcal X}_F^2={\sum_{j_1,j_2,\ldots j_N} \mathcal X_{j_1j_2\ldots j_N}^2}$,  NTF is the following optimization problem:
\begin{equation}\label{eq:NTF}
	\underset{a_{p}^{(i)}\geq 0, \; 1 \leq i \leq N,\; 1\leq p \leq r }{\min} \;   
	\norm{
		\mathcal{T}  - \sum_{p=1}^r \bigotimes_{i=1}^{N} a_p^{(i)}
	}_F^2,
\end{equation}
where $\bigotimes$ is a tensor product over $N$ real vector spaces $\mathbb{R}^{I_1}$, $\ldots$, $\mathbb{R}^{I_N}$ defined as follows:
\[
\left[\bigotimes_{i=1}^{N}{a^{(i)}_p}\right]_{j_1,j_2,\ldots,j_N} :=
\prod_{i=1}^{N}a^{(i)}_p(j_i),  \quad \text{ where } a^{(i)}_p \in \mathbb{R}^{I_i} \text{ for } i=1,2,\dots,N \text{ and } p=1,2,\dots,r.
\]
In other words, NTF is a low\footnote{Low means much smaller than the generic rank of tensors in the considered tensor space~\cite{Friedland2012generic,Chiantini2012generic}.} nonnegative rank approximation problem, since by definition any nonnegative rank $r$ tensor $\tX$ of order $N$ can be parameterized as
\[
\tX =
\sum_{p=1}^{r} \bigotimes_{i=1}^{N}
a^{(i)}_p \text{ where }
a^{(i)}_p \in \mathbb{R}_+^{I_i} \text{ for } i=1,2,\dots,N \text{ and } p=1,2,\dots,r.
\]

Intuitively, using NTF to approximate a tensor means using a part-based decomposition to summarize its content with a few ``simple'' rank-one tensors $a^{(1)}_p \otimes a^{(2)}_p \otimes \ldots \otimes a^{(N)}_p$  where the components $a^{(i)}_{p}$ are entrywise nonnegative, with $1 \leq p \leq r$ and $1 \leq i \leq N$.
This idea finds numerous applications in diverse areas, among which chemometrics or psychometrics are historical examples, see  \cite{CC1970,Harshman1970,Kiers2000}.

In this paper, we use the Frobenius norm to measure the error of approximation. It is arguably the most widely used measure, mostly because it has some nice properties (in particular, the subproblems in each block of variables is a convex quadratic problem; see below) and it corresponds to the maximum likelihood estimator in the presence of i.i.d.\ Gaussian noise.
NTF is a non-convex optimization problem.
Moreover, no closed-form solution is known to solve NTF; in fact, the problem is NP-hard already for the matrix case, that is, for $N=2$; see~\cite{vavasis2009complexity}.
Therefore, there has been a large amount of works dedicated to solving NTF using various optimization heuristics; see \S~\ref{sec:survey} for a  review of the state-of-the-art algorithms.
However, note that unlike unconstrained approximate tensor factorization, NTF is well-posed in the sense that there always exists an optimal solution; see \cite{LC2009}.
Moreover, a solution $\tX^{\ast}$ to NTF is almost always\,\footnote{The set of ``bad'' $\tT$ form an hypersurface.} unique for $N>2$, and the solution to~\eqref{eq:NTF}  also has exactly rank $r$; see~\cite{Qi2016Uniqueness}.

%\paragraph{Previous works}
\paragraph{Outline and contribution}
This paper focuses on computing solutions to NTF as fast as possible.  %for a given approximation error tolerance.
We derive new Block-Coordinate Descent (BCD) algorithms for NTF, that aim at being faster than existing BCD algorithms.
To achieve this empirical speed-up in convergence speed, an extrapolation scheme ``\`a la Nesterov'' is used every time a block has been optimized,  before switching to another block. The proposed Heuristic Extrapolation with Restarts (HER) algorithm consists of the following steps:
\begin{enumerate}
	
	\item Initialize $A^{(i)} = [a_1^{(i)},\ldots,a_r^{(i)}]$  and pairing variables $\hat{A}^{(i)}$ for $1 \leq i \leq N$. 
	
	\item Loop over the blocks $A^{(i)}$ ($1 \leq i \leq N$):
	\begin{enumerate}
		
		\item Update $A^{(i)}$ by minimizing~\eqref{eq:NTF} where the other blocks are fixed and take the value of the pairing variables $\hat{A}^{(j)}$ ($j \neq i$). For example, one can take a gradient step (see \S~\ref{sec:AOframework} for more sophisticated strategies).
		Keep the previous value of $A^{(i)}$ in memory as $A^{(i)}_{old}$.
		
		\item Update the pairing variable using extrapolation:
		\[
		\hat{A}^{(i)} = \max
		\left(0,
		A^{(i)} + \beta (A^{(i)} - A^{(i)}_{old})
		\right).
		\]
	\end{enumerate}
	
	\item If the reconstruction error $F$ has increased, reject the extrapolation and reset pairing variables $\hat A^{(i)}= A^{(i)}$ for $1 \leq i\leq N$; otherwise, update $A^{(i)}=\hat A^{(i)}$ for $1 \leq i\leq N$.
	
	\item Update the parameter $\beta$;  see \S~\ref{sec:HER:subsec:sequenceupdate} for the details.
	
	\item If convergence criterion is not met, go back to 2.
\end{enumerate}

This approach has been scarcely studied~\cite{Bro1998Multi,mitchell2018nesterov,ang2019accelerating}, while extrapolation is a rather well understood method to accelerate both convex and non-convex single-block descent algorithms; see for instance~\cite{Ghadimi2016Accelerated,Su2016Differential}.
The main novelty of this paper is to tackle a non-convex optimization problem using BCD with extrapolation between the block update, as opposed to inside each block update such as in~\cite{liavas2017nesterov} or after each outer loop as in~\cite{mitchell2018nesterov}.  
This in-between extrapolation comes at almost no additional computational cost.

Extrapolated BCD algorithms are shown to be considerably faster
than their standard counterparts in various difficult cases.
These algorithms were observed to be slower than existing BCD algorithms only for extremely sparse tensors.
Therefore, a contribution of this work is to experimentally show that using in-between block extrapolation allows to accelerate any BCD algorithm for dense NTF.
This opens interesting questions for other optimization problems usually solved by BCD, for which such an extrapolation scheme may be applicable.

\paragraph{Context}
Let us provide a brief historical note about tensor decomposition.
The origin of tensor decomposition can be traced back to the work of Hitchcock \cite{Hitchcock1927,Hitchcock1928}, whereas the idea of using multiway analysis is credited to the work of \cite{Cattell1944,Cattell1952}.
Since then, especially after the work of Tucker in the field of psychometrics \cite{Tucker1964,Tucker1966}, tensor decomposition has spread and become more and more popular in other fields such as chemometrics \cite{SBG2004}, signal processing \cite{LM1996,CML2015}, data mining \cite{BBB2008,Morup2011}, and many more.
We refer the readers to \cite{KB2009,KHP2014,K2008,Anandkumar2014tensor,sidiropoulos2017tensor} and references therein for comprehensive reviews of the applications of tensor decomposition.
It is important to note that NTF is just one of many tensor decomposition models.
Some other types of tensor decomposition or format include PARAFAC (that is, unconstrained approximate tensor factorization),  Tucker/HOSVD
\cite{Tucker1966,DeLathauwer2000multilinear}, and Tensor Train~\cite{Oseledets2011Tensor}, to name a few.
We focus on NTF in this paper.

%NMF
Nonnegative matrix factorization (NMF), a key problem in machine learning and data analysis, is a special case of NTF when $N=2$.
First introduced in \cite{PT1994}, it started becoming widely used after the seminal work of \cite{lee1999learning}, and NMF has since then been deeply studied and well documented with variety of applications such as document classification, image processing, audio source separation and hyperspectral unmixing;
see~\cite{cichocki2009nonnegative, Gillis2014, fu2018nonnegative} and the references therein for more details.
%why tensor?
On the other hand, there are classes of data for which being represented by tensors is more natural.
For example, a third-order tensor is preferably used to connect excitation-emission spectroscopy matrices in chemometrics \cite{SBG2004}, and RGB color images or 3D light field displays are generated as tensors \cite{Hirsch2012}; see \cite{cichocki2009nonnegative} for more examples.
NTF was first introduced in \cite{Carroll1989} for fitting the latent class model in statistics.
It has also been applied in model selection problem, sparse image coding in computer version \cite{SH2005}, sound source separation \cite{FCC2005}, image decomposition \cite{WW2001}, text mining \cite{chi2012tensors}, among others; see \cite{BD1997,cichocki2009nonnegative,FH2008,HPS2005,LC2009} and reference therein for more applications of NTF.

\paragraph{Notation}
Below we recall some important notations in tensor algebra.
First of all, the Kronecker product~\cite{Brewer1978Kronecker} of two matrices $A\in\mathds{R}^{I_1\times J_1}$ and $B\in\mathds{R}^{I_2\times J_2}$ is defined as follows:
\begin{equation}\label{eq:Kron}
	A \kron B =
	\left[
	\begin{array}{c|c|c}
		[A]_{11}B    & \cdots & [A]_{1J_1}B \\ \hline
		\vdots       & \ddots & \vdots      \\ \hline
		[A]_{I_1 1}B & \cdots & [A]_{I_1 J_1}B
	\end{array}
	\right].
\end{equation}
Moreover, the Kronecker product of several matrices can be deduced from the above definition by associativity.
The Khatri-Rao product $A\odot B$ is the columns-wise Kronecker product.
Setting $A = [a_1,\ldots,a_{J_1}]$ and $B = [b_1, \ldots, b_{J_1}]$,
\begin{equation}\label{eq:kr}
	A\odot B = \left[a_1\kron b_1,\ldots, a_{J_1}\kron b_{J_1} \right].
\end{equation}
The Hadamard product (element-wise product) is denoted $A\circledast B$.

\noindent \underline{Compact decomposition notations:}
There exist several complementary notations to parameterize a low-rank tensor.
In particular, grouping components $a^{(i)}_p$ as columns of factor matrices $A^{(i)} = [a^{(i)}_1, \dots, a^{(i)}_r]$, the following notations are equivalent:
\begin{eqnarray}
	\tX
	&=& ~~ \sum_{p=1}^{r} \bigotimes_{i=1}^{N} a^{(i)}_p \\
	\label{eq:kruskal}
	&=& ~~ \left[ A^{(1)},\ldots ,A^{(N)}\right] \\
	\label{eq:nmodeproduct}
	&=& ~~ \tI_r \times_1 A^{(1)}\times_2 \ldots \times_{N} A^{(N)}\\
	\label{eq:atensp}
	&\underset{\text{def}}{:=}& ~~ \left({\bigotimesa_{i=1~}^{N~}} A^{(i)} \right)\tI_r .
\end{eqnarray}
where $\bigotimesa$ is a tensor product of linear maps induced by the tensor product $\otimes$ of vectors.

Equation~\eqref{eq:kruskal} is the so-called Kruskal notation, equation~\eqref{eq:nmodeproduct} makes use of the n-mode product $\times_p$~(see \cite{KB2009}), and equation~\eqref{eq:atensp} uses the fact that linear applications on tensor spaces of finite dimensions also form a tensor space with tensor product $\left(A\otimes_a B\right) (x\otimes y) := Ax \otimes By$.
Because~\eqref{eq:atensp} exhibits this tensor product structure, we will make use of this compact formulation rather than the others. \vspace{1em}

\noindent \underline{Tensor unfoldings and useful formula}:
To derive partial derivatives of the NTF cost with respect to factors matrices, it is convenient to switch from a tensor formulation to a matrix description of the problem.
More precisely, the following relationships hold:
\begin{equation}\label{eq:unfold}
	\tX = \left( \bigotimesa_{i=1~}^{N~} A^{(i)} \right) \tI_r
	\; \qquad \equiv \qquad \;
	\forall i\in[1,N], \;
	X_{[i]} = A^{(i)}\left(\bigodot_{\substack{l\neq i \\ l=N}}^{1} A^{(l)}\right)^T,
\end{equation}
where unfoldings $X_{[i]}$ of a rank-one tensor $\tX$ are defined as follows:
\begin{equation}\label{eq:lrunfold}
	X_{[i]} :=
	a^{(i)} \otimes \left( 
	%\boxtimes_{\underset{l \neq i}{l=N}}^1 a^{(l)} 
	\mathop{\scalebox{2}{\raisebox{-0.2ex}{$\boxtimes$}}}_{\substack{l \neq i \\ l=N}}^1 a^{(l)} 
	\right) \;
	\in \mathbb{R}^\mathcal{D}, ~ \mathcal{D} = I_i \times \underset{l\neq i}{\prod}I_l.
\end{equation}
Unfoldings of a general tensor are obtained by linearity of the unfolding maps.
Note that several non-equivalent definitions are used in the tensor signal processing community; see~\cite{KB2009} and \cite{Cohen2015About}.

\section{The state-of-the-art algorithms for solving NTF}
\label{sec:survey}
Below, we provide an overview of various techniques to solve NTF, which can be reformulated as follows
\[
\min_{A^{(i)} \geq 0, 1 \leq i \leq N} \; F(A^{(1)},\ldots,A^{(N)}),
\]
where
\begin{equation}
	F(A^{(1)},\ldots,A^{(N)}) = \frac12
	\norm{
		\mathcal{T}  - \left( \bigotimesa_{i=1~}^{N~} A^{(i)} \right) \mathcal{I}_r
	}_F^2.
	\label{eqn:objFun}
\end{equation}
As a foreword, let us mention that there exist a wide range of algorithmic solutions for NTF (as for most of the tensor decomposition problems), that show different performances depending on the task at hand.

\paragraph{Algorithms for exact PARAFAC}
First of all, although the focus of this paper is approximate  decompositions, it is interesting to point out that several algebraic techniques based on eigendecompositions have been proposed to deal with exact unconstrained tensor factorizations~\cite{Lathauwer2004computation,Luciani2011semi,Domanov2014canonical} whenever such a factorization exists. However, by design, these algorithms are typically not robust to noise, or may even be numerically unstable; see~\cite{Beltran2019pencil}. These techniques are sometimes used for initialization. Moreover, to the best of our knowledge, there does not exist an algorithm that computes NTF exactly for any rank using algebraic techniques. 
In fact, even in the matrix case, that is, for nonnegative matrix factorization (NMF), such techniques do not exist. 
Therefore, in this paper, we do not discuss exact NTF algorithms; this is a direction for further research.  

\paragraph{Algorithms for approximate unconstrained tensor decomposition}
Because the Tensor Factorization (TF) model (that is, without nonnegativity constraints) has some interesting identifiability properties, it may occur that for well-conditioned tensors~\cite{Vannieuwenhoven2017condition}, approximate NTF can be computed with high precision by using a TF solver.
Solving the TF problem is however harder in theory (the tensor low-rank unconstrained approximation problem is ill-posed, see~\cite{Silva2008Tensor}) and not really easier in practice than solving NTF.  Actually, many algorithms that solve NTF are inspired from TF solvers and have similar complexity.
Therefore we do not discuss TF solvers in what follows, and assume the reader is interested in solving NTF with specific algorithms that make use of the properties of the NTF problem.

\paragraph{All-at-once optimization}
A first class of widely used methods to solve NTF are all-at-once gradient-based methods.
Indeed, it is quite straightforward to compute the gradients of $F$ with respect to each matrix $A^{(i)}$.
Let us denote
\begin{eqnarray}
	B^{(i)}=A^{(N)} \odot  \cdots  \odot A^{(i+1)} \odot A^{(i-1)}\odot \cdots \odot A^{(1)}.
	\label{def:B}
\end{eqnarray}
Then the gradient of $F$ with respect to $A^{(i)}$ is
\begin{equation}
	\label{def:grad}
	\nabla_{A^{(i)}} F=\Big( A^{(i)} \big( B^{(i)}\big)^T  - \mathcal T_{[i]}\Big) B^{(i)}.
\end{equation}

Therefore,  there is no obstacle to using any constrained gradient-based algorithm to (try to) find a stationary point of the non-convex NTF problem.
To the best of our knowledge, the oldest all-at-once algorithm for NTF is a Gauss-Newton approach~\cite{Paatero1997weighted}, but many approaches have been tested, including:
\begin{itemize}
	\item  Second-order optimization: using the fact that surrogates of the Hessian of $F$ are heavily structured, second-order information can be used to solve NTF at a reasonable cost~\cite{Vervliet2018compressed}. Limited-memory BFGS has also been employed when scalability is required~\cite{Acar2011ScalableA}.
	To enforce the nonnegativity constraints, one can for instance square the  variables, or use a variational approach (such a log-barrier).
	\item Primal-Dual optimization: the alternating direction method of  multipliers has been tested for NTF, with however less promising results than its block-coordinate counterpart discussed below, see \cite{Huang2016flexible}.
	\item Conjugate gradient:  it has been reported that conjugate gradient can be used to solve NTF by squaring the variables; see~\cite{Royer2012nonnegative}.
\end{itemize}

\paragraph{Block Coordinate Descent (BCD) Methods}
Other than the above mentioned algorithms, BCD has become a standard and efficient scheme for solving NTF, mainly because
(1) it essentially has cheap computation cost in each block update (BCD fixes all blocks except for one),
(2) BCD can make use of recent developments in convex constrained optimization to efficiently solve NTF with respect to each block,
and (3) under some suitable assumptions, many first-order BCDs and their accelerated versions have convergence guarantee in the context of general block-separable \emph{non-convex composite optimization} problem that subsumes NTF as a special case, see for example \cite{xu2013block} and \cite{Hien2019} and the references therein.
Below, we review several block coordinate methods for solving NTF; we list these algorithms in Table \ref{table:BCDalgorithm}.

\begin{table}
	\centering
	\stepcounter{table}
	\def\arraystretch{0.9}
	{
		%\begin{tabular}{lll}
		%Algorithms  & Section                         & Reference  \\ \hline
		%AO-AS       & \S\,\text{\ref{subsec:AO-AS}}   &\text{\cite{KHP2014}}     \\
		%AO-ADMM     & \S\,\text{\ref{subsec:ADMM}}    &\text{\cite{Huang2016flexible}}   \\
		%AO-Nesterov & \S\,\text{\ref{subsec:Nesterov}}&\text{\cite{Zhang2016fast}}        \\
		%A-HALS      & \S\,\text{\ref{subsec:AHALS}}   &\text{\cite{Gillis2014}}     \\
		%APG        & \S\,\text{\ref{subsec:APG}}     &\text{\cite{xu2013block}}                     \\
		%iBPG       & \S\,\text{\ref{subsec:iBPG}}    &\text{\cite{Hien2019}}  %\\
		%\end{tabular}
		\begin{tabular}{c|cccccc}
			Algorithms&AO-AS\cite{KHP2014}                        &AO-ADMM\cite{Huang2016flexible}                        &AO-Nesterov\cite{Zhang2016fast}                     &A-HALS\cite{Gillis2014}                       &APG\cite{xu2013block}                   & iBPG\cite{Hien2019} \\
			Section   &\S\,\text{\ref{subsec:AO-AS}}&\S\,\text{\ref{subsec:ADMM}}   &\S\,\text{\ref{subsec:Nesterov}}&\S\,\text{\ref{subsec:AHALS}}&\S\,\text{\ref{subsec:APG}}&\S\,\text{\ref{subsec:iBPG}}
		\end{tabular}
	}
	\addtocounter{table}{-1}%
	\caption{Several block coordinate methods for solving NTF}
	\label{table:BCDalgorithm}
\end{table}

\subsection{Alternating optimization (AO) framework}\label{sec:AOframework}
When solving NTF using BCDs, the blocks of variables that are alternatively updated must be chosen.
It turns out that $F$ is a quadratic function with respect to each matrix $A^{(i)} $ and therefore the optimization problem
\begin{equation}
	\label{eq:AO}
	\underset{A^{(i)}\geq 0}{\min }~ F(A^{(1)},\ldots,A^{(N)})
\end{equation}
is a linearly constrained quadratic programming problem referred to as Nonnegative Least Squares (NNLS).
In particular, it is strictly convex if and only if $B^{(p)}$ is full column-rank.
Therefore, it is quite natural to consider $A^{(i)}$ as the blocks in a BCD.
The AO framework, which is a standard procedure to solve NTF, alternatively (exactly/inexactly) solves \eqref{eq:AO} for each block.
We describe the AO framework in Algorithm~\ref{alg:ANLS}. Note that the objective function of AO methods decreases after each block update.
Depending on how the matrix-form NNLS problem \eqref{eq:NNLS} is solved, various implementations of AO algorithms can be obtained.
Some of them are very efficient for solving NTF, they are surveyed below.

%Alternating Nonnegative Least Squares (ANLS) and is a standard procedure to solve NTF. ANLS is .
\begin{minipage}{0.92\linewidth}
%\algsetup{indent=2em}
\begin{algorithm}[H]
\caption{Alternating optimization framework \label{alg:ANLS} }
\begin{algorithmic}[1]
\STATE Input: a nonnegative $N$-way tensor
\STATE Output: nonnegative factors $A^{(1)}, A^{(2)},\ldots, A^{(N)}$.
\STATE Initialization: $\big(A^{(1)}_{0}, \ldots, A^{(N)}_{0}\big )$. Set $k=1$.
\REPEAT
\FOR{$i=1,\ldots,N$}
\STATE Update $A_k^{(i)}$ as an exact/inexact solution of \vspace*{-3mm}
\begin{equation}
\label{eq:NNLS}
\min_{A^{(i)}\geq 0} F\Big(A_k^{(1)},\ldots,A_k^{(i-1)},A^{(i)},A_{k-1}^{(i+1)},\ldots,A_{k-1}^{(N)} \Big). \vspace*{-3.5mm}
\end{equation}
($A_{k-1}^{(i)}$ can be used as the initial point for the algorithm used to solve \eqref{eq:NNLS}.)
\ENDFOR
\STATE Set $k=k+1$.
\UNTIL{some criteria is satisfied}
\end{algorithmic}
\end{algorithm}
\end{minipage}

\subsubsection{AO-AS --  solving NNLS with Active Set}\label{subsec:AO-AS}
When $A_k^{(i)}$ is updated by an exact solution of the NNLS problem~\eqref{eq:NNLS}, we obtain an alternating nonnegative least squares algorithm, usually referred to as ANLS in the literature.
To solve exactly the NNLS subproblem~\eqref{eq:NNLS}, active set (AS) methods are usually rather effective and popular; see~\cite{KHP2014}.
We will refer to AO-AS as the ANLS algorithm where the NNLS subproblems are solved with AS.

\subsubsection{AO-ADMM -- solving NNLS with ADMM}\label{subsec:ADMM}
Designed to tackle a wide range of constrained tensor decomposition problems and various loss functions, AO-ADMM \cite{Huang2016flexible} applied to NTF boils down to using several steps of a primal-dual algorithm, the Alternating Direction Method of Multipliers (ADMM), to solve the cascaded nonnegative least squares problems.
Therefore, AO-ADMM for NTF problem~\eqref{eq:NTF} is a variant of the AO framework that solves \eqref{eq:NNLS} inexactly.

\subsubsection{AO-Nesterov}\label{subsec:Nesterov}
When Nesterov's accelerated gradient method is applied to solve the NNLS problem in \eqref{eq:NNLS}, we obtain AO-Nesterov; see~\cite{guan2012nenmf, Zhang2016fast}.

\subsubsection{A-HALS}\label{subsec:AHALS}
The hierarchical alternating least square (HALS) algorithm was introduced for solving the nonnegative matrix factorization (NMF) problem  $\underset{W\geq 0,H\geq 0}{\min} \norm{M-WH^T}_F^2$ (that is, NTF when $N=2$), and has been widely used for solving NMF as it performs extremely well in practice; see for example~\cite{cichocki2009nonnegative,Gillis2014}.

HALS cyclically updates each column of the factor matrix $A^{(i)}$ by solving an NNLS problem with respect to that column while fixing the others.
The optimal solution of this NNLS subproblem can be written in closed form.
A-HALS, which is short for accelerated HALS, was proposed in~\cite{Gillis2012} to accelerate HALS.
A-HALS repeats updating each factor matrix several times before updating the other ones.
Hence A-HALS can be considered as a variant of the AO framework where each NNLS is inexactly solved itself by a BCD with closed-form updates.
Let us briefly describe A-HALS for solving NTF.
The NNLS problem \eqref{eq:NNLS} of A-HALS  is inexactly solved by repeating cyclically updating the columns of $A_{k-1}^{(i)}$.
In particular, let $M=X_{(i)}$, $W=A_{k-1}^{(i)}$, and $H=A_{k-1}^{(N)} \odot \ldots A_{k-1}^{(i+1)} \odot A_k^{(i-1)} \ldots \odot A_k^{(1)}.$
The $j$-th column of $A_{k-1}^{(i)}$ is updated by
\begin{equation*}
	W_{:,j}=\frac{   \max\left( ~0~,~ MH_{j,:}^T-\sum_{l\ne j} W_{:,l}H_{l,:}H_{j,:}^T \right) }{\norm{H_{j,:}}^2}.
\end{equation*}
It is worth noting that A-HALS for NTF has subsequential convergence guarantee (that is, every limit point is a stationary point of the objective function), see \cite[Section 7]{Razaviyayn2013}.

\subsection{Block proximal gradient type methods}\label{subsec:proximal}
The NNLS problem \eqref{eq:NNLS} does not have a closed-form solution.
From Equation \eqref{def:grad}, we see that $F$, when restricted to $A^{(i)}$, is a $L^{(i)}$-smooth function, that is, the gradient $\nabla_{A^{(i)}} F$ is Lipschitz continuous with the constant $L^{(i)}=\norm{\big(B^{(i)}\big)^T B^{(i)}}$, where $B^{(i)}$ is defined in (\ref{def:B}).
This property can be employed to replace the objective function in the NNLS problem \eqref{eq:NNLS} by its quadratic majorization function, that leads to a new minimization problem which has a closed-form solution.
This minimization-majorization approach, in the literature of block-separable composite optimization problem with the block-wise $L$-smooth property, is known as proximal gradient block coordinate descent method (see e.g., \cite{Hien2019}). 
Considering the NTF problem,  the closed-form solution of minimizing the majorization function is a projected gradient step.
Applying Nesterov-type acceleration for the proximal gradient step improves the convergence of the BCD algorithm.
Below we review the two recent accelerated proximal gradient BCD methods that were proposed for solving the general composite optimization problem.

\subsubsection{APG --  An Alternating Proximal Gradient method for solving NTF}\label{subsec:APG}
APG was proposed by Xu and Yin~\cite{xu2013block}; see the Appendix of \cite{2020arXiv200104321A} and \cite[Section 3.2]{xu2013block} for the algorithm pseudocode.
APG \emph{cyclically} update each block (a.k.a each factor matrix) $A^{(i)}$ by calculating an extrapolation point $\hat  A_{k-1}^{(i)}= A_{k-1}^{(i)} + w_{k-1}^{(i)} \Big( A_{k-1}^{(i)}-A_{k-2}^{(i)} \Big)$ (here $ w_{k-1}^{(i)}$ is some extrapolation parameter)  and embedding this point in a projected gradient step
\[
A_{k}^{(i)}=\max \left(0,\hat  A_{k-1}^{(i)}- \frac{1}{L_{k-1}^{(i)}}  \Big (\hat A_{k-1}^{(i)} \big(B_{k-1}^{(i)}\big)^T- \mathcal T_{[i]}\Big) B_{k-1}^{(i)}\right).
\]
After all blocks are updated, APG needs a restarting step, that is, if the objective function has increased  then the projected gradient step would be re-done by using the previous values of all blocks instead of using the extrapolation points.

\subsubsection{iBPG -- An inertial Block Proximal Gradient Method}\label{subsec:iBPG}
Recently proposed in \cite{Hien2019}, iBPG computes two different extrapolation points $\hat A_{k-1}^{(i,1)}$ and $\hat A_{k-1}^{(i,2)}$: one is for evaluating the gradient and the other one for adding inertial force. iBPG updates one matrix factor using a projected gradient step
\[
A_{k}^{(i)}=\max \left(0,\hat  A_{k-1}^{(i,2)}- \frac{1}{L_{k-1}^{(i)}}  \Big (\hat A_{k-1}^{(i,1)} \big(B_{k-1}^{(i)}\big)^T- \mathcal T_{[i]}\Big) B_{k-1}^{(i)}\right),
\]
see the Appendix of \cite{2020arXiv200104321A} for the algorithm pseudocode.
Furthermore, similarly to A-HALS, iBPG allows updating each matrix factor some times before updating another one -- this feature would help save some computational costs since some common expressions can be re-used when repeating updating the same block. iBPG does not require a restarting step which make it suitable for solving large-scale NTF problems where evaluating the objective functions is costly.

\section{Making BCD significantly faster with HER}\label{sec:exBCD}
With modern machine learning applications of NTF in mind, for which input tensor sizes can be extremely large and NTF should be provided as a low-level routine, there would be a definite economical and scientific gain to speeding up NTF algorithms.
Radically different approaches exist in the literature to speed up existing algorithms for solving NTF, such as parallel computing~\cite{ravindran2014memory,Ballard2018parallel}, compression and sketching~\cite{Vervliet2019exploiting,battaglino2018practical}.
The combinations and relationships between these methods is poorly understood.
In this paper, we focus on the acceleration of BCD using extrapolation.

Extrapolation inside BCD algorithms such as the workhorse ALS algorithm has been studied in the tensor decomposition literature as an empirical trick to speed up computations and avoid swamps, see Section~\ref{subsec:Bro}. However, in a recent work on rank-one approximations of rank-two tensors, Gong, Mohlenkamp and Young~\cite{gong2018optimization} used gradient flow to study transverse stability, and  provided a deeper analysis of the optimization landscape of tensor low-rank approximation which further supports the use of extrapolation.

As reviewed in \S \ref{subsec:proximal}, we have seen that APG and iBPG accelerate block proximal gradient methods by using extrapolation points in the projected gradient step to update each factor matrix. In another line of works, AO (Algorithm~\ref{alg:ANLS}) was   accelerated by using extrapolation \underline{between} each block update (rather than inside the block update as in APG and iBPG); in other words, each factor matrix is updated by the extrapolation between previous updated factors. In the literature of tensor decomposition, the second type of extrapolation has been used to accelerate alternating least squares algorithms for solving CPD. Those works will be reviewed in \S~\ref{subsec:Bro}. In the following, we introduce HER - a novel extrapolation scheme that can be categorized into the class of accelerated AO algorithms using extrapolation between block update.

\subsection{Heuristic Extrapolation with Restarts (HER)}\label{subsec:HER}
HER was first proposed for solving NMF in \cite{ang2019accelerating}, and found to be extremely effective on NTF in a {preliminary work} \cite{ang2019accelerating2}. The sketch of HER was given in the introduction and its pseudo-code is given in Algorithm \ref{alg:HER}. In the following, we elaborate on HER with more details.

%\algsetup{indent=2em}
\begin{algorithm}[t]
	\caption{HER \label{alg:HER}}
	\begin{algorithmic}[1]
	\STATE Input: a nonnegative $N$-way tensor%, and a switch parameter  $\mathbf{p} \in\{0,1\}$
	\STATE Output: nonnegative factors $A^{(1)}, A^{(2)},\ldots, A^{(N)}$.
	\STATE Initialization:
	Choose $\beta_0 \in (0,1)$, $\eta \geq \bar{\gamma} \geq \gamma \geq 1$ and 2 sets of initial factor matrices
	$\big(A^{(1)}_{0}, \ldots, A^{(N)}_{0}\big )$ and $\big(\hat{A}^{(1)}_{0}, \ldots, \hat{A}^{(N)}_{0}\big)$.
	Set $\bar{\beta}_0 =1$ and $k=1$.
	\REPEAT
	\FOR{$i=1,\ldots,N$}
	\STATE \textbf{Update step}
	Let $A_k^{(i)}$ be an exact/inexact solution of
	\begin{minipage}{0.96\linewidth}
	\begin{equation}
	\min_{{A}^{(i)}\geq 0} F\left(\hat{A}_k^{(1)},\ldots,\hat{A}_k^{(i-1)},A^{(i)},\hat{A}_{k-1}^{(i+1)},\ldots,\hat{A}_{k-1}^{(N)} \right).                   \label{NNLS_her_exact}
	\end{equation}
	\end{minipage}

	\STATE \textbf{Extrapolation step}
	\begin{minipage}{0.96\linewidth}
	\begin{equation}	\label{extproj_her}
     \hat{A}_k^{(i)} = \max \left ( 0, A^{(i)}_{k} + \beta_{k-1} ( A^{(i)}_{k} - A_{k-1}^{(i)}) \right ).
	\end{equation}
	\end{minipage}

	\ENDFOR
	\STATE Compute $\hat{F}_k := F\left(\hat{A}_k^{(1)},\hat{A}_k^{(2)},\ldots,\hat{A}_k^{(N-1)},A_k^{(N)} \right). $
	 \IF{$\hat{F}_k > \hat{F}_{k-1}$}
	\STATE Set $\hat{A}_{k}^{(i)} = A^{(i)}_k$, $i = 1, ..., N$  \hfill \% \textit{abandon the sequence $\hat{A}_k^{(i)}$}
	\STATE Set $\bar{\beta}_k = \beta_{k-1}$, ~ $\beta_k = \beta_{k-1}/\eta $. \hfill	\% \textit{Update $\bar{\beta}$, decrease $\beta$}
	 \ELSE
	\STATE Set $A^{(i)}_k  = \hat{A}_{k}^{(i)}$, $i = 1, ..., N$. \hfill \% \textit{keep the sequence $\hat{A}_k^{(i)}$}
	 \STATE Set $\bar{\beta}_k = \min\{ 1, \bar{\beta}_{k-1} \bar{\gamma}  \}$, ~ $\beta_k = \min\{ \bar{\beta}_{k-1}, \beta_{k-1} \gamma \}$.
	 \hfill \% \textit{Increase $\bar{\beta}$ and $\beta$}
	\ENDIF
	\STATE Set $k=k+1$.
	\UNTIL{some criteria is satisfied}
\end{algorithmic}
\end{algorithm}
%\end{minipage}

\subsubsection{Update step -- line 6}
It is clear that Algorithm \ref{alg:HER} has the form of an alternating optimization framework in which the key optimization sub-problem \eqref{NNLS_her_exact} is a NNLS problem.
As reviewed in \S~\ref{sec:survey}, some efficient algorithms for the NNLS problem \eqref{NNLS_her_exact} include AS, ADMM, Nesterov's accelerated gradient, or A-HALS.
The main difference between AO and HER is that HER does not use the latest values of the other blocks $A^{(j)}$ ($j\ne i$) but employs the latest values of their extrapolation $\hat A^{(j)}$ ($j\ne i$).
For convenience, we refer to $\left\{ \hat A^{(i)}_k,i=1,\ldots,N\right\}_{k\geq 0}$ as the extrapolation sequence.

\subsubsection{Extrapolation step -- line 7}
After the update of $A_{k}^{(i)}$, the same block of the extrapolation sequence $\hat{A}^{(i)}_k$ is updated by extrapolating $A^{(i)}_k$ along the direction $A_k^{(i)} - A^{(i)}_{k-1}$, see (\ref{extproj_her}).
Note that $\hat{A}^{(i)}_k$ produced by (\ref{extproj_her}) is always feasible.
It is possible to remove the projection in (\ref{extproj_her}), but we do not consider such approach in this work.
Note that, regarding feasibility, $A^{(i)}_k$ produced by line 6 of Algorithm \ref{alg:HER} is always feasible regardless of the feasibility of $\hat{A}^{(i)}_k$.

\subsubsection{The restart mechanism -- lines 9-16}\label{sec:HER:subsec:restartmechanism}
After the update-extrapolate process on all the blocks, a restart procedure is carried out to decide whether or not we replace $A^{(i)}$ ($1 \leq i \leq N$) with the extrapolation sequence.  The command in line 14 of Algorithm \ref{alg:HER} has the same spirit with the update in \eqref{eqn:update2_y}  \color{black} where the factor matrices are updated by \emph{the extrapolation between block update}.

It may raise a question why $F\big({A}_{k-1}^{(1)},\ldots,A_{k-1}^{(N)} \big)$ does not appear in the restarting condition -- line 10.
The answer is due to the practicality of the algorithm.
As stated in \cite{ang2019accelerating}, using $F$ as the restart criterion is computationally much more expensive than using the approximate $\hat{F}$.
When computing $\hat{F}$, no explicit computation is required; instead, one may reuse already computed components from the updates of $A^{(N)}$ and $\hat{A}^{(N)}$.
This creates an important reduction of computational complexity.
For example, let us consider an order-$N$ NTF problem with factor matrices $\big \{ A^{(i)} \big \}_{i=1,2,\dots,N}$ with size $I_1 \times r$, $I_2 \times r$, $\dots$ up to $I_N \times r$.
Reusing already computed components (such as gradient) in the update of the last block $(A^{(N)}$, $\hat{A}^{(N)})$, we can compute $\hat{F}(\hat{A}^{{1}},\dots, \hat{A}^{{N-1}},A^{{N}})$ under $I_N r^{N-1}$ flops.
However, if we compute $F(A^{{1}},\dots,A^{{N}})$, it takes $\prod_{i=1}^{N} I_i$ flops.
If $r^N \ll \prod_i  I_i$, then such reduction in complexity from $\prod_{i=1}^{N} I_i$ to $I_N r^{N-1}$ is significant even when $N$ is low. Furthermore, we can even rotate the tensor such that $I_N$ is the mode with the smallest size among all the modes.
In fact, computing the cost function naively can be as costly as one block update, and thus using $\hat{F}$ instead of $F$ as the restart criterion is important, since restart using $F$ requires computing the cost function at each iteration, while restart using $\hat{F}$ is much cheaper.

Moreover, note that if the iterates sequence is converging, then the extrapolated sequence also converges to the same limit point. Therefore, since $F$ is a continuous map, if convergence of the iterates is observed then the surrogate cost $\hat{F}$ will asymptotically converge to the same final value as $F$. Although we did not characterize how fast this convergence happens, this justifies to use $\hat{F}$ as a surrogate at least near a stationary point.

\subsubsection{The extrapolation parameters in lines 9-16}\label{sec:HER:subsec:sequenceupdate}
The extrapolation weight $\beta_k$ is computed within the restart mechanism of lines 9-16 of Algorithm~\ref{alg:HER}, and it is updated using four parameters; see Table~\ref{table:HER_params}.
\begin{table}[ht!]
	\centering
	\stepcounter{table}% for \thetable
	\def\arraystretch{1.01}
	\begin{tabular}{cllll}
		Symbol          & Name                         & Setting                                     & Range                & Requires tuning?   \\ \hline
		$\beta_k$       & Extrapolation weight         & update as (\text{\ref{HER_eq:beta_update}})        &$[0,1]$               & Yes for $\beta_0$  \\
		$\gamma$        & Growth rate of $\beta$       & constant                                    &$[\bar{\gamma},\eta]$	& Yes                \\
		$\eta$          & Decay rate of $\beta$        & constant       					         & $[\gamma, \infty)$	& Yes               \\
		$\bar{\gamma}$  & Growth rate of $\bar{\beta}$ & constant                                    & $[1, \gamma]$		& Yes \\
		$\bar{\beta}_k$ & Upper bound for $\beta$           & update  as (\text{\ref{HER_eq:beta_bar_update}})   &	$[\beta_k, 1]$      & No, $\bar{\beta}_0 = 1$
	\end{tabular}
	\addtocounter{table}{-1}%
	\caption{Parameters in the HER scheme}
	\label{table:HER_params}
\end{table}

In the initialization stage, we set the upper bound for $\beta$ as $\bar{\beta}_0 = 1$, pick $\beta_0 \in (0,1)$, and select $\eta, \gamma$ and $\bar{\gamma}$ such that $1 < \bar{\gamma} \leq \gamma \leq \eta$.
The parameter $\bar{\beta}$, which is initialized as 1, is called the upper bound parameter for $\beta$.
This parameter is used to limit the growth of $\beta$; see below for more details.
The parameter $\gamma$ is called the (multiplicative) growth rate of $\beta$: when the error decreases, $\beta$ is updated with $\gamma \beta$.
Similarly $\bar{\gamma}$ is the (multiplicative) growth rate of $\bar{\beta}$.
Finally, $\eta$ is called the decay rate of $\beta$.
This value is used to update $\beta$ with $\beta/\eta$ when the error increases.
The parameters $(\gamma, \bar{\gamma},\eta)$ are fixed constants, while $\beta$ and $\bar{\beta}$ are updated depending on the restart condition.

\paragraph{The update of $\beta$} HER updates $\beta_k$ as
\begin{eqnarray}
	\beta_{k+1} & = &
	\left\{\begin{matrix}
		\beta_k / \eta & \text{if } \hat{F}_{k+1} > \hat{F}_k \\
		\min\{ \gamma \beta_k, \bar{\beta}_k \} &  \text{if } \hat{F}_{k+1} \leq \hat{F}_k
	\end{matrix}\right. ,
	\label{HER_eq:beta_update}
\end{eqnarray}
which is explained as follows :
\begin{itemize}
	\item If restart occurs, that is, if $\hat{F}_{k+1} > \hat{F}_k$, we assume it is caused by an over-sized $\beta_k$ (recall that, for $\beta_k = 0$, decrease is guaranteed by the update in line 6) and we shrink the value of $\beta$ for the next iteration using the decay parameter $\eta$ as in (\ref{HER_eq:beta_update}).
	
	\item Otherwise, $\hat{F}_{k+1} \leq \hat{F}_k$, and we assume $\beta_k$ can safely be increased.
	We grow $\beta$ for the next iteration as $\gamma \beta$. To prevent $\beta$ grow indefinitely, we use an upper bound $\bar{\beta}$ as in (\ref{HER_eq:beta_update}).
\end{itemize}

\paragraph{The update of $\bar{\beta}$} HER updates $\bar{\beta}_k$ as follows
\begin{eqnarray}
	\bar{\beta}_{k+1}
	& = &
	\left\{\begin{matrix}
		\beta_k                                 & \text{if } \hat{F}_{k+1} > \hat{F}_k  \\
		\min\{ \bar{\gamma} \bar{\beta}_k, 1 \} & \text{if } \hat{F}_{k+1} \leq \hat{F}_k
	\end{matrix}\right. .
	\label{HER_eq:beta_bar_update}
\end{eqnarray}
The explanations are as follows :
\begin{itemize}
	\item If there is no restart, that is, if $\hat{F}_{k+1} \leq \hat{F}_k$, $\bar{\beta}$ is increased if $\bar{\beta}$ is smaller than 1.
	
	\item Otherwise $\hat{F}_{k+1} > \hat{F}_k$ and $\bar{\beta}_{k+1}$ is set to $\beta_k$ to prevent $\beta_{k+1}$ growing larger than $\beta_k$ too fast in the future. In fact, $\beta_k$ indicates a too large value for $\beta$ since the error has increased.
\end{itemize}

\begin{center}\begin{figure}[ht!] % *** Unbalanced
		\centering
		\begin{subfigure}[b]{0.322\textwidth}
			\includegraphics[width=\linewidth]{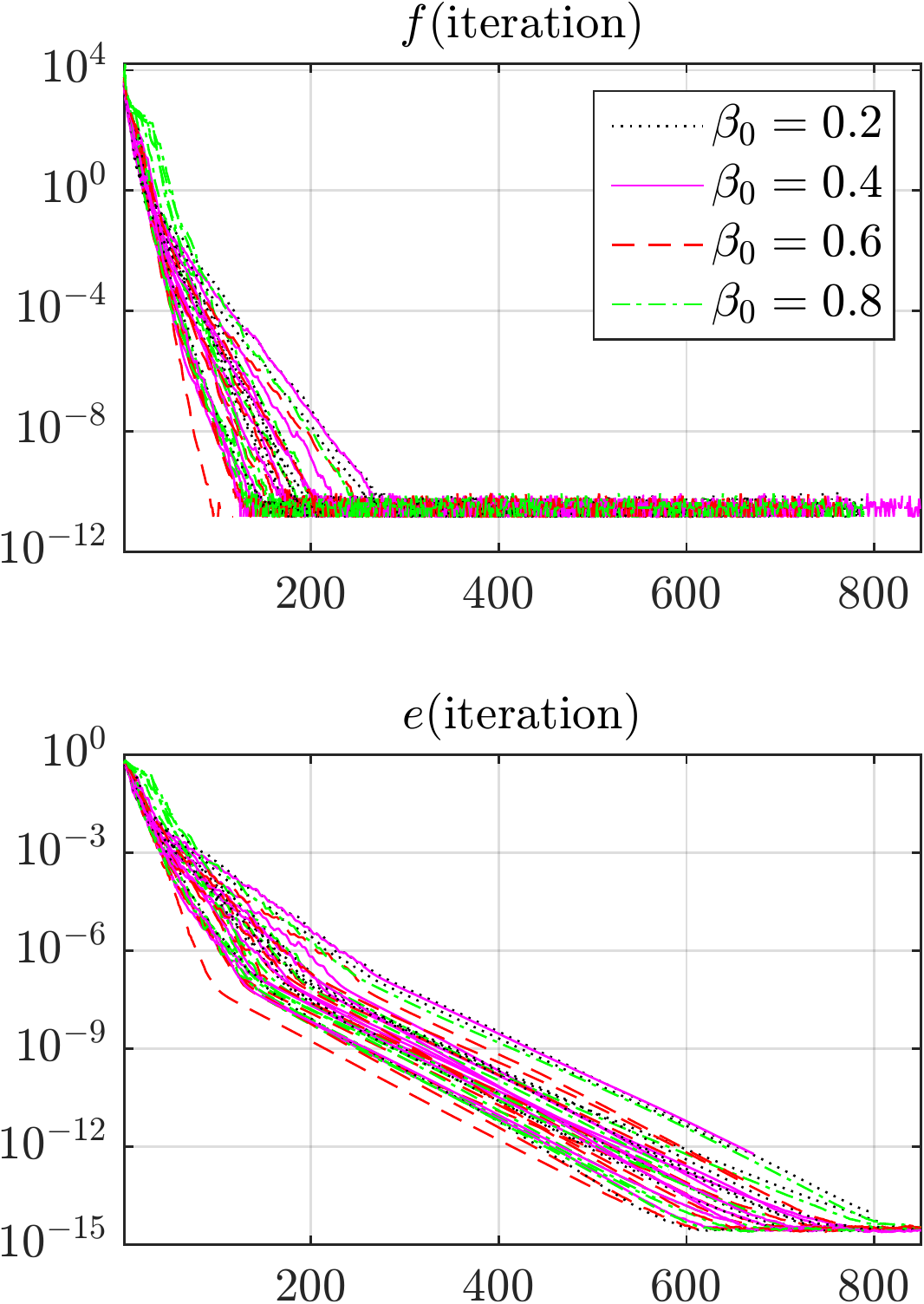}
			\caption{On fix $\gamma, \bar{\gamma}, \eta$.}
		\end{subfigure}
		\hfill
		\begin{subfigure}[b]{0.322\textwidth}
			\includegraphics[width=\linewidth]{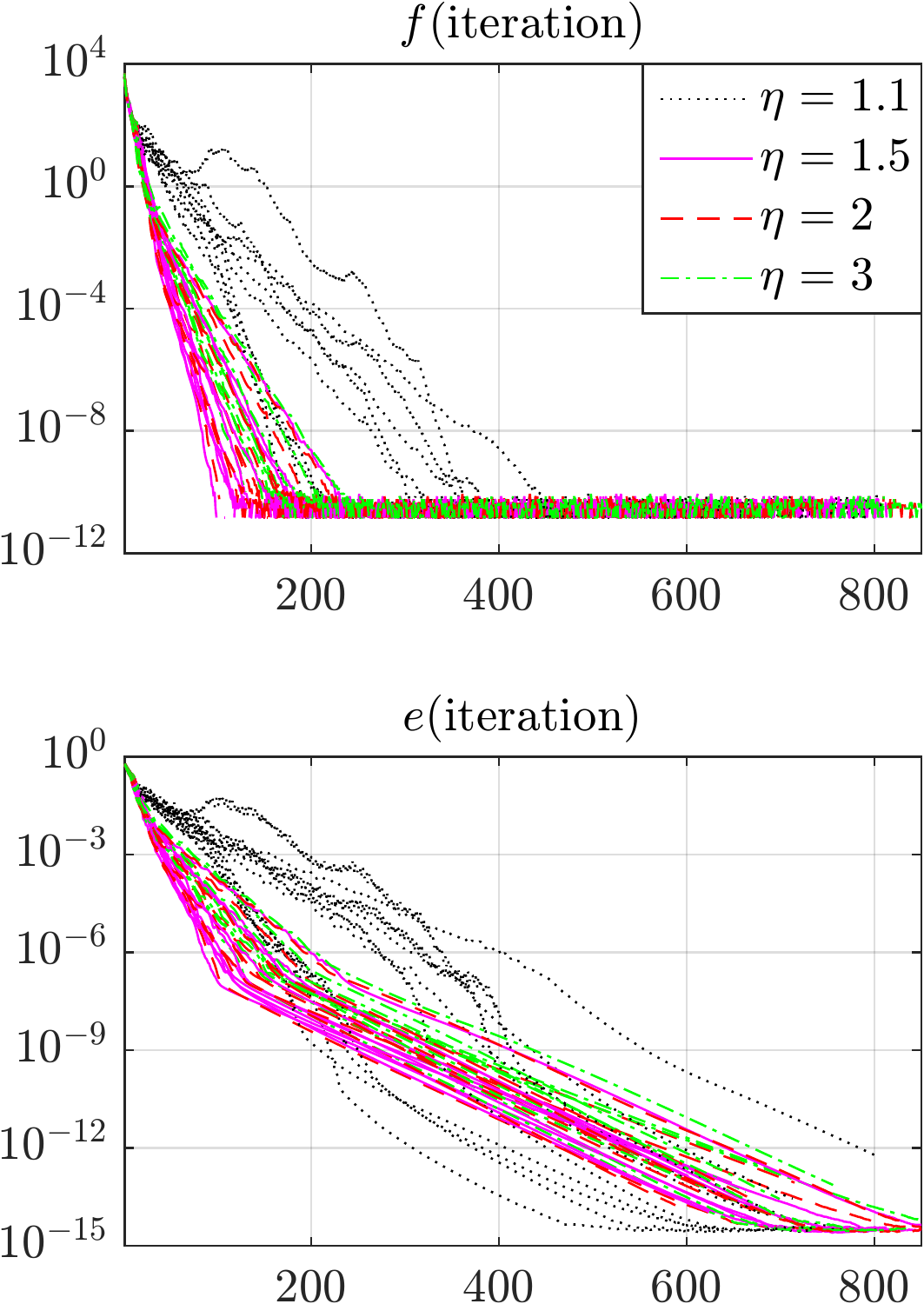}
			\caption{On fix $\gamma, \bar{\gamma}, \beta_0$.}
		\end{subfigure}
		\hfill
		\begin{subfigure}[b]{0.322\textwidth}
			\centering
			\includegraphics[width=\linewidth]{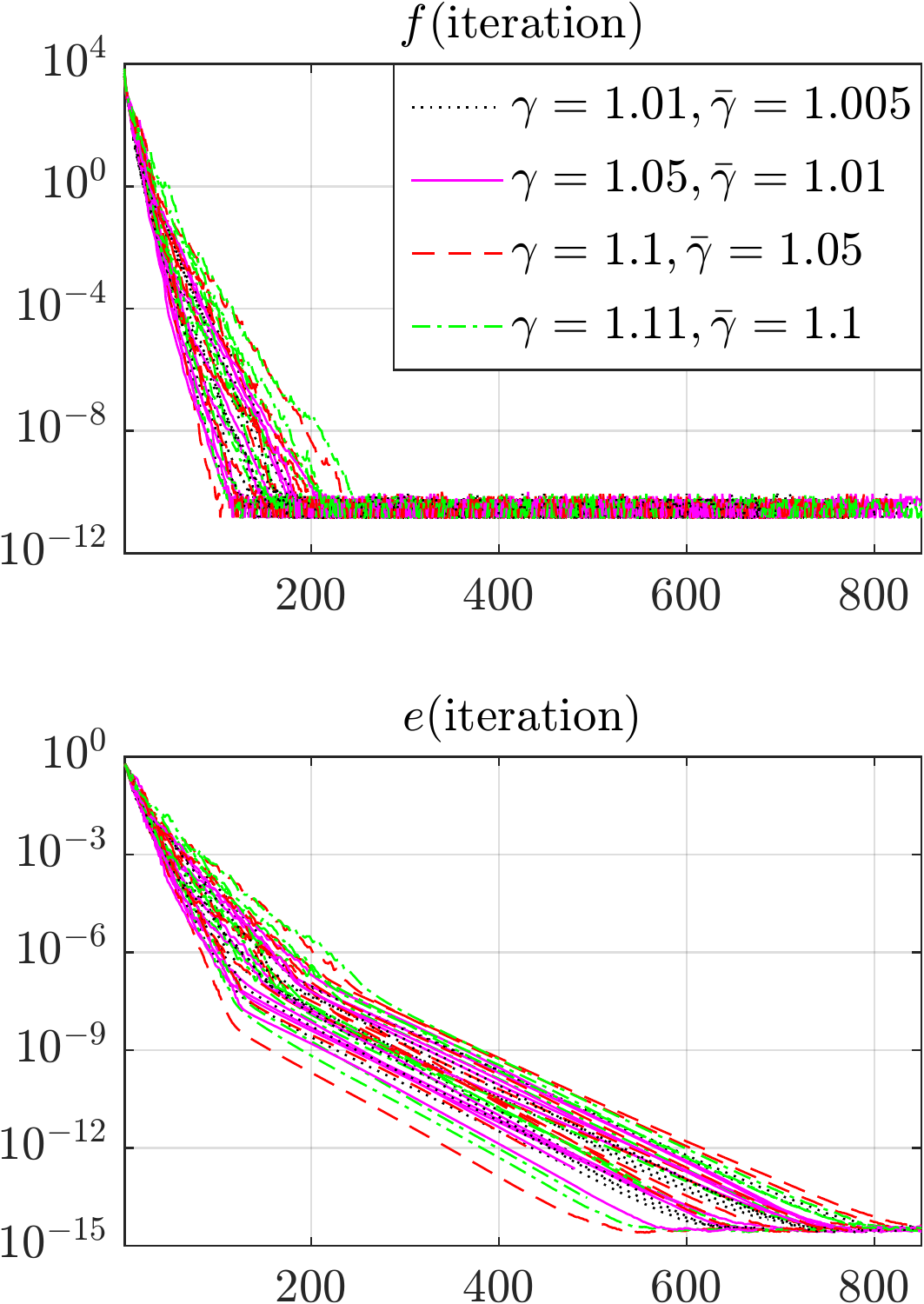}
			\caption{On fix $\beta_0, \eta$.}
		\end{subfigure}
		\caption{Comparison of HER with different parameters on the same NTF problems:
			a rank-$10$ factorization on noiseless tensors  generated by random with size $50 \times 50 \times 50$.
			For each set of parameters,
			the decomposition is repeated 10 times over 10 different data tensors and initializations;
			see \S\,\ref{sec:exp:setup} for more details.
			The top plots  representing $f$ display the error of the approximation, and the bottom plots representing $e$ display the distance to the ground truth factors; see~\eqref{err_true_def} and~\eqref{err_fact_def}.
			The default set of hyper-parameters are $[\beta_0 = 0.5, \gamma = 1.05, \bar{\gamma} = 1.01, \eta = 1.5]$.
			The results here showed that HER is not very sensitive to its parameters as all the curves are not deviating away from each other, except for the case $\eta = 1.1$, suggesting that $\eta$ should not be too small.
		}\label{fig:exp:HER_param}
\end{figure}\end{center}

Let us make a few remarks:

\begin{itemize}
	\item The relationships between the parameters in HER is as follows:
	\begin{equation}
		0 < \beta_k \leq \bar{\beta}_k \leq 1 <  \bar{\gamma} \leq \gamma \leq \eta < \infty.
		\label{HER:relation}
	\end{equation}
	By construction, $\beta_k \leq \bar{\beta}_k \leq 1$, while $\bar{\gamma} \leq \gamma$ ensures that $\bar{\beta}$ increases slower than $\beta$, while $\gamma \leq \eta$ ensures that $\beta$ is decreased faster.
	
	\item We have observed that HER is more effective if the NNLS subproblems \eqref{NNLS_her_exact} are solved with relatively high precision.
	Empirically Fig.~\ref{fig:exp:pg_apg_herpg} suggests to use HER with repeated projected gradient steps rather than just a single step.
	The suffix 50 after the algorithms' name in Fig.~\ref{fig:exp:pg_apg_herpg} means that we run 50 iterations for the algorithms to solve \eqref{NNLS_her_exact}.
	
	\item A drawback of the HER approach is the parameter tuning.
	There are 4 parameters to tune: $\beta_0, \gamma, \bar{\gamma}, \eta$.
	However HER is not too sensitive for reasonable values of the parameters; see Fig.~\ref{fig:exp:HER_param} for an illustration.
	Therefore, all the experiments in this paper are executed with no parameter tuning, even in difficult cases when data are ill-conditioned or rank is very high; namely we will use
	$\beta_0=0.5$,
	$\gamma=1.05$,
	$\bar{\gamma}=1.01$ and
	$\eta=1.5$.
	
	\item In the implementation, we initialize $\hat{A}^{(i)}_0 = A^{(i)}_0, i=1,\ldots,N$.
\end{itemize}

\subsection{Discussion on convergence} \label{subsec:convHER}
Unfortunately, the HER framework to accelerate BCD methods for NTF cannot be guaranteed to converge using current existing proofs in the literature.
The main distinction of HER that makes it very efficient but difficult to prove convergence is its \emph{dynamic and flexibility}  in the update of extrapolation parameters. Specifically, HER performs a very aggressive extrapolation strategy: it extrapolates after each update of a block variable and uses that extrapolated point for the evaluation of the next block. Moreover, as long as the objective function decreases, it increases the extrapolation parameter $\beta$ (as discussed in the previous section). The accelerated block proximal gradient methods using extrapolation inside the block update such as APG and iBPG do not have such strategy. Hence, although these methods have strong convergence guarantee, the way they choose extrapolation parameters is more conservative. Our extensive experiments strongly confirm the efficacy of dynamic and flexibility in choosing extrapolation parameters of HER when extrapolation between each block update is used and monotonicity of the objective function is taken into account. Studying convergence guarantees for the HER framework would be very challenging and is an important direction of further research.

\subsection{Related works}\label{subsec:Bro}
Extrapolated AO algorithms can be traced back to a seminal work by Harshman \cite{Harshman1970}, in which extrapolation was seen as a way to speed up the convergence of ALS.
In this subsection we review two works on extrapolated AO algorithm: an older one \cite{Bro1998Multi} and a recent one \cite{mitchell2018nesterov}.
To differentiate the extrapolation parameter used in HER (denoted as $\beta_k$), we denote the extrapolation parameter used in these works as $\omega_k$.

\paragraph{General description}
Before we provide the details about the algorithms,
let us point out a very important observation:
unlike Algorithm~\ref{subsec:HER}, these algorithms use a global extrapolation parameter shared among all the blocks, and the extrapolation is conducted after all the blocks have been updated.
That is, they first update all the blocks variables, then stack all the blocks together to form a vector $x$, and extrapolate it as
\begin{equation}
	x_{k+1} = x_{k+\frac{1}{2}}  + \omega_k ( x_{k+\frac{1}{2}} - x_{k} ),
	\label{ext_BroGRLS_org}
\end{equation}
where $x_{k+\frac{1}{2}}$ is the vector obtained by stacking all the blocks $A^{(1)}_k,\dots,A^{(N)}_k$ after all of them have been updated.
The algorithms\cite{Bro1998Multi,mitchell2018nesterov} follow \eqref{ext_BroGRLS_org}, and they differ in the way $\omega_k$ is computed, which we discuss below.

\paragraph{Extrapolated AO algorithms with Bro's sequence}
Bro revisited and optimized the seminal work by Harshman~\cite{Harshman1970}, and came up with an extrapolation scheme with convincing empirical speed-ups for computing CPD:
the extrapolation parameter is tuned as
\[
\omega_k = k^{\frac{1}{h(k)}}-1,
\]
where $h(k)$ is a recursive function so that $h(k+1) = h(k)$ if the error has not increased for more than four iterations, $h(k+1) = 1+ h(k)$ otherwise, and $h(1) = 3$. Moreover, no extrapolation is performed in the first few iterations because of stability issues.
Furthermore, when the error increases, the extrapolation is not performed, that is, the extrapolated sequence is abandoned, as for HER.

Note that there is no particular modification of the Bro extrapolation scheme for aNCPD.
In this paper we implement Bro-AHALS, Bro-ADMM  and Bro-Nesterov -- the three versions of Bro's accelerated methods in which we respectively use the same strategy using A-HALS (see \S\,\ref{subsec:AHALS}), ADMM and Nesterov's accelerated gradient method for solving the NNLS problem \eqref{eq:NNLS} inexactly.

\paragraph{Extrapolated AO algorithms with gradient ratio and line search}

Recently two heuristic approaches similar to Bro's have been proposed \cite{mitchell2018nesterov}, and compute $\omega_k$ as follows
\begin{eqnarray}
	\textrm{Gradient ratio (GR):} &~& \omega_k = \nabla_{x} F_k(x) \big|_{x = x_{k+\frac{1}{2}}}  \;  \Big/ \; \nabla_x F_{k-1}(x) \big|_{x = x_{k-1}},
	\label{GRbeta}
	\\ \nonumber ~\\
	\textrm{Linear search (LS):}  &~& \omega_k = \underset{\omega}{\argmin} \, F\Big(x_{k + \frac{1}{2}} + \omega ( x_{k + \frac{1}{2}} - x_{k-1}) \Big),
	\label{LSbeta}
\end{eqnarray}
where $ \nabla_x F$ is the gradient of $F$ with respect to $x$.
As for Bro's accelerated algorithms, we implement in this paper GR-AHALS, GR-ADMM, GR-Nesterov, LS-AHALS, LS-ADMM and LS-Nesterov where we correspondingly use the same strategy as for A-HALS (see  \S\,\ref{subsec:AHALS}), ADMM and Nesterov's accelerated gradient method for \eqref{eq:NNLS}.

Note that GR and LS are designed for aCPD but not aNCPD, and similar to Bro's approach \cite{Bro1998Multi}, there is no modification of the GR and LS schemes for the nonnegative decomposition case.

\paragraph{The modified implementation of Bro, GR and LS}
The algorithms Bro, GR and LS use a global extrapolation parameter shared among all the blocks, which is different from the extrapolation parameter $\beta_k$ used in HER that is tuned independently for each block.
Preliminary tests have showed that HER is always speeding up BCD algorithms much faster than Bro, GR and LS (see Fig.~\ref{fig:exp:BCDs} in the next section).
Such superiority can be explained in part by the fact that the block-wise tuning of $\beta_k$ in HER gives HER much more degrees of freedom than Bro, GR and LS.
Hence, to make a fair comparison between the different extrapolation strategies,  we make the following modifications so that Bro , GR and LS have the same algorithmic structure as HER.

First, the update of $x_k$ is performed block wise, that is, one $A^{(i)}$ at a time.
Next, we extrapolate the blocks right after they have been updated, using the same extrapolation coefficients as described in Bro, GR and LS.
It is important to note that, in the original algorithms, all the block variables are extrapolated with the same ``global" extrapolation coefficient.
That is, the extrapolation coefficients for every block $A^{(i)}$ in the original algorithm are the same.
In the modification here, we ``split" the global extrapolation coefficient into block-specific extrapolation coefficient.
For example, in GR, the update-then-extrapolate step is performed for all $i$ as
\begin{eqnarray}
	\textrm{Update:} &~~~&
	\hspace*{-4mm}
	A^{(i)}_{k+\frac{1}{2}} \text{ using~\eqref{eq:NNLS}},
	\label{eqn:update2_x} \\ \hspace*{-1em}
	\textrm{Extrapolate:} &~~~&
	\hspace*{-4mm}
	A^{(i)}_{k+\frac{1}{2}}  + \dfrac{  \| \nabla_{A^{(i)}} F_k \| }{ \| \nabla_{A^{(i)}} F_{k-1} \| }  (A^{(i)}_{k+\frac{1}{2}} - A^{(i)}_{k}), %\vspace*{-2mm}
	\label{eqn:update2_y}
\end{eqnarray}
where $\nabla_{A^{(i)}} F_k$ is the gradient of $F$ with respect to block $A^{(i)}$ at iteration $k$ (see \eqref{def:grad}), and $A^{(i)}_{k+\frac{1}{2}}$ is the block $A^{(i)}$ at iteration $k$ just after the update, that is, we extrapolate the block right after it has been updated, as in HER.
Moreover, (\ref{eqn:update2_y}) uses the ratio between the norm of the gradient of the current block $A$ and the norm of the gradient of the same block in the last iteration.
We use the same strategy on splitting the global extrapolation coefficient into block-specific extrapolation coefficient in Bro and LS.
That is, the term $\frac{  \| \nabla_{A^{(i)}} F_k \| }{ \| \nabla_{A^{(i)}} F_{k-1} \| }$ in equation (\ref{eqn:update2_y}) is replaced by $\big(k^{\frac{1}{h(k)}}-1 \big)$ in Bro, and for LS, the extrapolation parameter $\omega_k$ is computed by solving a minimization subproblem: consider the update of the $i$th block at iteration $k$, we have
\begin{eqnarray}
	\omega_k &=& \argmin_{\omega} F\Big( A^{(1)}_k , \dots, A^{(i-1)}_k, A^{i}_k + \omega ( A^{i}_k - A^{i}_{k-1}), A^{(i+1)}_{k-1}, \dots,  A^{(N)}_{k-1} \Big ). %\vspace*{-2mm}
	\label{eqn:LS}
\end{eqnarray}
By expanding $F$ in terms of $\omega$,  (\ref{eqn:LS}) can be expressed as a second-order
polynomial in $\omega$, and hence a closed-form solution for $\omega$ exists.

\paragraph{Computational cost on Bro, GR and LS compared with HER}
The per-iteration cost in both GR and LS schemes is much larger than that of Bro.
Both Bro, GR and LS have restart, but Bro's extrapolation parameter is basically a scalar computation, while GR has multiple matrix-matrix multiplications and LS even has to solve a minimization sub-problem:
\begin{itemize}
	\item Here we solve \eqref{eqn:LS}, which is a second-order polynomial in $\omega$, exactly.
	\item For \eqref{LSbeta}, for a 3rd-order NCPD problem, we need to minimize a sixth-order polynomial in $\omega$.
	\item In general, for a $N$-order NCPD problem, we need to minimize a $2N$-order polynomial in $\omega$.
	Let $x = [ A^{(1)}, \dots ,A^{(N)} ] $ denotes the stacking of the block into vector\footnote{Here it is not the Kruskal notation \eqref{eq:kruskal}.}, then
	\begin{eqnarray}
		\beta_k & = &
		\underset{\omega}{\arg\min}
		F \Big( x_{k + \frac{1}{2}} + \omega ( x_{k + \frac{1}{2}} - x_{k-1} ) \Big) \nonumber \\
		& = &
		\underset{\omega}{\arg\min}
		F \Bigg(  \big [ A^{(1)}_{k + \frac{1}{2}}, \dots ,A^{(N)}_{k + \frac{1}{2}} ]
		+ \omega \Big( \big[   A^{(1)}_{k + \frac{1}{2}}, \dots ,A^{(N)}_{k + \frac{1}{2}} \big ] - [A^{(1)}_{k-1}, \dots, A^{(N)}_{k-1}] \Big) \Bigg) \nonumber \\
		& \overset{\eqref{eq:atensp}}{=} &
		\underset{\omega}{\arg\min}
		\dfrac{1}{2} \Bigg \|
		\mathcal T -
		\left({\bigotimesa_{i=1~}^{N~}}  	\Big( A^{(i)}_{k + \frac{1}{2}} + \omega \big( A^{(i)}_{k + \frac{1}{2}} - A^{(i)}_{k-1} \big) \Big)   \right)\tI_r
		\Bigg \|_F^2.
		\nonumber
	\end{eqnarray}
	We can see the cost of computing the coefficients for the polynomial can potentially be very high.
	Due to such reason, in the original paper, the $\omega_k$ in \eqref{LSbeta} is solved approximately using cubic line search in the Poblano toolbox.
	
	\item As pointed out in \cite{ang2019accelerating}, exact line search has bad performance in NMF, so we believe this is the same for LS, for both the original form and the modified form.
\end{itemize}
In general, the per-iteration cost of the extrapolation step in Bro's extrapolation is negligible, while for GR is higher than one ALS, and for LS it is much higher than one ALS.
For this reason, we only run the modified LS in the numerical tests, i.e., \eqref{eqn:LS}, and we will see that, often LS performs the worst in the experiments, see for example Fig.~\ref{fig:exp:BCDs}.

\paragraph{Remarks on Bro, GR and LS compared with HER}
In the numerical experiments, we will compare the original form of Bro, GR, LS, as well as the modified version.
There are several remarks on Bro, GR and LS.
\begin{itemize}
	\item There are two sequences $A^{(i)}$ and $\hat{A}^{(i)}$ used in HER, while there is no auxiliary sequence in Bro, GR and LS: that is, the extrapolation in these approaches is conducted on the same block $A^{(i)}$.
	
	\item %The implementations of the modified Bro, GR and LS in this paper are different from the original one proposed where Bro, GR and LS.
	By splitting of the extrapolation coefficient into block extrapolation coefficients, the original GR and LS are improved as the $\omega_k$ in the new GR and LS are more adapted to each block variable.
	
	\item As Bro, GR and LS are designed for aCPD but not aNCPD, and similar to Bro's approach~\cite{Bro1998Multi}, there is no modification of the GR and LS schemes for aNCPD.
	This means that there is no guarantee on feasibility of the iterates produced by these methods for aNCPD.
\end{itemize}

\section{Experiments}
\label{sec:exp}
In this section, we empirically prove the efficacy of HER by extensively test its performance on a rich set of synthetic data sets as well as real data sets. As presented in \S\,\ref{subsec:HER}, HER is a scheme to accelerate AO algorithms by using extrapolation between block update; and as such, by using HER, we can derive several different algorithms corresponding to the solver we use for the NNLS problem \eqref{NNLS_her_exact}. We stress out that HER can be used in combination with any BCD algorithm to make it faster.
In this section, we combined it with the most well-known algorithms to tackle NTF, namely AS, ADMM, Nesterov accelerated gradient and AHALS for solving  \eqref{NNLS_her_exact}, which we denote by HER-AS, HER-ADMM, HER-Nesterov and HER-AHALS, respectively.
We call HER-AO the set of these algorithms. 
Table~\ref{table:allalgs} lists the algorithms that we implement and test in our experiments. 
Note that our goal is not to compare these various algorithms, but to show that HER can accelerate all of them significantly.

\begin{table}
	\centering
	\stepcounter{table}% for \thetable
	\def\arraystretch{1.02}
	{
		\rowcolors{1}{red!0!green!0}{red!5!green!5}
		\begin{tabular}{ll}
			Algorithms   & Reference \\ \hline
HER-AS, HER-ADMM, HER-Nesterov, HER-AHALS &\S\,\text{\ref{subsec:HER}}  \\
AO-AS, AO-ADMM, AO-Nesterov, AHALS & \S\,\text{\ref{sec:AOframework}} \\
Bro-ADMM, Bro-Nesterov, Bro-AHALS & \S\,\text{\ref{subsec:Bro}}\\
GR-ADMM, GR-Nesterov, GR-AHALS & \S\,\text{\ref{subsec:Bro}}\\
LS-ADMM, LS-Nesterov, LS-AHALS & \S\,\text{\ref{subsec:Bro}}\\
APG, iBPG         & \S\,\text{\ref{subsec:proximal}} \\
		\end{tabular}
	}
	\addtocounter{table}{-1}%
	\caption{Algorithms for solving NTF}
	\label{table:allalgs}
\end{table}

All experiments are run with MATLAB (v.2015a) on a laptop with 2.4GHz CPU and 16GB RAM.
The code is available from \url{https://angms.science/research.html}.

\vspace*{2mm}
\noindent\textbf{Remark} ~~
In this paper, we focus on dense NTF problems, for which the input tensor has mostly positive entries. However, the HER framework can also be applied to sparse tensors. This was in fact done for NMF in~\cite{ang2019accelerating} with similar conclusions, that is, HER can accelerate algorithms significantly for sparse data sets as well. 
The problem with sparse data is not the algorithm itself, but rather its implementation. Handling sparse data also means dealing with extremely large dataset, which we cannot deal with our current implementation in Matlab. This could be fixed by integrating the proposed HER framework within a toolbox which features efficient sparse tensor manipulations and contractions, see for instance\cite{smith2015tensor} and the references therein. 
We leave an efficient implementation of the HER framework for very sparse and large tensors for future works.

\subsection{Set up}\label{sec:exp:setup}
\paragraph{Performance measurement}
Two important factors in the evaluation of the  performance of an algorithm are the data fitting error and the factor fitting error that are computed as follows. We use the value of the objective function
\begin{eqnarray}
	f_k := F\left(A_k^{(1)},A_k^{(2)},\ldots,A_{k}^{(N-1)},A_{k}^{(N)} \right)
	\label{err_true_def}
\end{eqnarray}
to represent the data fitting error. Supposing the ground truth factor matrices $A^{(i)}_\text{true}$, $i=1\ldots,N$ are available, then we compute the factor fitting error $e_k$  as
\begin{eqnarray}
	e_k := \dfrac{1}{N} \sum_{i = 1}^N \dfrac{ \norm{ \text{normalize}(A^{(i)}_\text{true}) - \text{normalize}(A^{(i)}_k) \Pi }_F  }{ \norm{\text{normalize} (A^{(i)}_\text{true}) }_F }.
	\label{err_fact_def}
\end{eqnarray}
Here $\text{normalize}(\cdot)$ is the column-wise normalization step (i.e., the $i$-th column of $\text{normalize}(A)$ is set to   $\frac{A(:,j)}{\| A(:,j)\|_2}$), and $\Pi$ is the permutation matrix computed through the Hungarian algorithm.
The use of $\Pi$ is to remove the permutation degree of freedom for matching the columns of $A^{(i)}$ to the column of $A^{(i)}_\text{true}$, and the use of normalization is to remove the scaling degree of freedom for matching the columns of $A^{(i)}$ to the column of $A^{(i)}_\text{true}$.

\paragraph{Generate a synthetic data}
To generate a synthetic tensor, we first generate ground truth factor matrices $A^{(i)}_\text{true} \in \mathbb{R}_+^{I_i \times r}, i = 1,\dots N$ whose entries are sampled from i.i.d.\ uniform distributions in the interval $[0,1]$. 
The tensor $\mathcal T^\text{clean} \in \mathbb{R}_+^{I_1 \times \dots \times I_N}$ is then constructed from  $A^{(i)}_\text{true}$, $i = 1,\dots N$. 
Finally, we form a synthetic data  $\mathcal T$ by adding some noise to $\mathcal T^\text{clean}$, $\mathcal T = \max( 0, \mathcal T^\text{clean} + \sigma \mathcal E)$, where $\sigma \geq 0$ is the noise level, and $\mathcal E \in \mathbb{R}^{I_1 \times \dots \times I_N}$ is a tensor whose entries are sampled from a unitary centered normal distribution.

\paragraph{Initialization, number of runs and plots}
For each run of an algorithm, we use a random initialization, i.e., the initial factor matrices $A^{(i)}_0$, $i=1,\ldots,N$, are generated by sampling uniform distributions in [0,1].  Note that, testing a specific data tensor, we use the same initialization in one run of all algorithms.  
We run all the algorithms 20 times with 20 different initializations. We stop one run of an algorithm when the maximum time (which is chosen before running the algorithms) is reached.

In presenting the results, we plot $f - f_{\min}$; and if the ground truth is known, we also  plot $e - e_{\min}$. Here $f_{\min}$ and $e_{\min}$ are respectively the minimal value of all the data fitting errors and the factor fitting errors obtained across all algorithms on all runs. In noiseless settings $(\sigma = 0)$, exact factorization is possible, so we set $f_{\min} = 0$.
In order to have a better observation of the performance of the algorithms, we plot the curves with respect to both time and iterations \footnote{We do not report the
	number of MTTKRP (Matricized tensor times Khatri-Rao product) as all the algorithms in the experiments (except for AS) share the same number of MTTKRP (which is $N$ for an tensor with order $N$), so the performance in terms of number of MTTKRP is contained implicitly in the plot with respect to the iterations.}. 
We remark that, ``an iteration'' for AO algorithms means the counter $k$ of the outer loop after all blocks being updated. Regarding the time evaluation, we record the time stamp for each iteration, and then perform a linear interpolation to synchronize the time curves.
Note that such linear interpolation does not reflect 100\% truly the real convergence behaviour as it is just an linear estimate, but we consider such estimate to be accurate enough.

In our experiment, we emphasize on plotting the median curves of the 20 runs (which are the thick curves in the upcoming figures), because there may be large deviations between different runs.

\paragraph{Solving the NNLS problem \eqref{eq:NNLS} and \eqref{NNLS_her_exact}}
When using AS, ADMM, Nesterov's accelerated gradient descent algorithm or AHALS to solve  \eqref{eq:NNLS} or \eqref{NNLS_her_exact}, in our implementation, we terminate the solver when the number of iterations reaches 50 or when  $\| A^{(i)}_{s} - A^{(i)}_{s-1} \| \leq 10^{-2} \| A^{(i)}_{0} - A^{(i)}_{1} \|$, where $s$ is the iteration counter of the solver.

\paragraph{Parameter set up for HER}
We use the following set of parameters for HER-AO (unless otherwise specified):
$ \beta_0 = 0.5, \gamma = 1.05, \bar{\gamma} = 1.01, \eta = 1.5 $, and (\ref{extproj_her}) is used for the extrapolation point.

\paragraph{List of experiments}
Table \ref{table:exp} lists the figures that report our diverse experiments on synthetic data and real data sets.
All the experiments have $N=3$ and the input tensor is dense.

\paragraph{Complete numerical experimental results}
We refer the interested reader to the appendix of the arXiv version of the current paper \cite{2020arXiv200104321A} for all the numerical experiments as reported in Table \ref{table:exp}. The conclusions remain the same: HER significantly accelerate the convergence of BCD algorithms, while the HER-BCD outperforms both iBPG and APG.

\begin{table}[ht!]
	\centering \stepcounter{table}% for \thetable
	\def\arraystretch{1.0}
	{\rowcolors{1}{red!0!green!0!blue!0}{red!5!green!5!blue!5}
		\begin{tabular}{clc}
			Fig.  & Test description    & $[I_1,I_2,I_3,r,\sigma]$   \\ \hline
			\multicolumn{3}{c}{Synthetic data} \\ \hline
			\text{\ref{fig:exp:bal_and_unbal}}&Cube size, low rank, noiseless & $[50,50,50,10,0]$\\
			\text{\ref{fig:exp:bal_and_unbal}}&Unbalanced size, low rank, noiseless & $[150,10^3,50,12,0]$\\
			\text{\ref{fig:exp:highrank}}     &Unbalanced size, larger rank, noiseless & $[150,10^3,50,25,0]$\\
			\text{\ref{fig:exp:BigNoisy}}     &Large cube size, low rank, noisy & $[500,500,500,10,0.01]$\\
			\text{\ref{fig:exp:cond3_150_1000_50_12_n001}} &Unbalanced size, low rank, noisy, ill-condition & $[150,10^3,50,12,0.001]$  \\
			\text{\ref{fig:exp:pg_apg_herpg}} &HER-AO-gradients compared with APG and iBPG & $[150,10^3,50,10,0.01]$ \\
			\text{\ref{fig:exp:BCDs}}         & Comparing $\{$HER,Bro,GR,LS$\}$-AHALS  &
			$\begin{matrix}
			[50,50,50,10,0]	\\
			[150,10^3,50,12,0.01]	\\
			[150,10^3,50,25,0.01]
			\end{matrix}$ \\ \hline
			\multicolumn{3}{c}{Real data}  \\ \hline
			\text{\ref{fig:exp:HSI}} & Two HSI images : PaviaU and Indian Pine &
			$\begin{matrix}
			[610,340,103,10]	\\
			[145,145,200,15]
			\end{matrix}$ \\
			\text{\ref{fig:exp:real_video}} & Big data : black-and-white video sequence  &
			$[153,238,1.4\times 10^4,\{10,20,30\}]$
			\vspace*{-2mm}
	\end{tabular}}
	\addtocounter{table}{-1}%
	\caption{List of experiments on NTF.} \label{table:exp}
\end{table}

%\subsection{Experiments on synthetic data: effectiveness of HER}
\subsection{Experiments on synthetic data sets}
\label{sec:exp:exp}
%\paragraph{On AO vs HER-AO}
As listed in Table \ref{table:exp}, the experiments on synthetic data sets are designed to simulate different kinds of situations that may occur in real applications, which includes : low rank, larger rank, noiseless, noisy, tensor with balanced size (cubic tensor), tensor with unbalanced size (rectangular tensor), and  ill-conditioned tensor.

\begin{center}\begin{figure}[ht]
		\centering
		\includegraphics[width=\linewidth]{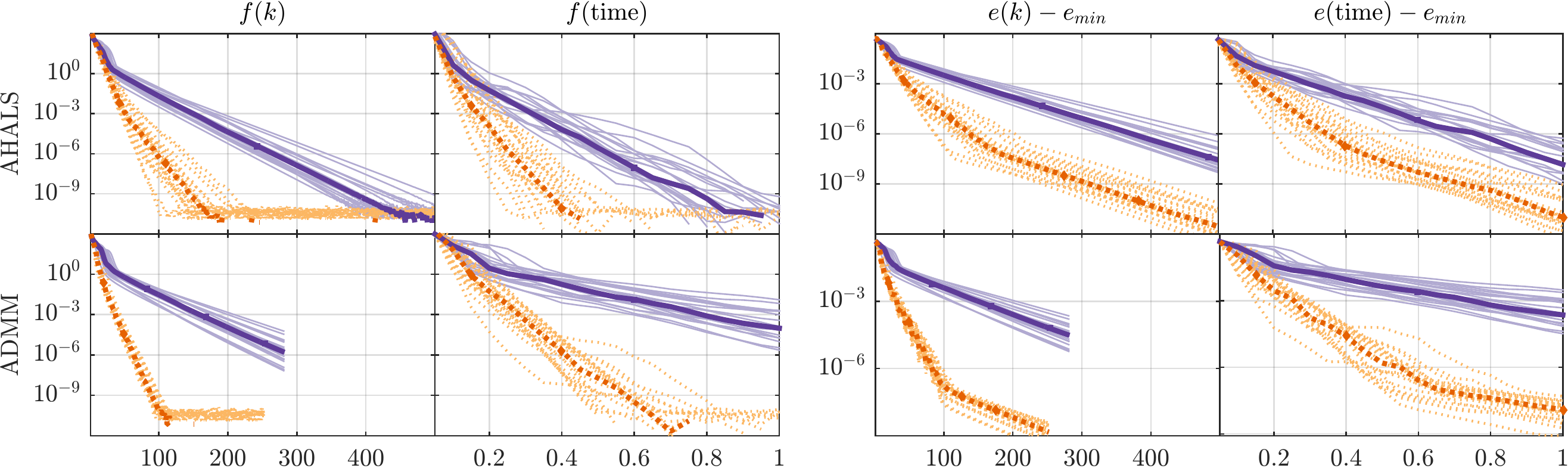}
		\includegraphics[width=\linewidth]{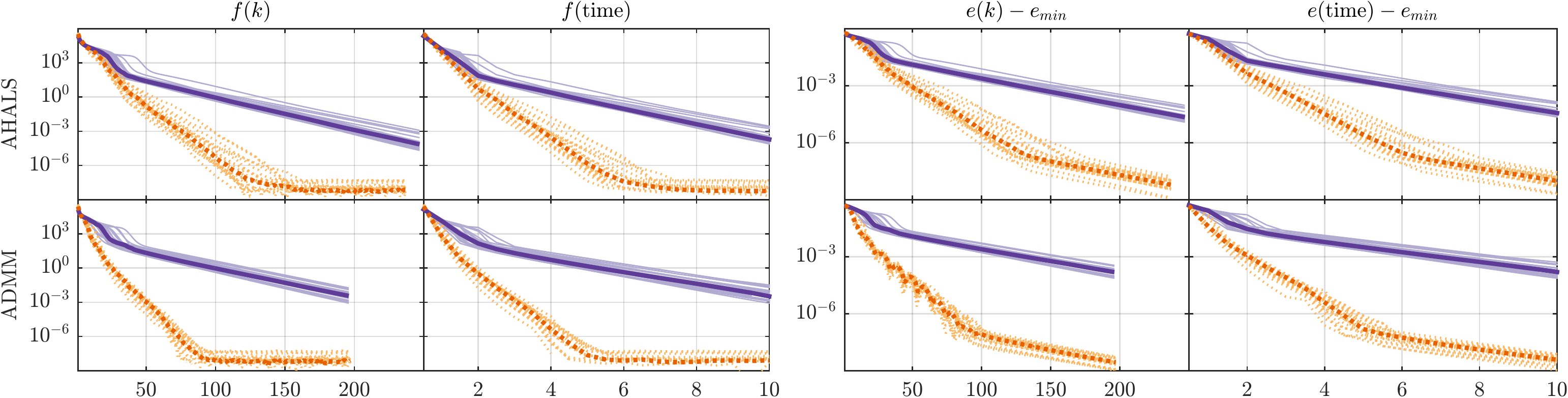}\vspace*{-2mm}
		\caption{Convergence of algorithms : A-HALS and AO-ADMM without HER (solid purple) and with HER (dotted orange), on standard test case (top) : $[I_1,I_2,I_3,r] = [50,50,50,10]$ and
			unbalanced sizes (bottom) $[I_1,I_2,I_3,r] = [150,10^3,50,12]$.
			The results show that HER improves the convergence significantly, the convergence in both $f$ and $e$ for HER-accelerated methods are already multiple-order of magnitude better than the un-accelerated algorithms.
			Notice that due to a higher per-outer-iteration cost, ADMM-based algorithms (AO-ADMM and HER-AO-ADMM) run fewer number of outer-iteration than the AHALS-based algorithms.
			See Appendix~\ref{appendix:allexps} for the results on other algorithms where we observe a similar behaviour.
		}
		\label{fig:exp:bal_and_unbal}
\end{figure}\end{center}

\begin{center}\begin{figure}[p]
		\centering\includegraphics[width=\linewidth]{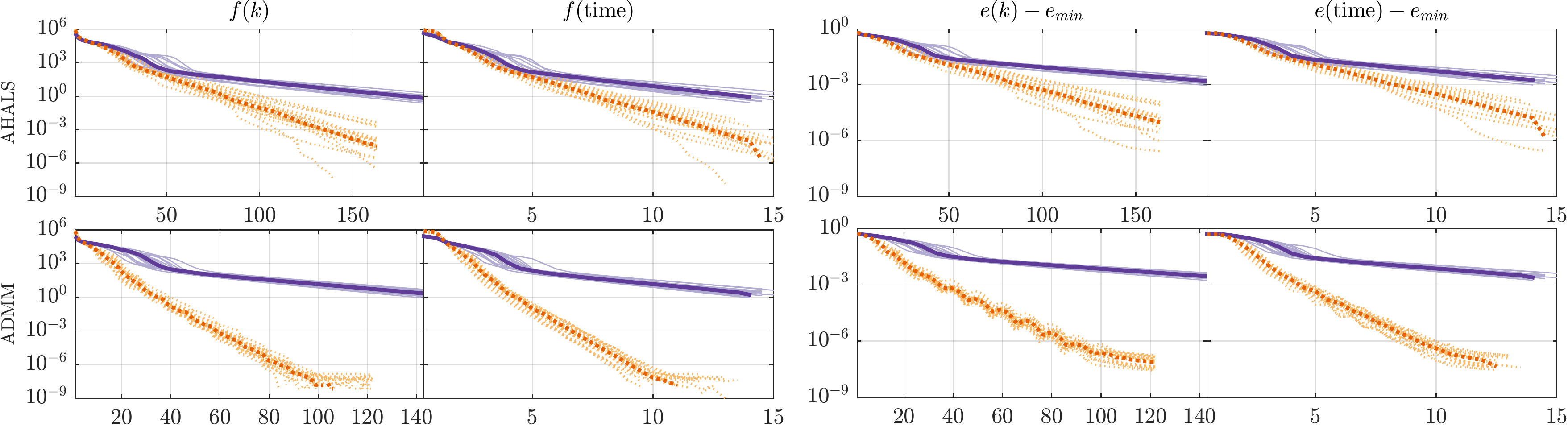}\vspace*{-2mm}
		\caption{On large rank $[I_1,I_2,I_3,r]=[150,10^3,50,25]$.
			For the plot set up, see Fig. \ref{fig:exp:bal_and_unbal}.
			Results show HER improves the convergence speed significantly.
			See Fig. \ref{fig:exp:bal_and_unbal} for the plot set up, and Appendix~\ref{appendix:allexps} for the results on other algorithms.
		}\label{fig:exp:highrank}
\end{figure}\end{center}

\begin{center}\begin{figure}[p]
		\centering\includegraphics[width=\linewidth]{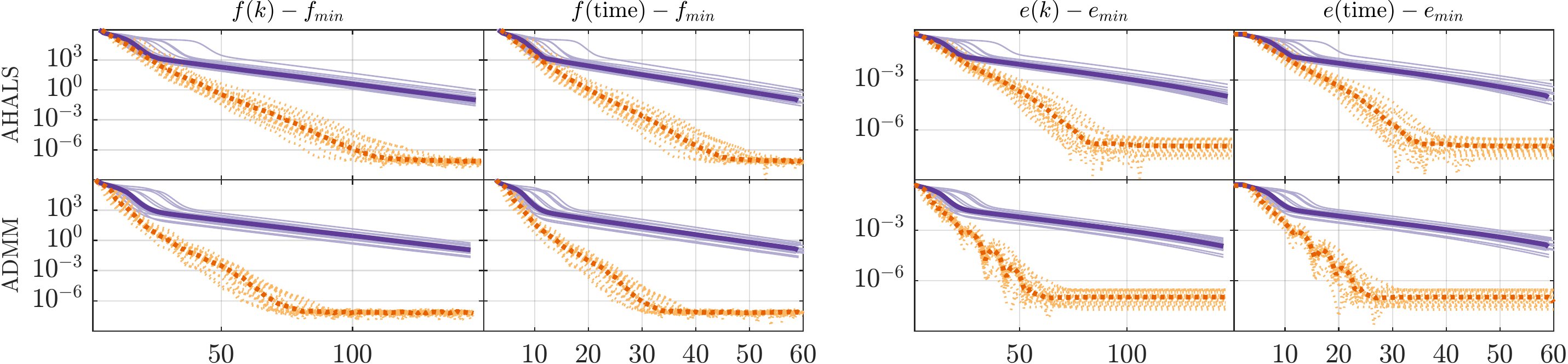}\vspace*{-2mm}
		\caption{On big and noisy case $I_1=I_2=I_3=500, [r,\sigma] = [10,0.01]$.
			Results show HER improves the convergence speed significantly.
			See Fig. \ref{fig:exp:bal_and_unbal} for the plot set up, and Appendix~\ref{appendix:allexps} for the results on other algorithms.
		}\label{fig:exp:BigNoisy}
\end{figure}\end{center}

\begin{center}\begin{figure}[p]
		\centering\includegraphics[width=\linewidth]{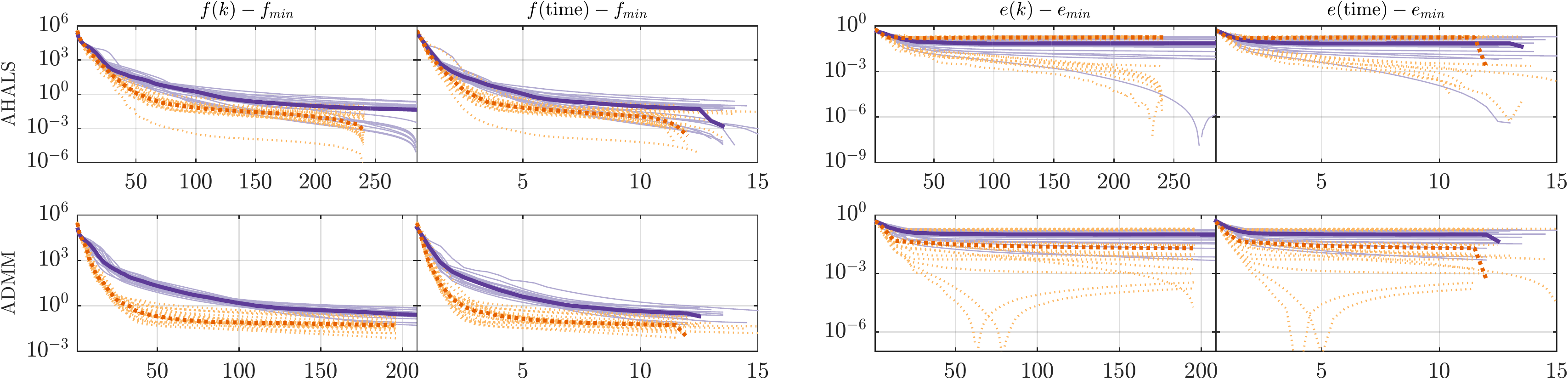}\vspace*{-2mm}
		\caption{On ill-conditioned case
			$[I_1,I_2,I_3,r,\sigma] = [150,10^3,50,12,0.01]$, where $A_i(:,1) = 0.99 A_i(:,2) + 0.01 A_i(:,1)$ for $i\in \{1,2,3\}$.
			For the plot set up, see Fig. \ref{fig:exp:bal_and_unbal}.
			Results show HER improves the convergence speed.
			See Fig. \ref{fig:exp:bal_and_unbal} for the plot set up, and Appendix~\ref{appendix:allexps} for the results on other algorithms.
		}\label{fig:exp:cond3_150_1000_50_12_n001}
\end{figure}\end{center}

Figure \ref{fig:exp:bal_and_unbal}, \ref{fig:exp:highrank} \ref{fig:exp:BigNoisy} and \ref{fig:exp:cond3_150_1000_50_12_n001} strongly affirm that HER-ADMM and HER-AHALS significantly outperform their counterparts AO-ADMM and AHALS in term of both $f_k$ and $e_k$. We stress that the improvement is often of several orders of magnitude (at least $10^4$  in most cases). We observe the same result for HER-AS and HER-Nesterov vs AO-AS and AO-Nesterov. The full experiments can be viewed in Appendix~\ref{appendix:allexps}.

%In general, the results show that HER-AO is has a much better convergence performance than the un-accelerated counterparts.

%Furthermore, HER-AO often has better convergence performance, compared with APG and iBPG.
%Note that the better convergence performance on HER-AO are reflect in both $f_k$ and $e_k$ on both outer-iteration $k$ and computational time, and the improvement is often in many order of magnitude (says $10^4 - 10^12$).
%Due to space limit, we only plot the results on HER-AHALS and HER-ADMM here.
%For the complete plots listing all the results including APG and iBPG, see appendix.

\begin{center}\begin{figure}[p]
		\centering\includegraphics[width=\linewidth]{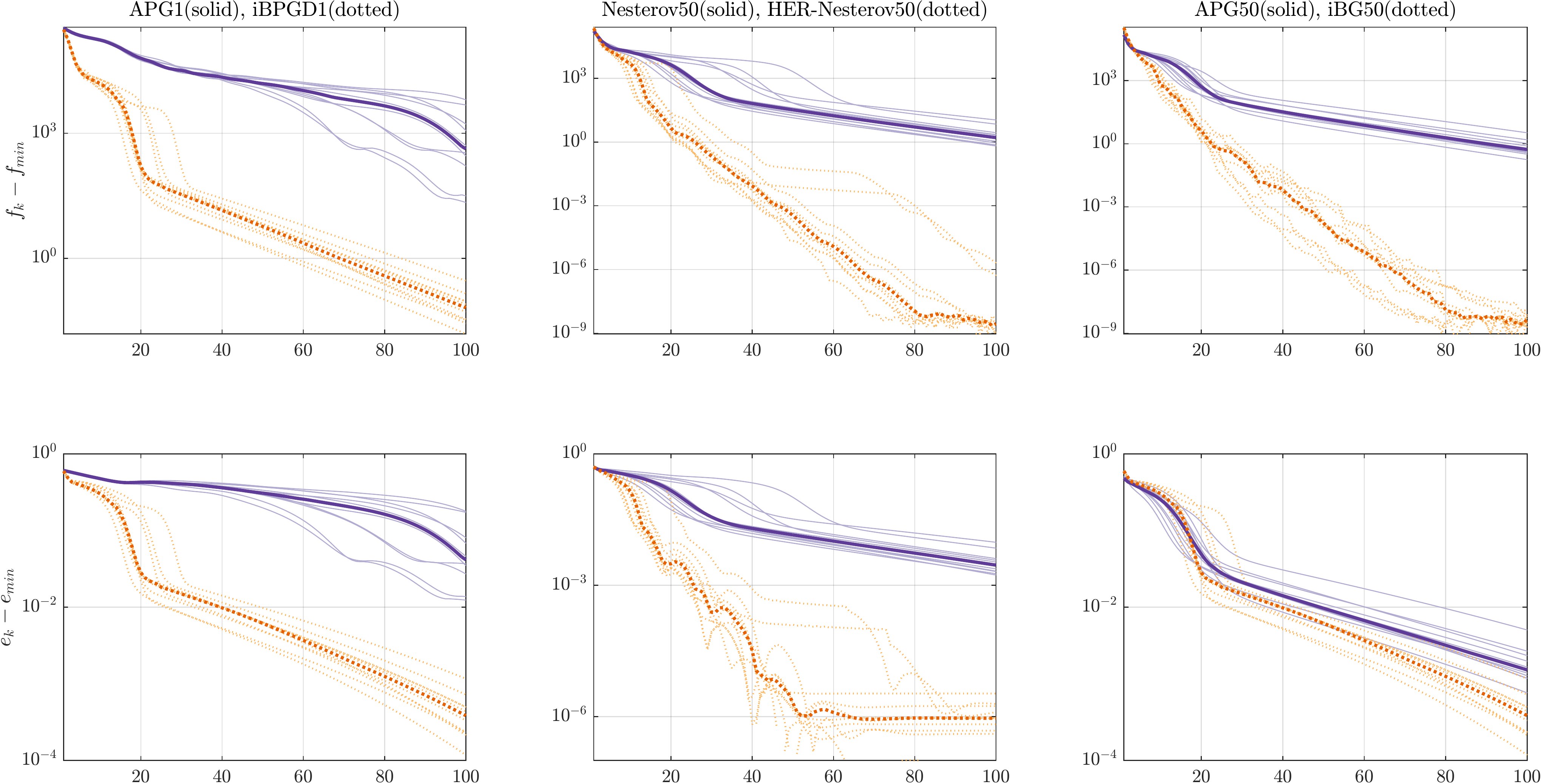}\vspace*{-2mm}
		\caption{Comparing gradient algorithms on $[I_1,I_2,I_3,r,\sigma] = [150, 10^3, 50, 10, 0.01]$.
			Suffix number denotes the maximum number of inner iterations.
			Result shows HER works for both inexact and exact BCD using gradient.
			Here HER-Nesterov50 and HER-PGD50 are the best algorithms in both $f$ and $e$.
			We do not plot the time plot here as they are similar to the iteration plot.}
		\label{fig:exp:pg_apg_herpg}
\end{figure}\end{center}

%\paragraph{On gradient-based methods}
Compared with APG and iBPG, we observe from Fig.~\ref{fig:exp:pg_apg_herpg} that HER-Nesterov outperforms both APG and iBPG in term of $f$ and significantly outperforms them in term of $e$. From extensive experiments (see Appendix~\ref{appendix:allexps} for more results),
we observe that HER, the scheme that makes use of the extrapolation between block update scheme, shows much better performance than APG and iBPG, the accelerated block proximal gradient methods that use Nesterov-type extrapolation inside each block update.

\begin{figure}[ht!]
	\centering
	\begin{subfigure}[b]{0.326\textwidth}
		\centering\includegraphics[width=\textwidth]{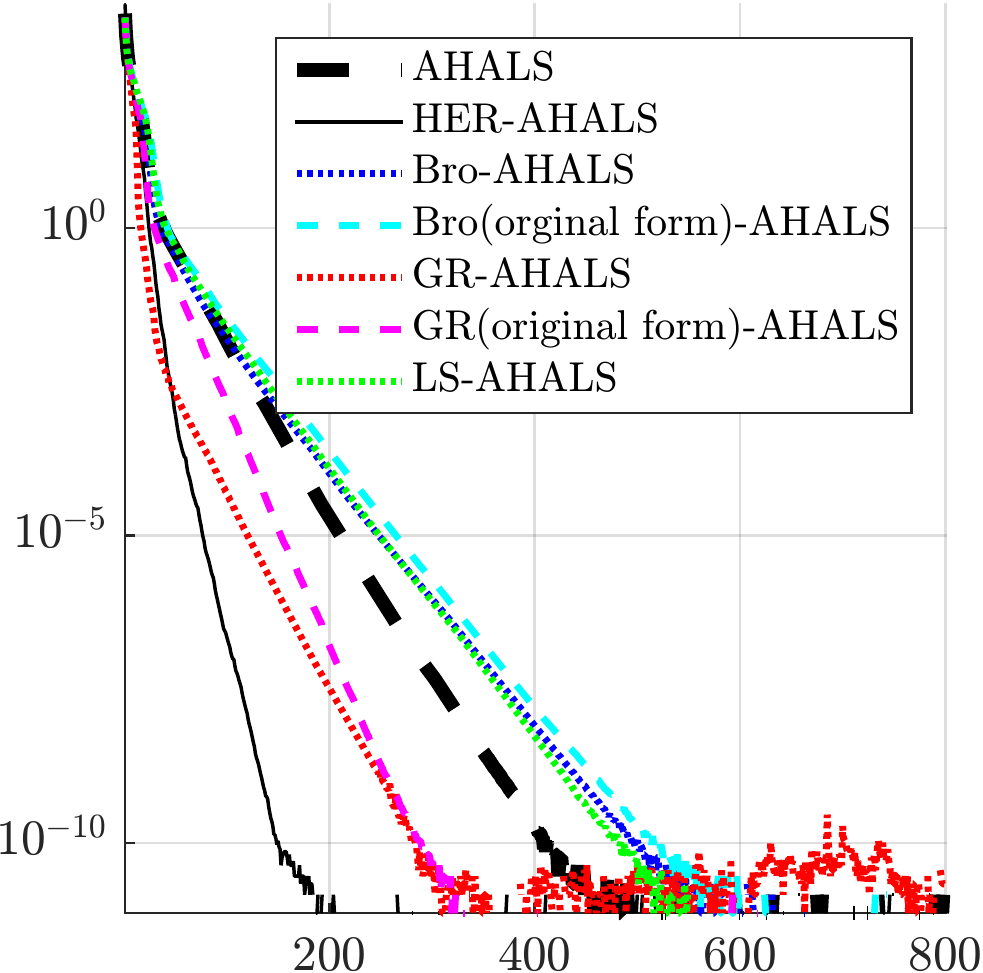}
		\caption{$[50,50,50,10,0]$}\label{fig_expBH50_50_50_10_n0}
	\end{subfigure}	\hfill
	\begin{subfigure}[b]{0.326\textwidth}
		\centering\includegraphics[width=\textwidth]{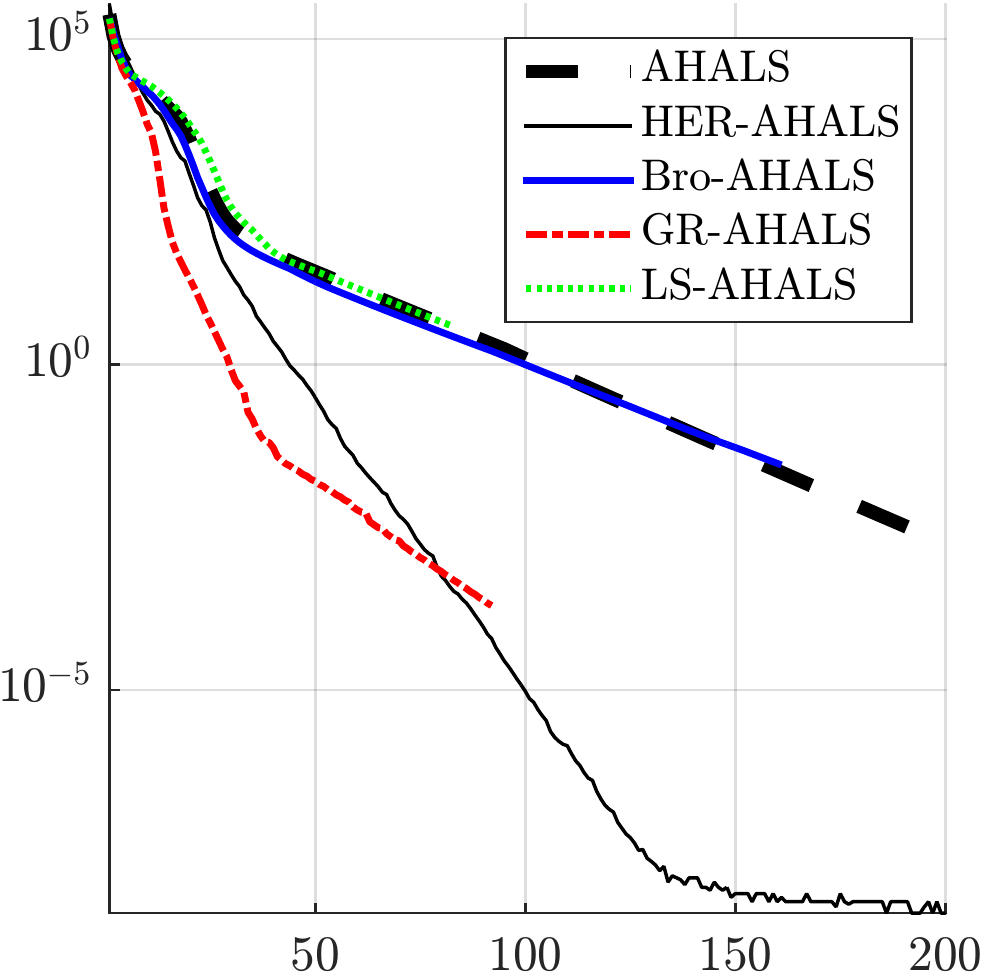}
		\caption{$[150,10^3,50,12,0.01]$}\label{fig_expBH150_1000_50_12_n001}
	\end{subfigure}	\hfill
	\begin{subfigure}[b]{0.326\textwidth}
		\centering\includegraphics[width=\textwidth]{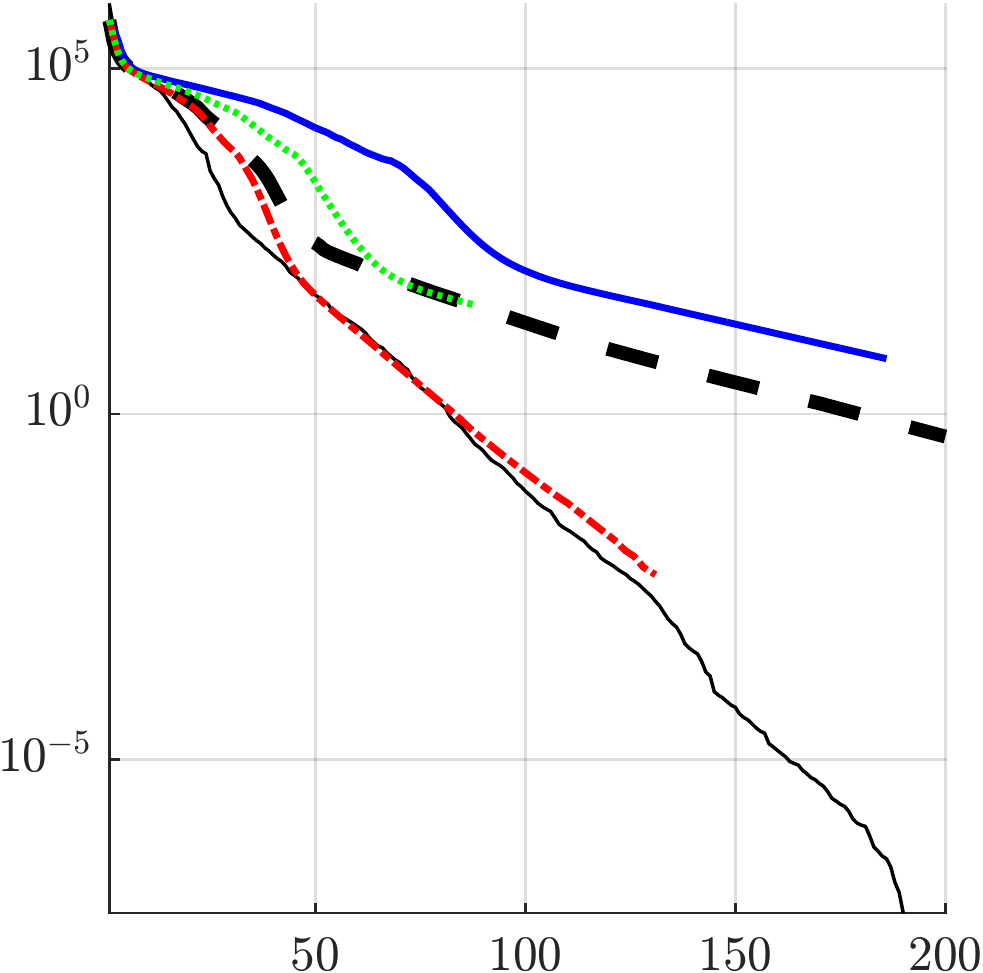}
		\caption{$[150,10^3,50,25,0.01]$}\label{fig_expBH150_1000_50_25_n001}
	\end{subfigure}	\hfill \vspace*{-2mm}
	\caption{Comparing AHALS with different acceleration frameworks on synthetic datasets on 3 setting of $[I_1,I_2,I_3,r,\sigma]$.
		The curves are the median in $f(k)-f_{\min}$.
		The x-axis is the number of iteration, and all algorithm run with same run time limited.
		These results show that that (1) HER-AHALS performs better than all other algorithms, and (2) the modified Bro and GR algorithms perform better than their original counterpart; see subplot(a), 
		while 	we do not show the result by Bro and GR in their original form in other subplots because they perform worse.  
		Note that LS and GR run less number of iterations due to their larger per-iteration cost.
		Bro's approach has lower per-iteration cost, but it is even slower than vanilla AHALS.	
		See the Appendix for more results.
	}
	\label{fig:exp:BCDs}
\end{figure}

Compared with Bro-AHALS, GR-AHALS and LS-AHALS, Fig.\,\ref{fig:exp:BCDs} shows that our HER-AHALS performs the best in the three experimental settings (only median is plotted here).
%Results show HER performances better than Bro, GR and LS.
Note that since the acceleration frameworks Bro, GR and LS are not designed for NTF, it is possible the iterates produced by these frameworks are infeasible.
Here we only compare HER-AHALS with Bro- AHALS, GR- AHALS and LS- AHALS; the comparison of these methods where AHALS is replaced with AO-ADMM and AO-ADMM are available in Appendix~\ref{appendix:allexps}, and similar conclusions are drawn, namely that HER outperforms the other accelerations.

%\hien{Hien: why don't we compare with other Bro,GS, LS - Nesterov and ADMM?} %all iterates (main sequence and the sequence of the extrapolation point) are feasible due to (\ref{extproj_her}).
%\andersen{Then why don't we compare with Bro,GS,LS,-AS? There will be many many figures, I think the use of AHALS here is representative enough. In fact, I did plot the comparisons on Bro,GS,LS-Nesterov, and also on ADMM, but Nicolas suggest to remove them.}
% * * *  * * *  * * *  * * *  * * *
% Balanced and Inbalanced Size
% * * *  * * *  * * *  * * *  * * *
% * * *  * * *  * * *  * * *  * * *
% unbalanced high rank noisless
% * * *  * * *  * * *  * * *  * * *
% * * *  * * *  * * *  * * *  * * *
% Big Balanced Noisy
% * * *  * * *  * * *  * * *  * * *
% * * *  * * *  * * *  * * *  * * *
% Unbalanced illcondition
% * * *  * * *  * * *  * * *  * * *
% * * *  * * *  * * *  * * *  * * *
%  Gradients : PGD vs Nesterov vs HER-PGD vs HER-Nesterov, APG, iBPG
% * * *  * * *  * * *  * * *  * * *
% %%%%%%%%%%%%%%%%%%%%%%%%%%%%%%%%%%%%%
% HER VS Bro GR LS
% %%%%%%%%%%%%%%%%%%%%%%%%%%%%%%%%%%
% % % % % % % % % % % % % % % % % % % % % % % % % % % % % % % % % %
% Real data
% % % % % % % % % % % % % % % % % % % % % % % % % % % % % % % % % %
\subsection{On real data}\label{sec:exp:exp_real}
\paragraph{Two hyper-spectral images}
We test the performance of the algorithms on two hyperspectral images (HSI)  PaviaU and Indian Pines\,\footnote{Data available from \url{http://www.ehu.eus/ccwintco/index.php?title=Hyperspectral_Remote_Sensing_Scenes}}.
They are nonnegative 3rd-order tensor; PaviaU has size $[610,340,103]$ with $r=10$ and Indian Pines has size $[145,145,200]$ with $r=15$. The $r$ chosen are commonly used in practice.

\begin{center}\begin{figure}[ht!]
		\centering
		\begin{subfigure}[b]{0.49\textwidth}
			\centering\includegraphics[width=\linewidth]{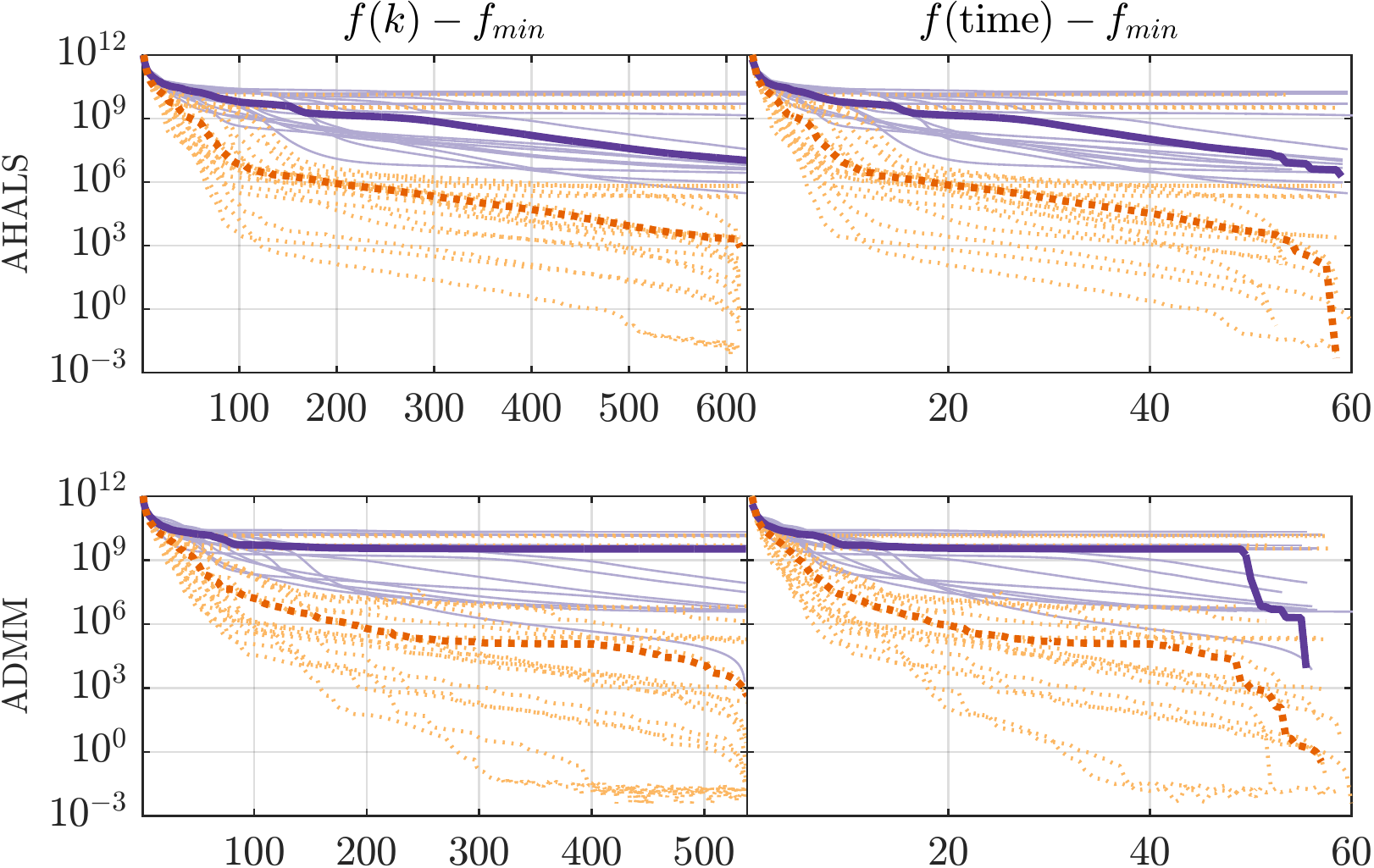}
			\caption{On PaviaU dataset.}\label{fig:exp:real_PaviaU10}
		\end{subfigure}
		\hfill
		\begin{subfigure}[b]{0.495\textwidth}
			\centering\includegraphics[width=\linewidth]{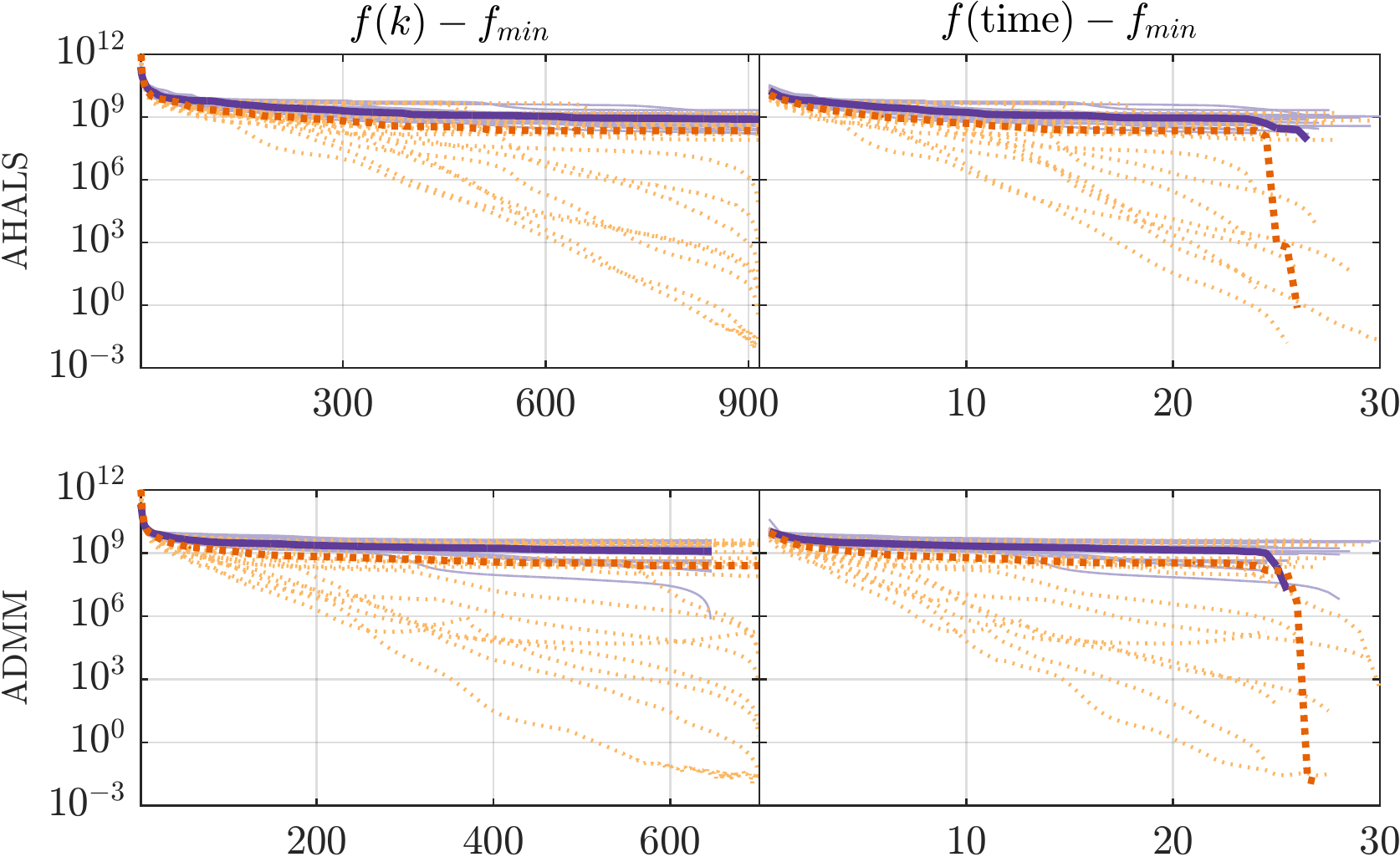}
			\caption{On Indian Pines dataset.}\label{fig:exp:real_IndianPines15}
		\end{subfigure}
		\caption{The results on HSI data.
			For the plot set up, see Fig. \ref{fig:exp:bal_and_unbal}.
			Results show HER improve convergences.
			See the Appendix of \cite{2020arXiv200104321A} for more results.}\label{fig:exp:HSI}
\end{figure}\end{center}

We perform minimal pre-processing on the raw data : NaN or negative values (if any) are replaced by zero.
Hence, it is possible the pre-processed data contains many zeros and  being ill-conditioned. Figure \ref{fig:exp:HSI} reports the performance of HER-AHALS, HER-ADMM and their counterparts AO-AHALS and AO-ADMM on the two data sets. As there are no ground truth factors, we only show $f$ in the results.

We observe that there are multiple swamps, which are common for real datasets as the data are highly ill-conditioned (the condition numbers of the metricized pre-processed data tensor along all modes are $[593,642,1009]$ for Indian Pines and $[944,462,8083]$ for PaviaU).
Nevertheless, considering the ``best case'' among the trials, HER-AHALS and HER-ADMM provide solutions with error $10^8 - 10^{10}$ times smaller than the best case of their un-accelerated counterparts. To compare with other algorithms, the readers can view the results in the Appendix of \cite{2020arXiv200104321A}. We observe that iBPG, APG and the AO (AO-AHALS and AO-ADMM) algorithms accelerated by GR, Bro and LS schemes are much slower than our AO  (AO-AHALS and AO-ADMM) algorithms accelerated by HER. GR-AO and Bro-AO  (for AO being AO-AHALS or AO-ADMM) even sometimes diverge.

\begin{center}\begin{figure}[p]
		\centering
		\begin{subfigure}[b]{0.328\textwidth}
			\centering\includegraphics[width=\textwidth]{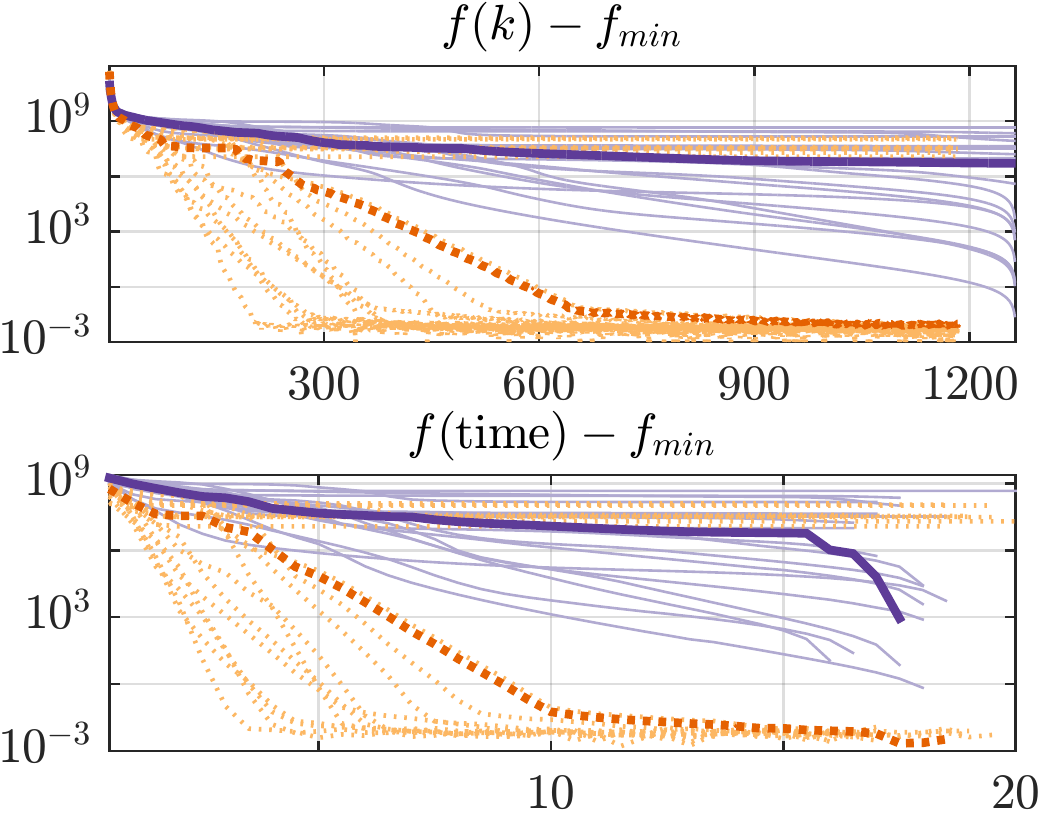}
			\caption{Case $r=10$.}
			\label{fig:real_video10}
		\end{subfigure}
		\hfill
		\begin{subfigure}[b]{0.328\textwidth}
			\centering\includegraphics[width=0.96\textwidth]{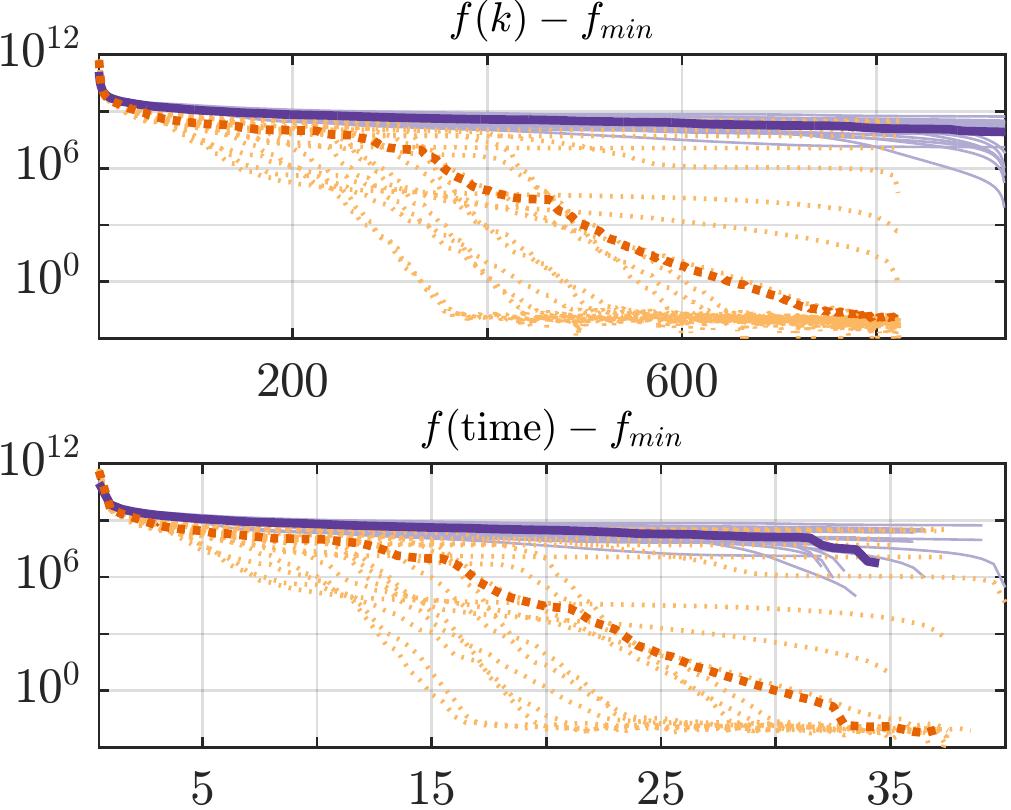}
			\caption{Case $r=20$.}
			\label{fig:real_video20}
		\end{subfigure}
		\hfill
		\begin{subfigure}[b]{0.328\textwidth}
			\centering\includegraphics[width=\textwidth]{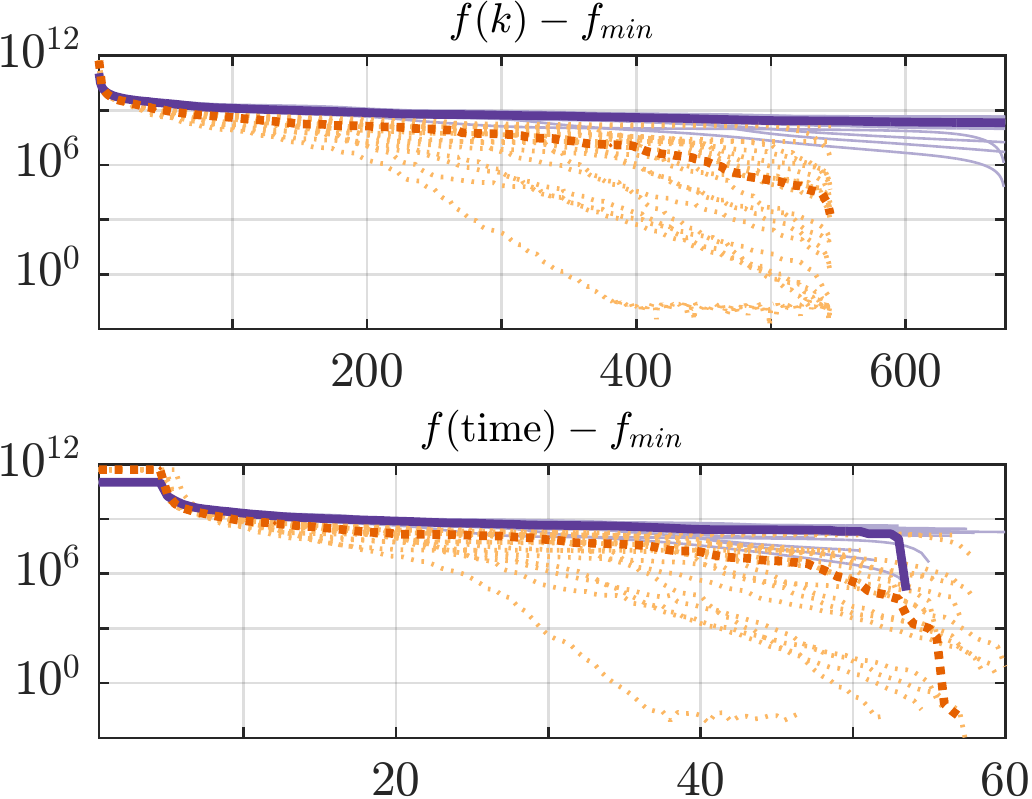}
			\caption{Case $r=30$.}
			\label{fig:real_video30}
		\end{subfigure}
		\caption{On video data $[153, 238, 14000]$ for three values of $r$.
			Results show HER improve convergence and works well with Tucker-based compression.}
		\label{fig:exp:real_video}
\end{figure}\end{center}

\paragraph{On big data : video sequences}
We test HER-AHALS on the video data of the UCSD Anomaly Dataset \cite{mahadevan2010anomaly}.
Constructed by combining all the frame images of 70 surveillance video in the dataset, we have a tensor with sizes $153 \times 238 \times 14000$, where the first two modes are the screen resolution and the third mode is the number of frame.
No pre-processing is performed on the raw data.
Data of such size is too big to store in our computer memory, so we perform compression using Tucker decomposition, based on the built-in function from the Tensor toolbox~\cite{Bader2007efficient}.
We compare AHALS and HER-AHALS with $r\in\{10,20,30\}$.
Results in Fig.\,\ref{fig:exp:real_video} shows that HER-AHALS performs much better than AHALS. For the details on how HER works with Tucker compression, see the Appendix of \cite{2020arXiv200104321A}.
% HSI : PaviaU(Left) IndianPines(Right)
%
% Video r 10 20 30
%
As a conclusion for this section, we give some remarks on HER-AO. From our extensive experiments, we observe that HER-AS has inferior performance than others when the data is either big in size, high rank, or ill-conditioned.  When the data has small size, all HER-AO algorithms have similar performance, and they all outperform their un-accelerated counterpart algorithms in term of both time and iteration. Among HER-AO algorithms, we highly recommend HER-AHALS for NTF as it shows good performance in all experiments.

\section{Discussion and conclusion}\label{sec:concl}
In this paper, we have proposed an extrapolation strategy in-between block updates, referred to as heuristic extrapolation with restarts (HER), for improving the empirical convergence speed of block-coordinate descent algorithms for approximate nonnegative tensor factorization (NTF).
HER significantly accelerates the empirical convergence speed of most existing block-coordinate algorithms for dense NTF, in particular for challenging computational scenarios, while requiring a negligible additional computational budget.
The core of HER is to apply a special extrapolation-restart mechanism that aims to reduce the computational cost of restart while making sure the restart criterion follows the standard function restarts.
The performance of HER was verified by the experiments reported in this paper.
In all scenarios, HER-AHALS provides among the best results hence we recommend its use in practice.

Future works include deriving theoretical convergence  for HER, and to apply it on other challenging applications.

\footnotesize
\bibliographystyle{vancouver}
\bibliography{NTF}
\clearpage

\section*{Appendix} 
\appendix
\section{Efficient tensor compression via Tucker format} \label{sec:compression}
Although there is a long history of using the Tucker model as a compression tool to pre-process big dataset, only recently has been formally discussed that compression does not actually imply transforming the large dataset into a smaller tensor~\cite{Vervliet2019exploiting}.
Given the tensor $\tT$, its Tucker format is expressed as :
\begin{equation}
\tT = \left( \bigotimesa_{p=1}^N U^{(p)} \right) \tG
\end{equation}
where $U^{(p)} \in \mathbb{R}^{n_p\times r_p}$, $\tG\in\mathbb{R}^{r_1\times
	\ldots\times r_N}$ and $\{r_p\}_{p\leq n}$ are inputs integer parameters of the format, sometimes called Tucker ranks~\cite{Hackbusch2012tensor}.
This representation is not unique but still offers a compressed expression of $\tT$ thus the name format rather than decomposition.

A typical situation is that of a tensor $\tT$ too big to fit in memory, since either too large and dense, or extremely large but sparse.
%In both cases, there is an extensive amount of literature on how to compute the Tucker format for
%such large tensors, see the introduction of \S~\ref{sec:exBCD}.
Therefore, a third party may instead provide the data directly in a compact format such as the Tucker format.
As Tucker format is in practice an approximation of the real data, the cost function of the aNCPD problem is modified as follows:
\begin{equation}
F_t(A^{(1)},\ldots,A^{(N)}) = \frac12
\norm{
	\left( \bigotimesa_{p=1}^N U^{(p)} \right) \tG
	- \left( \bigotimesa_{p=1~}^{N~} A^{(p)} \right) \mathcal{I}_r
}_F^2.
\label{eqn:objFunTucker}
\end{equation}

On top of the storage gain, there is a huge computational burden ease in using structured representations of the data when computing the MTTKRP.
Indeed, the gradient of $F_t$ wrt say $A^{(1)}$ is obtained as follows:
\begin{equation}
\nabla_{A^{(1)}} F_t = - U^{(1)} G_{[1]} \left[ \bigodot_{p=2}^{N}
(U^{(p)})^T A^{(p)} \right] + A^{(1)}\left[
\mathop{\scalebox{2}{\raisebox{-0.2ex}{$\circledast$}}}_{p=2}^{N}
(A^{(p)})^TA^{(p)} \right]
\end{equation}
where $\circledast $ is the Hadamard product.
The equation involves only ``cheap" products if Tucker ranks $r_p$ are small compared to the data tensor dimensions $n_p$.
In \S~\ref{sec:exp:exp_real}, we check that indeed herBCD is compatible with accelerating the aNCPD using the Tucker format, and this actually opens the door to many problems that could not be tackled with simply herBCD, while enhancing at no cost the convergence speed of BCD algorithms for minimizing $F_t$.
This contrasts with usual developments of fast techniques to solve aNCPD that typically do not consider other kind of acceleration in conjunction.

%\newpage

\section{Algorithm pseudocodes for APG and iBPG}\label{sec:algo_supp}
In this section, we provide the details for the implementations of APG and iBPG.

%****************************************************
\begin{minipage}{0.95\linewidth}
	%\algsetup{indent=2em}
\begin{algorithm}[H]
	\caption{APG \label{alg:APG_supp} }
	\begin{algorithmic}[1]
		\STATE Input: nonnegative $N$-way tensor $\mathcal{T}$
		\STATE Output: nonnegative factors $A^{(1)}, A^{(2)},\ldots, A^{(N)}$.
		\STATE Initialization: Choose $\delta_w<1$, $t_0=1$, and a set of initial factor matrices $\big(A^{(1)}_{0}, \ldots, A^{(N)}_{0}\big )$. Set $k=1$.
		\REPEAT
		\FOR{i=1,\ldots,N}
		\STATE Compute  $t_k=\frac12\big( 1 + \sqrt{1 + 4 t_{k-1}^2}\big)$, $\hat w_{k-1}=\frac{t_{k-1}-1}{t_k}$ and \vspace*{-2mm}
		$$w_{k-1}^{(i)}=\min \left( \hat w_{k-1}, \delta_w \sqrt{\frac{L_{k-2}^{(i)}}{L_{k-1}^{(i)}}}\right). $$ \vspace*{-2mm}
		\STATE Compute an extrapolation point \vspace*{-2mm}
		$$\hat  A_{k-1}^{(i)}= A_{k-1}^{(i)} + w_{k-1}^{(i)} \Big( A_{k-1}^{(i)}-A_{k-2}^{(i)} \Big).$$  \vspace*{-2mm}
		\STATE Update $A_{k}^{(i)}$ by projected gradient step: \vspace*{-2mm}
		\begin{equation}
		\label{projgrad}
		A_{k}^{(i)}=\max \left(0,\hat  A_{k-1}^{(i)}- \frac{1}{L_{k-1}^{(i)}}  \Big (\hat A_{k-1}^{(i)} \big(B_{k-1}^{(i)}\big)^T- \mathcal T_{[i]}\Big) B_{k-1}^{(i)}\right) \vspace*{-2mm}
		\end{equation}
		\ENDFOR
		\IF{$F\big( A_k\big)>F\big( A_{k-1}\big) $}
		\STATE Update $A_{k}^{(i)}$ by the projected gradient step \eqref{projgrad} with $\hat  A_{k-1}^{(i)} = A_{k-1}^{(i)}$.
		\ENDIF
		\STATE Set $k=k+1$.
		\UNTIL{some criteria is satisfied}
	\end{algorithmic}
\end{algorithm}
\end{minipage}

%****************************************************
\begin{minipage}{\linewidth}
	%\algsetup{indent=2em}
\begin{algorithm}[H]
	\caption{iBPG \label{alg:iBPG_supp} }
	\begin{algorithmic}[1]
		\STATE Input: a nonnegative $N$-way tensor $\mathcal{T}$
		\STATE Output: nonnegative factors $A^{(1)}, A^{(2)},\ldots, A^{(N)}$.
		\STATE Initialization: Choose $\delta_w=0.99$, $\beta=1.01 $, $t_0=1$, and 2 sets of initial factor matrices $\big(A^{(1)}_{-1}, \ldots, A^{(N)}_{-1}\big )$ and $\big(A^{(1)}_{0}, \ldots, A^{(N)}_{0}\big )$. Set $k=1$.
		\STATE Set $A_{\rm{prev}}^{(i)}= A^{(i)}_{-1}$, $i=1,\ldots,N$. \% \textit{$A_{\rm{prev}}^{(i)}$ is to save the previous value of block $i$. }
		\STATE Set $A_{\rm{cur}}^{(i)}= A^{(i)}_{0}$, $i=1,\ldots,N$. \% \textit{$A_{\rm{cur}}^{(i)}$ is to save the current value of block $i$. }
		\REPEAT
		\FOR{i=1,\ldots,N}
		\STATE Compute  $t_k=\frac12\big( 1 + \sqrt{1 + 4 t_{k-1}^2}\big)$, $\hat w_{k-1}=\frac{t_{k-1}-1}{t_k}$ and  \vspace*{-2mm}
		$$w_{k-1}^{(i)}=\min \left( \hat w_{k-1}, \delta_w \sqrt{\frac{L_{k-2}^{(i)}}{L_{k-1}^{(i)}}}\right). $$  \vspace*{-3mm}
		\REPEAT
		\STATE Compute two extrapolation points  \vspace*{-2mm}
		$$\hat  A^{(i,1)}= A_{\rm{cur}}^{(i)} + w_{k-1}^{(i)} \Big( A_{\rm{cur}}^{(i)}-A_{\rm{prev}}^{(i)} \Big),$$ \vspace*{-2mm}
		and
		$$\hat A^{(i,2)}= A_{\rm{cur}}^{(i)} +\beta w_{k-1}^{(i)} \Big( A_{\rm{cur}}^{(i)}-A_{\rm{prev}}^{(i)} \Big)$$ \vspace*{-2mm}
		\STATE Set $A_{\rm{prev}}^{(i)}=A_{\rm{cur}}^{(i)}$.
		\STATE Update $A_{\rm{cur}}^{(i)}$ by projected gradient step:
		\begin{equation*}
		%\label{projgrad2}
		A_{\rm{cur}}^{(i)}=\max \left(0,\hat  A^{(i,2)}- \frac{1}{L_{k-1}^{(i)}}  \Big (\hat A^{(i,1)} \big(B_{k-1}^{(i)}\big)^T- \mathcal T_{[i]}\Big) B_{k-1}^{(i)}\right). \vspace*{-2mm}
		\end{equation*}
		\UNTIL{some criteria is satisfied}
		\STATE Set $A_{k}^{(i)}=A_{\rm{cur}}^{(i)}$.
		\ENDFOR
		\STATE Set $k=k+1$.
		\UNTIL{some criteria is satisfied}
	\end{algorithmic}
\end{algorithm}
\end{minipage}
%\clearpage

\section{Full experimental results}\label{appendix:allexps}
In this section, we provide more plots for the experimental set up presented in \S\,\ref{sec:exp}\,.
This includes other algorithms, and other dimensions of the input tensors.

\begin{center}\begin{figure}[ht!] % *** Balanced
\centering\includegraphics[width=\linewidth]{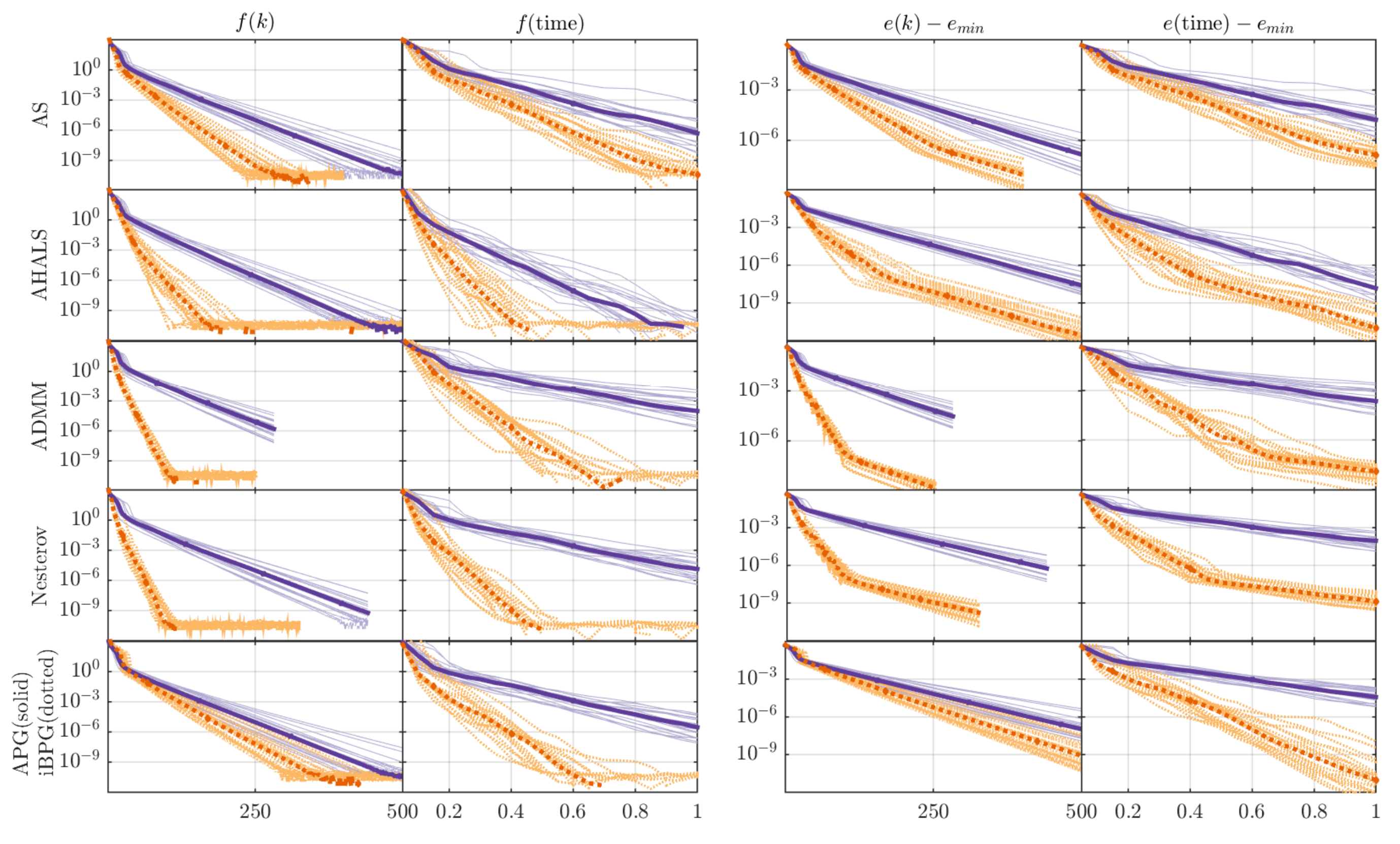}
\caption{Experiments on synthetic data with $[I_1,I_2,I_3,r] = [50,50,50,10]$ on algorithms without HER (solid purple curves) and with HER (dotted orange curves).
		Thick curves are median of the 20 thin curve with the same color.
		The results show that HER improves BCD algorithms, and that HER-BCD has a better performance than APG and iBPG.
		}\label{fig:exp:balance_supp}
\end{figure}\end{center}
\begin{center}\begin{figure}[ht!] % *** Unbalanced
\centering\includegraphics[width=\linewidth]{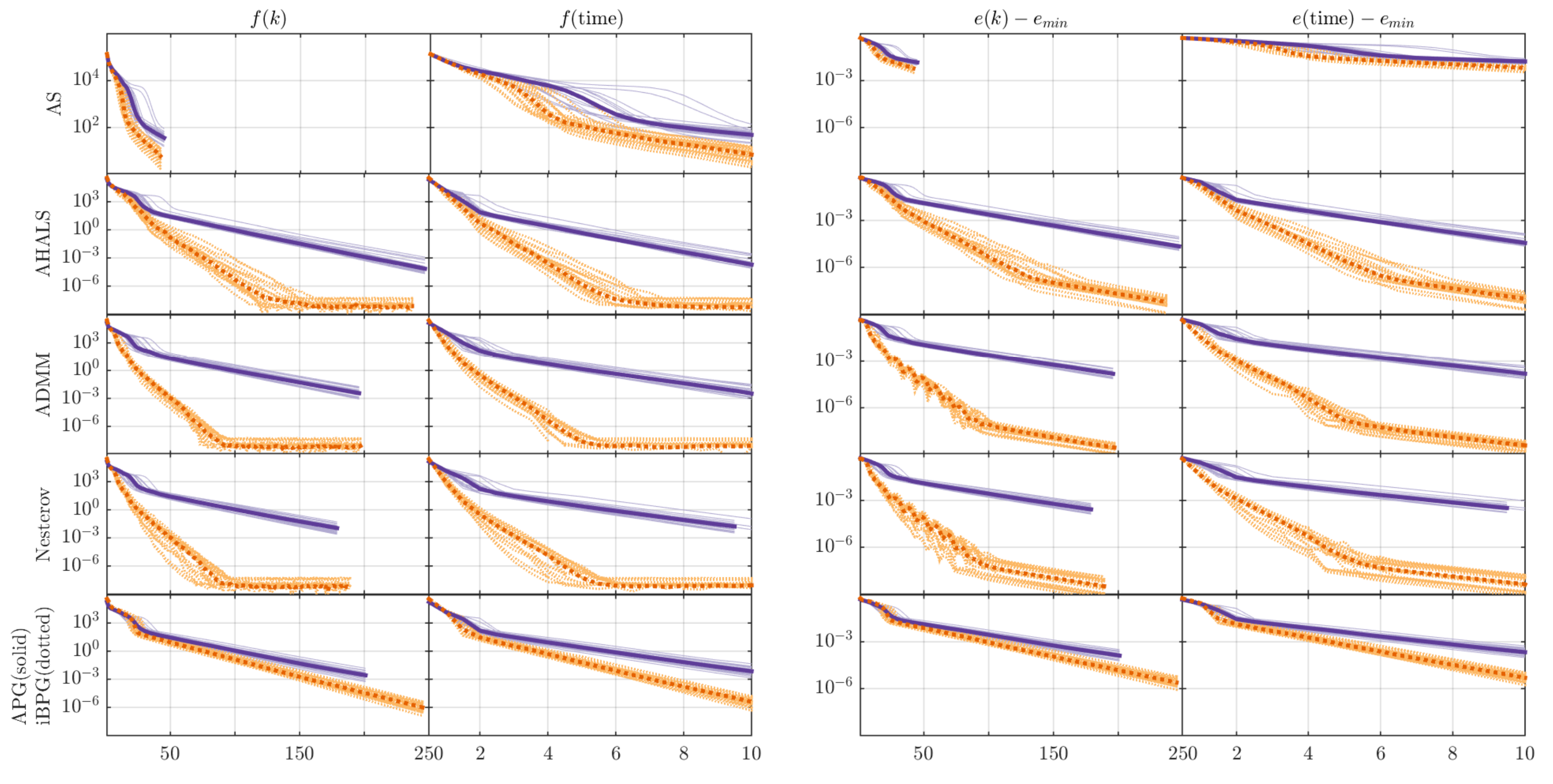}
\caption{Experiments on synthetic data with unbalanced tensors with $[I_1,I_2,I_3,r] = [150,10^3,10,10]$.
The same setting as in Figure \ref{fig:exp:balance_supp} is used.
The results show that HER improves BCD algorithms, and that HER-BCD has a better performance than APG and iBPG.
Note that here AS only ran approximately 40 iterations in 10 seconds due to high per-iteration cost.
}\label{fig:exp:unbalance_supp}
\end{figure}\end{center}
% * * *  * * *  * * *  * * *  * * *
% unbalanced high rank noisless
% * * *  * * *  * * *  * * *  * * *
\begin{center}\begin{figure}[ht!]
\centering\includegraphics[width=\linewidth]{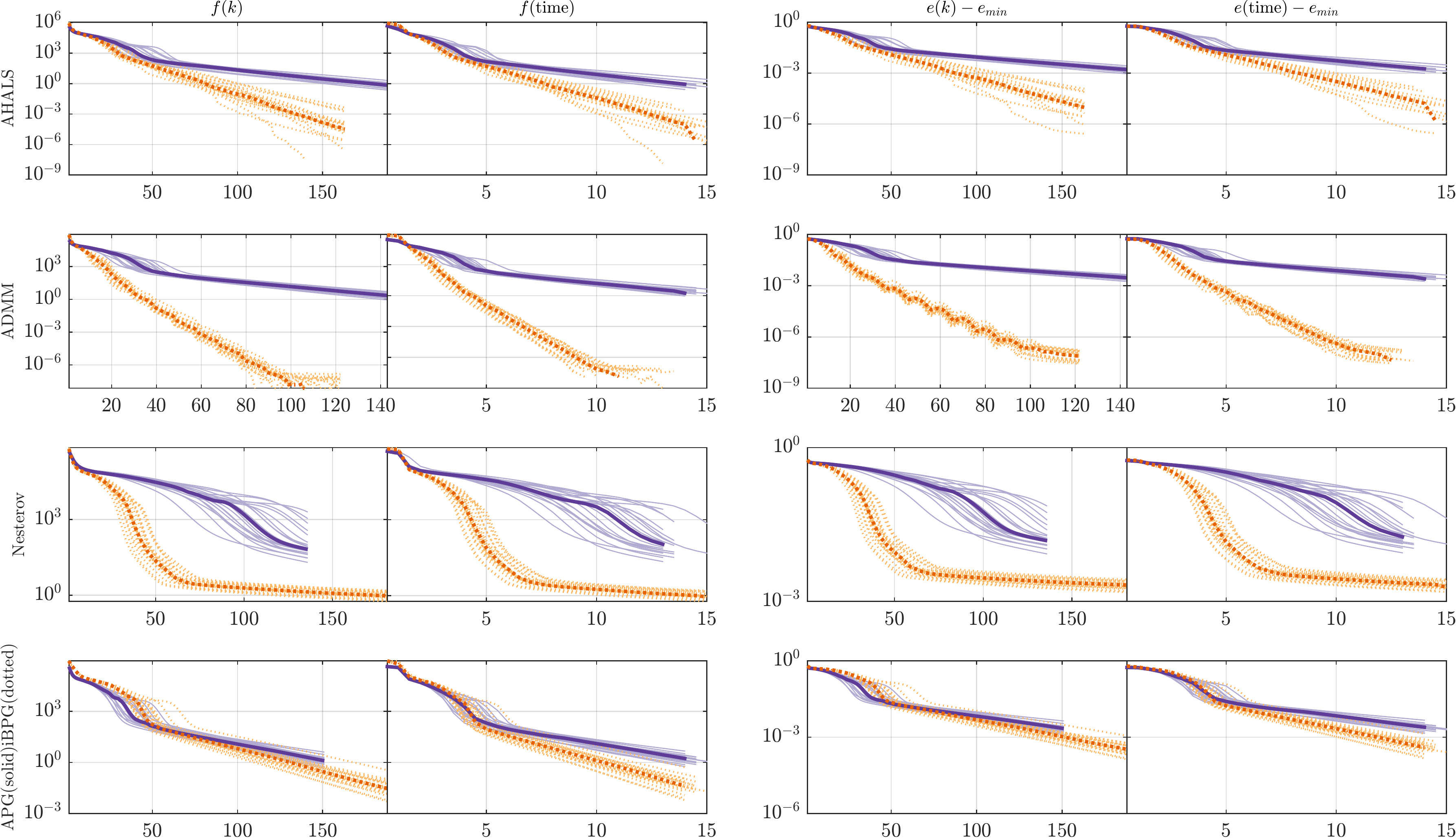}
\caption{Experiments on synthetic data with  $[I_1,I_2,I_3,r]=[150,10^3,50,25]$.
The same setting as in Figure~\ref{fig:exp:balance_supp} is used.
The results show that HER improves BCD algorithms, and that  HER-BCD has a better performance than APG and iBPG.
}\label{fig:exp:highrank_supp}
\end{figure}\end{center}
% * * *  * * *  * * *  * * *  * * *
% big noisy
% * * *  * * *  * * *  * * *  * * *
\begin{center}\begin{figure}[ht!] % *** Big Noisy
\centering\includegraphics[width=\linewidth]{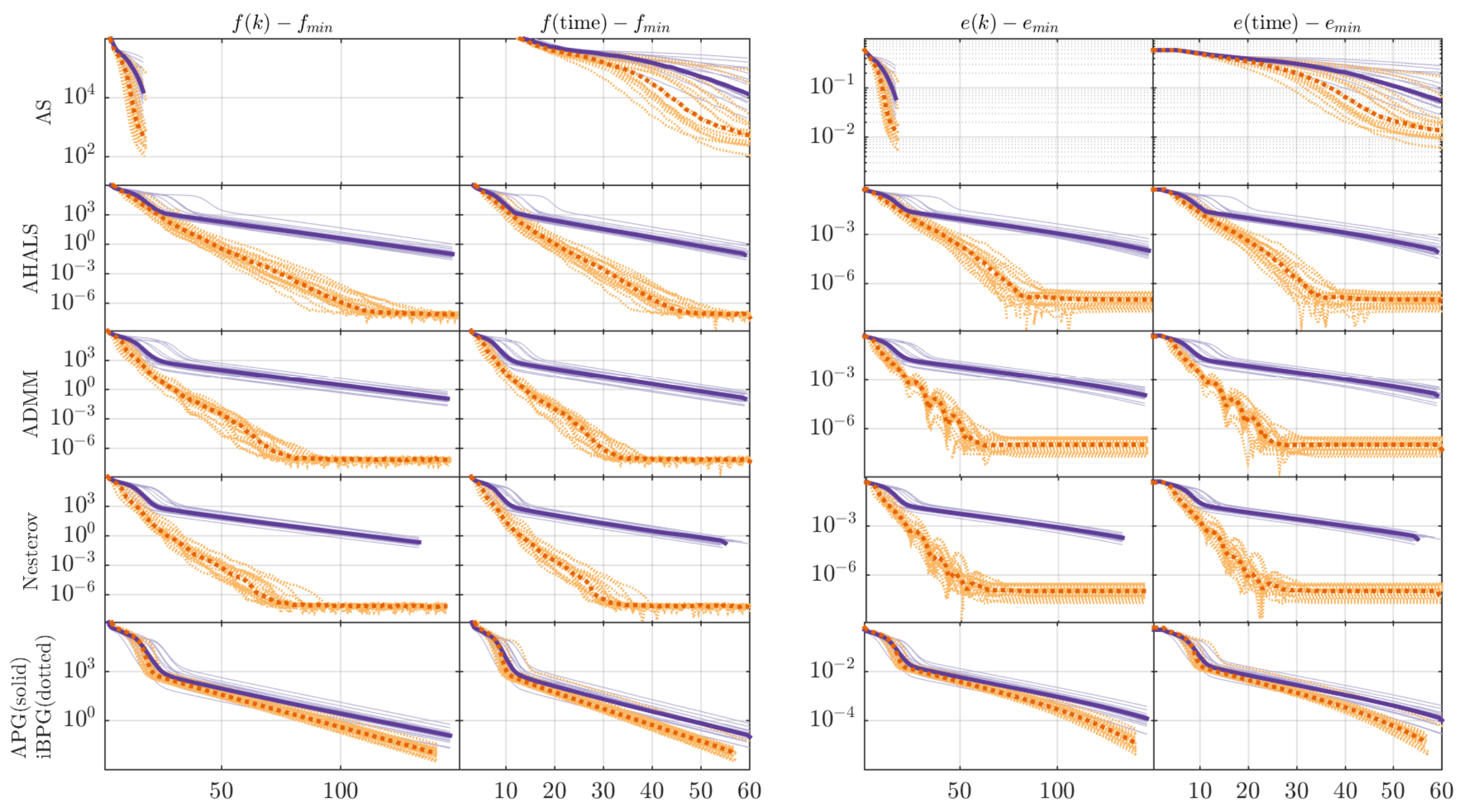}
\caption{Experiments on synthetic data with $[I_1,I_2,I_3,r,\sigma] = [500,500,500,10,0.01]$ on algorithms without HER (solid purple curves) and with HER (dotted orange curves).
The same setting as in Figure \ref{fig:exp:balance_supp} is used.
The results show that HER improves BCD algorithms, and that HER-BCD has a better performance than APG and iBPG. Here AS only ran approximately 40 iterations in 10 seconds due to high per-iteration cost.
}\label{fig:exp:unbalance_supp}
\end{figure}\end{center}
% * * *  * * *  * * *  * * *  * * *
% Ill condition
% * * *  * * *  * * *  * * *  * * *
\begin{center}\begin{figure}[ht!] % *** Big Noisy
\centering\includegraphics[width=\linewidth]{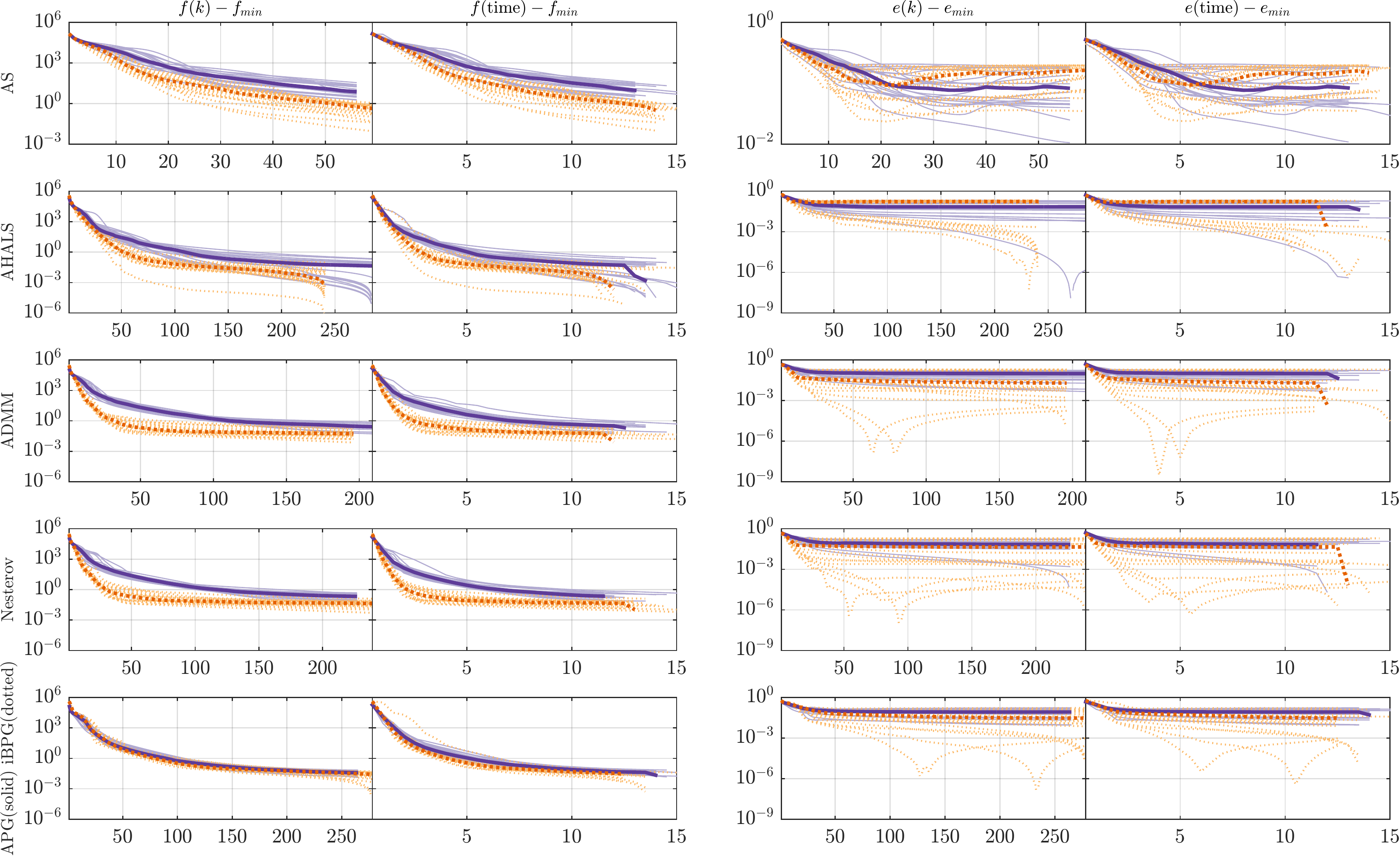}
\caption{Experiments on synthetic data with $[I_1,I_2,I_3,r,\sigma] = [150,10^3,50,12,0.01]$ and  with ill-conditioned tensors; ground truth mode factors $A_i(:,1) = 0.99 A_i(:,2) + 0.01 A_i(:,1)$ for $i\in \{1,2,3\}$.
In terms of best case performance, the results show that HER improves BCD algorithms, and that HER-BCD has a better performance than APG and iBPG on $f$, and a better performance than APG on $e$.
}\label{fig:exp:unbalance_supp}
\end{figure}\end{center}
\begin{center}\begin{figure}[ht!] % *** Big Noisy
\centering\includegraphics[width=\linewidth]{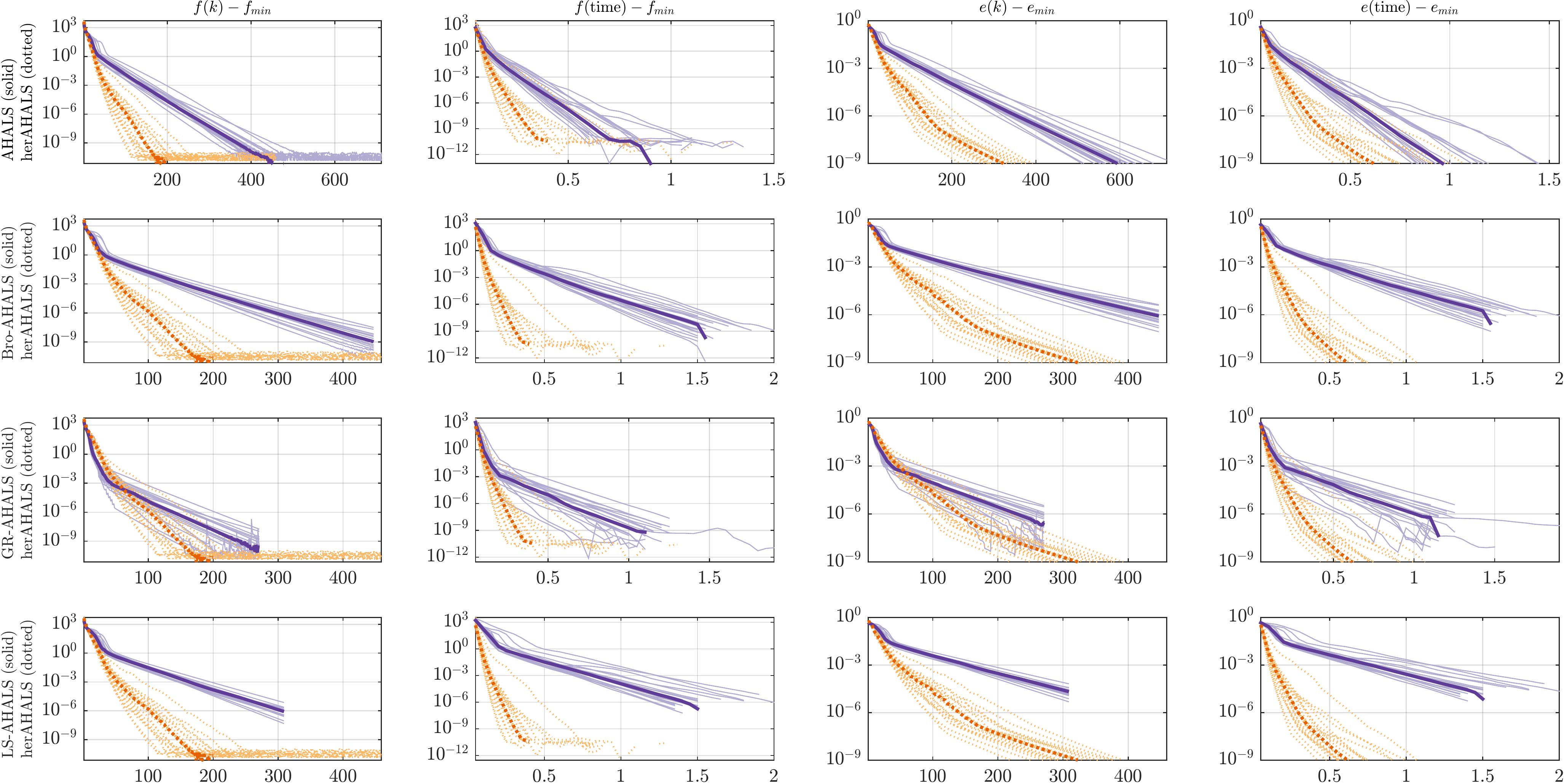}
\caption{Experiments on synthetic data with
		$[I_1,I_2,I_3,r] = [50,50,50,10]$ for different AHALS algorithms.
		Results show HER-AHALS has better performance than all other AHALS algorithms.
}\label{fig:exp:ahals_BH50_50_50_10_n0}
\end{figure}\end{center}
\begin{center}\begin{figure}[ht!] % *** Big Noisy
	\centering\includegraphics[width=\linewidth]{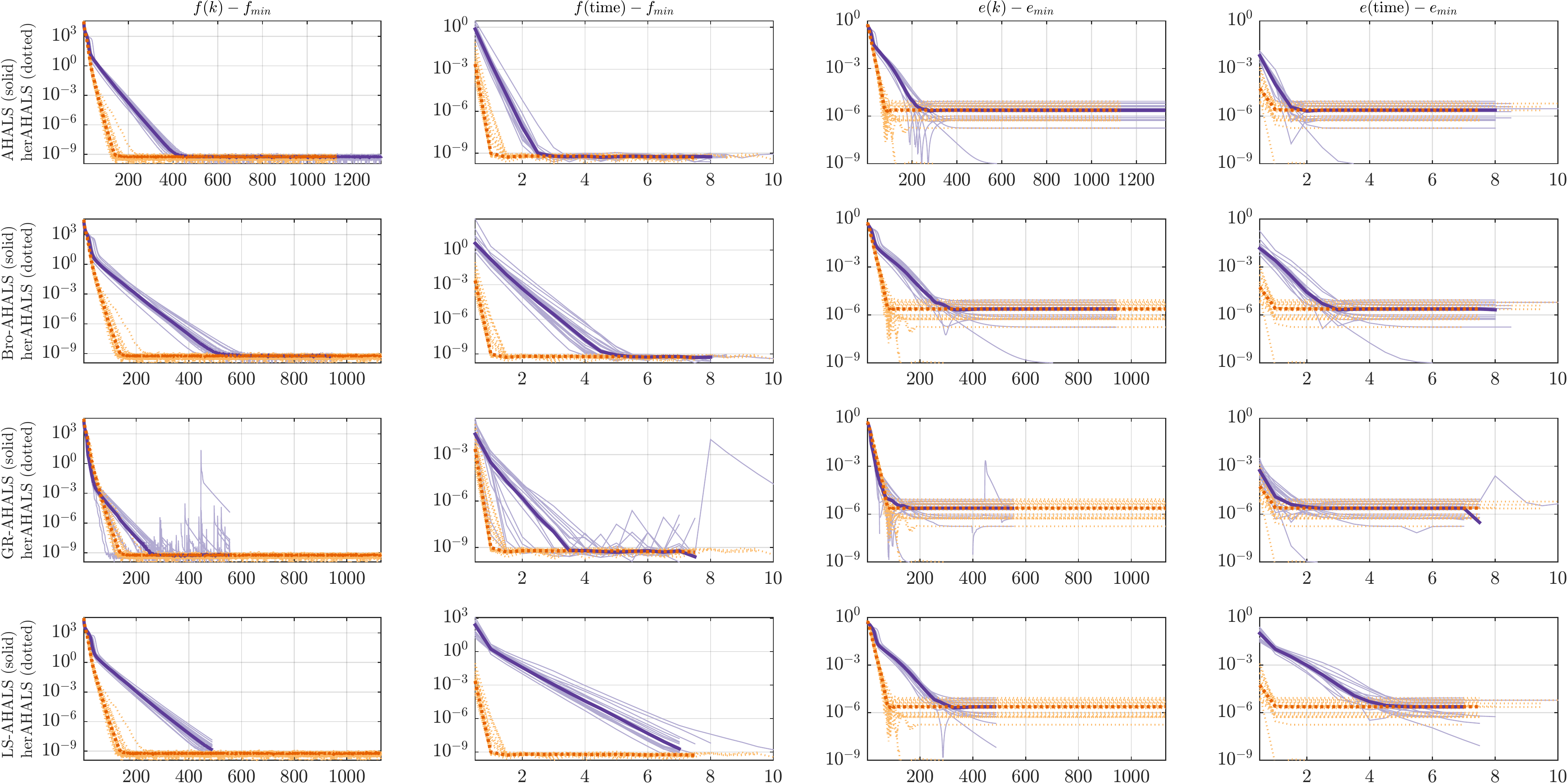}
	\caption{Experiments on synthetic data with
			$[I_1,I_2,I_3,r,\sigma] = [150,1000,50,12,0.01]$ for different AHALS algorithms.
			Results show HER-AHALS has better  performance than all other AHALS algorithms.
			}\label{fig:exp:ahals_BH150_1000_50_12_n001}
	\end{figure}\end{center}
\begin{center}\begin{figure}[ht!] % *** Big Noisy
\centering\includegraphics[width=\linewidth]{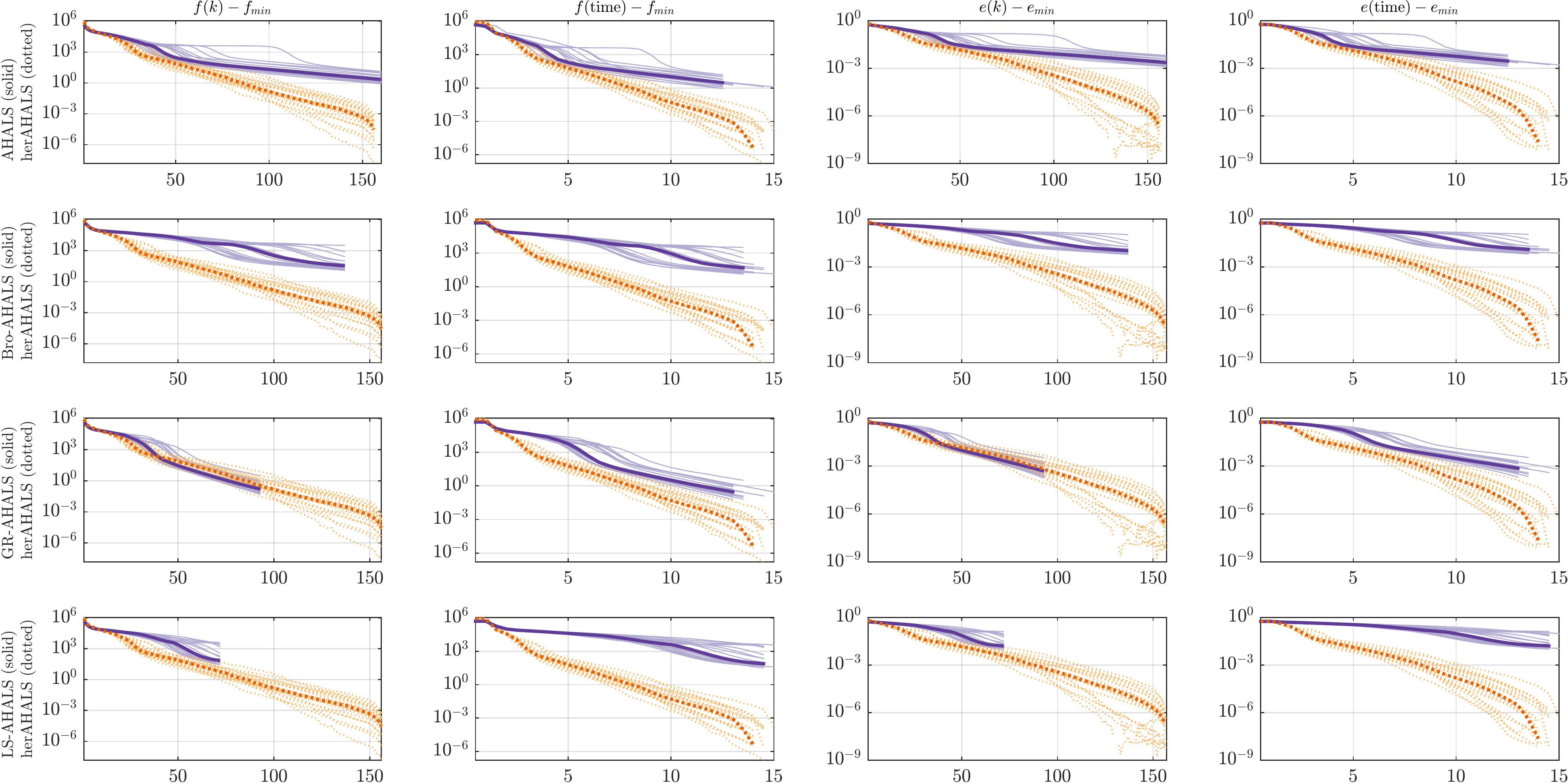}
\caption{Experiments on synthetic data with
     	$[I_1,I_2,I_3,r,\sigma] = [150,1000,50,25,0.01]$ for different AHALS algorithms.
     	Results show HER-AHALS has better performance than all other AHALS algorithms.
		}\label{fig:exp:ahals_BH150_1000_50_25_n001}
\end{figure}\end{center}
\begin{center}\begin{figure}[ht!] % *** Big Noisy
		\centering\includegraphics[width=\linewidth]{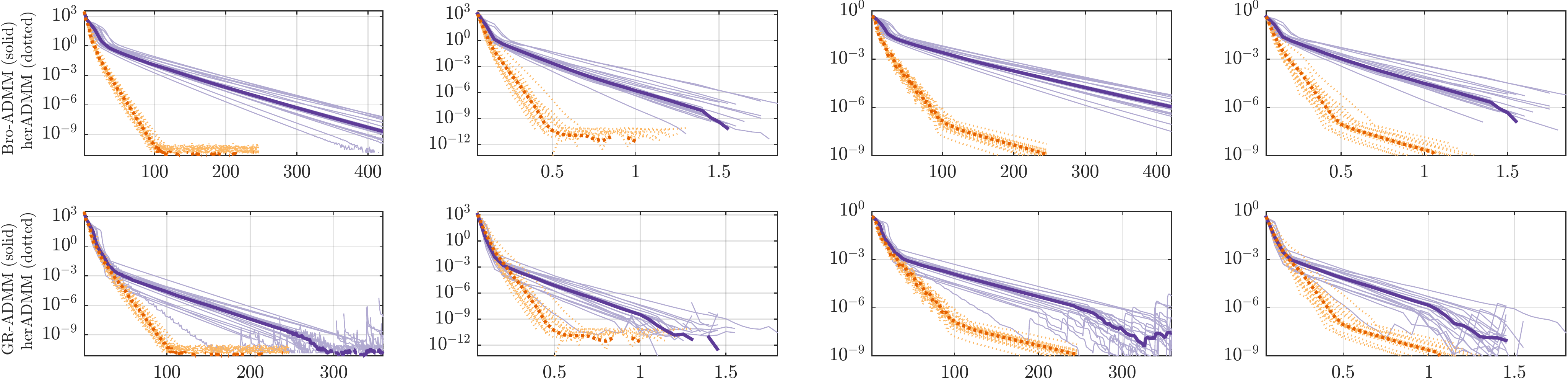}
		\caption{Experiments on synthetic data  with
		$[I_1,I_2,I_3,r] = [50,50,50,10]$ over different AO-ADMM algorithms.
		Results show HER-ADMM has better performance than all other AO-ADMM algorithms.
		}\label{fig:exp:admm_BH50_50_50_10_n0}
\end{figure}\end{center}
\begin{center}\begin{figure}[ht!] % *** Big Noisy
			\centering\includegraphics[width=\linewidth]{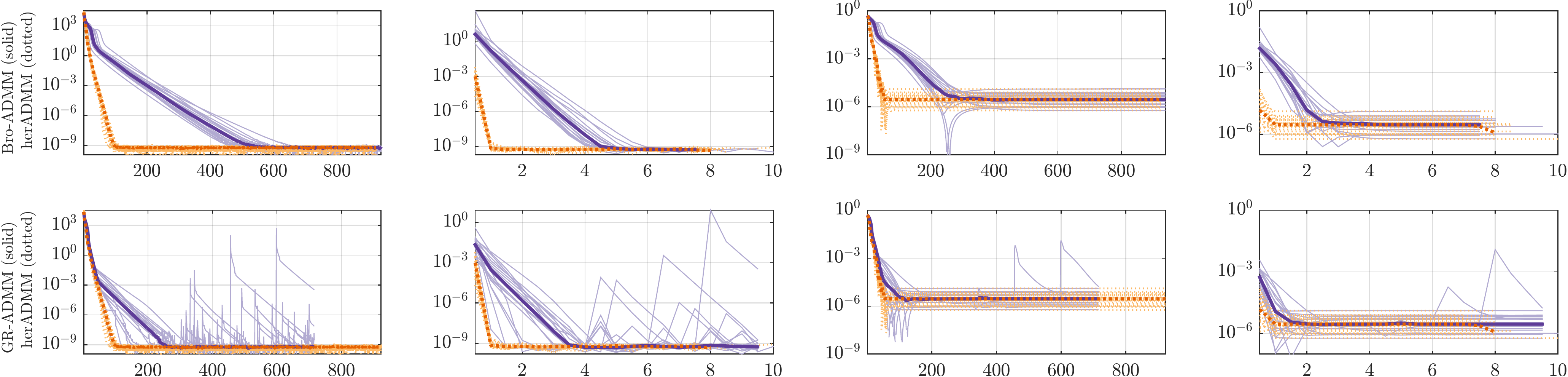}
			\caption{Experiments on synthetic data  with 				$[I_1,I_2,I_3,r,\sigma] = [150,1000,50,12,0.01]$ for different AO-ADMM algorithms.
				Results show HER-ADMM has better  performance than all other AO-ADMM algorithms.
			}\label{fig:exp:admm_BH150_1000_50_12_n001}
	\end{figure}\end{center}
\begin{center}\begin{figure}[ht!] % *** Big Noisy
		\centering\includegraphics[width=\linewidth]{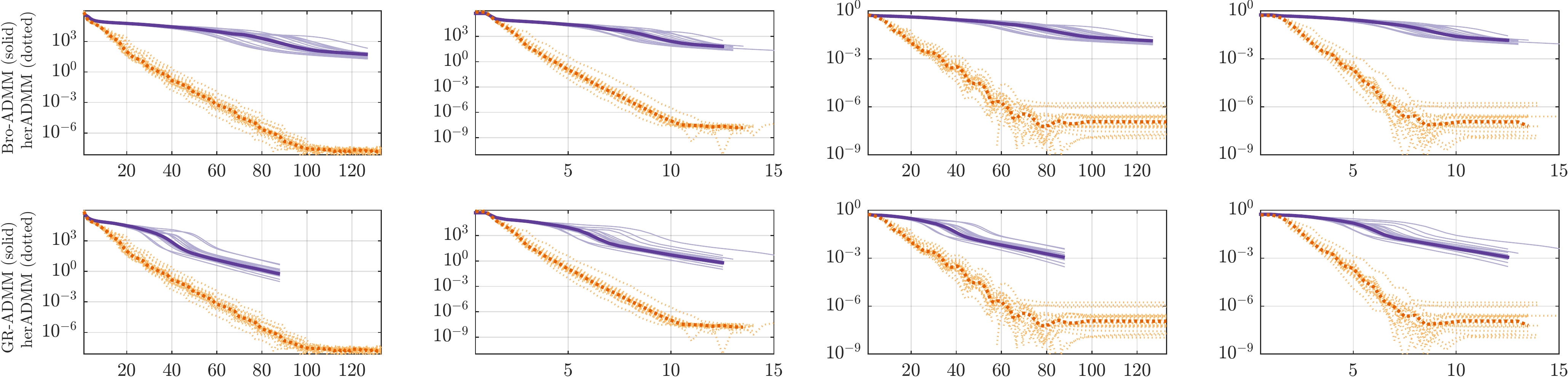}
		\caption{Experiments on synthetic data  with
			$[I_1,I_2,I_3,r,\sigma] = [150,1000,50,25,0.01]$ for different AO-ADMM algorithms.
			Results show HER-ADMM has better  performance than all other AO-ADMM algorithms.
		}\label{fig:exp:admm_BH150_1000_50_25_n001}
\end{figure}\end{center}
\begin{center}\begin{figure}[ht!] % *** Big Noisy
		\centering\includegraphics[width=\linewidth]{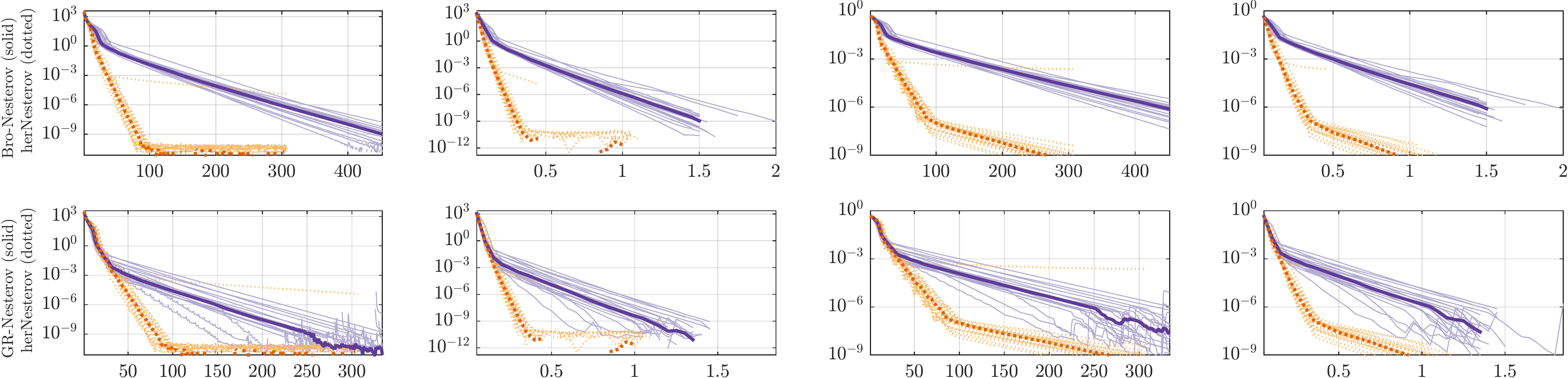}
		\caption{Experiments on synthetic data with
			$[I_1,I_2,I_3,r] = [50,50,50,10]$ over different AO-Nesterov algorithms.
			Results show HER-Nesterov has better convergence performance than all other AO-Nesterov algorithms.
		}\label{fig:exp:nest_BH50_50_50_10_n0}
\end{figure}\end{center}
\begin{center}\begin{figure}[ht!] % *** Big Noisy
			\centering\includegraphics[width=\linewidth]{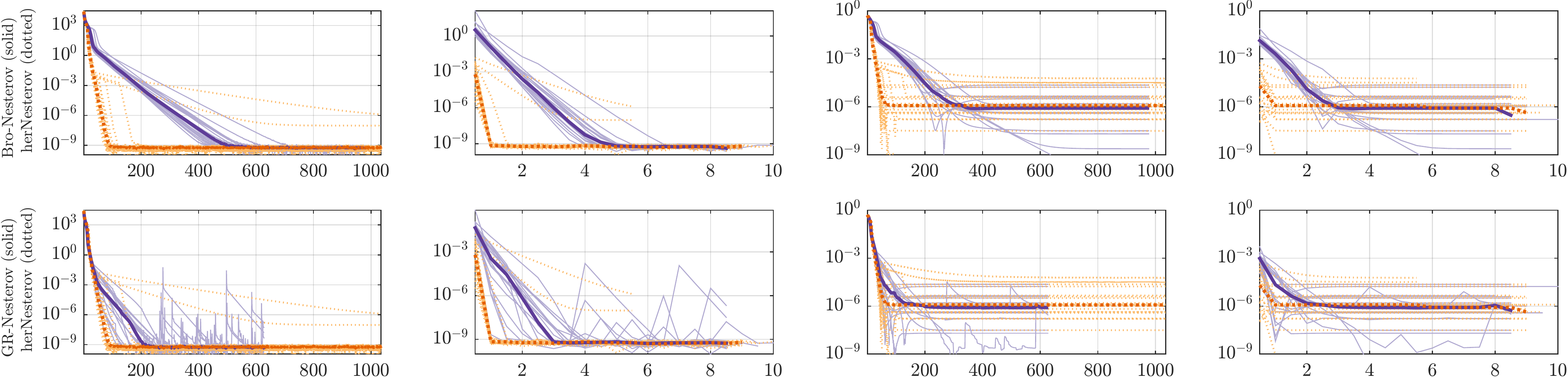}
			\caption{Experiments on synthetic data  with
				$[I_1,I_2,I_3,r,\sigma] = [150,1000,50,12,0.01]$ for different AO-Nesterov algorithms.
				Results show HER-Nesterov has better  performance than all other AO-Nesterov algorithms.
			}\label{fig:exp:nest_BH150_1000_50_12_n001}
	\end{figure}\end{center}
\begin{center}\begin{figure}[ht!] % *** Big Noisy
		\centering\includegraphics[width=\linewidth]{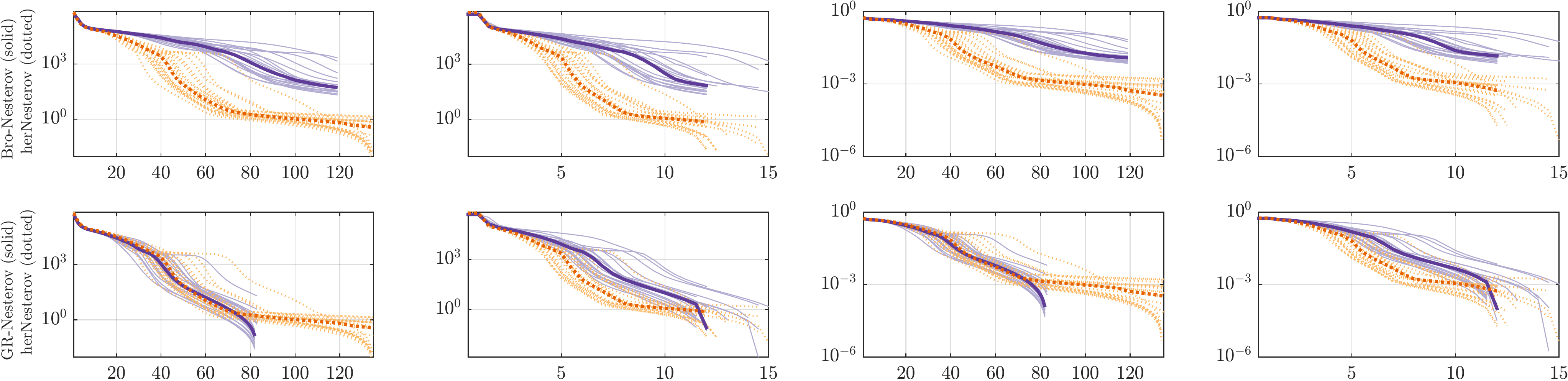}
		\caption{Experiments on synthetic data  with
			$[I_1,I_2,I_3,r,\sigma] = [150,1000,50,25,0.01]$ for different AO-Nesterov algorithms.
			Results show HER-Nesterov has better  performance than all other AO-Nesterov algorithms,  except GR-Nesterov.
		}\label{fig:exp:nest_BH150_1000_50_25_n001}
\end{figure}\end{center}

% * * *  * * *  * * *  * * *  * * *
% HER VS BRO VS GR VS LS
% * * *  * * *  * * *  * * *  * * *
\clearpage
We now shows the full results on comparing HER-BCD with Bro-BCD, GR-BCD (using Gradient ratio, in the modified form) and LS-BCD (line search, in the modified form), with BCD being AHALS, AO ADMM or AO Nesterov.
While Fig.\ref{fig:exp:ahals_BH50_50_50_10_n0} to Fig.\ref{fig:exp:nest_BH150_1000_50_25_n001} will show the full result, the table below shows the summary of theses results: here each curve represents the average over 20 trials.
In general, HER-BCD has the best performance among all the extrapolated BCD methods.

\begin{table}[ht!]
\label{tbe:BH_summary}
\begin{tabular}{ccc}
 $[50,50,50,10,0]$ & $[150,10^3,50,12,0.01]$ & $[150,10^3,50,25,0.01]$  \\ \hline
\includegraphics[width=0.32\linewidth]{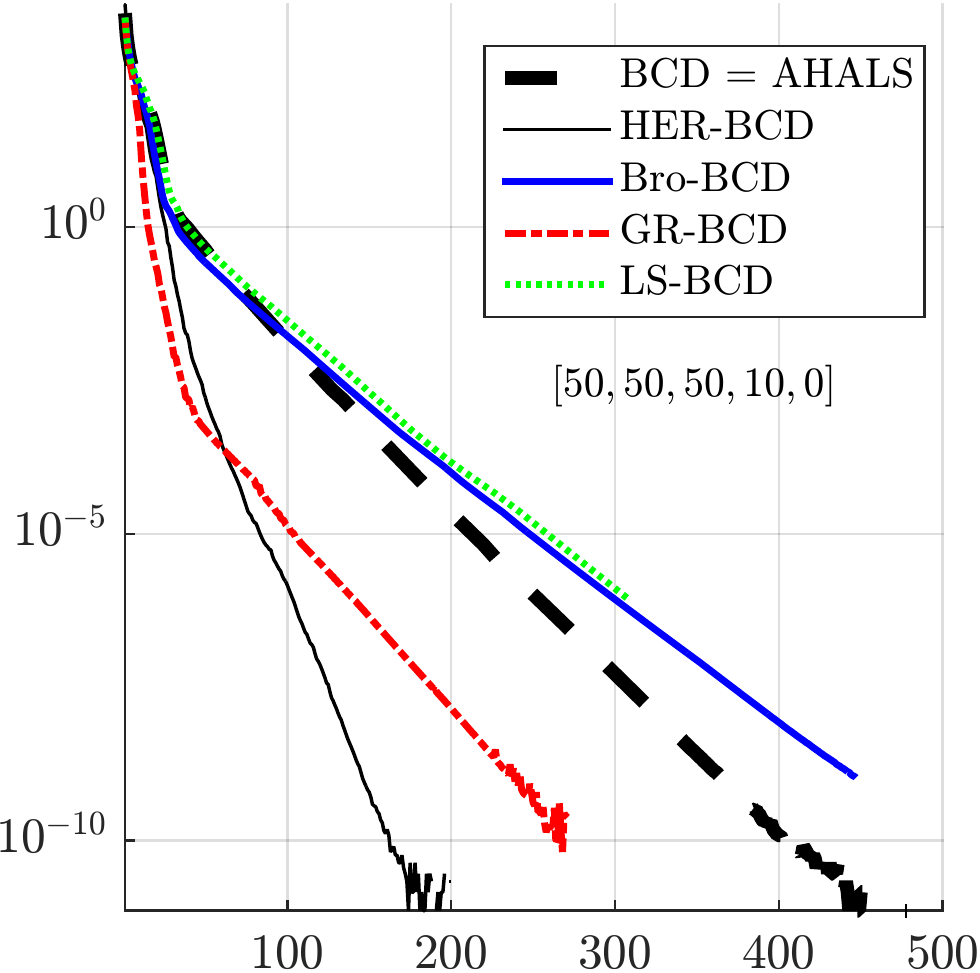}  
& \includegraphics[width=0.32\linewidth]{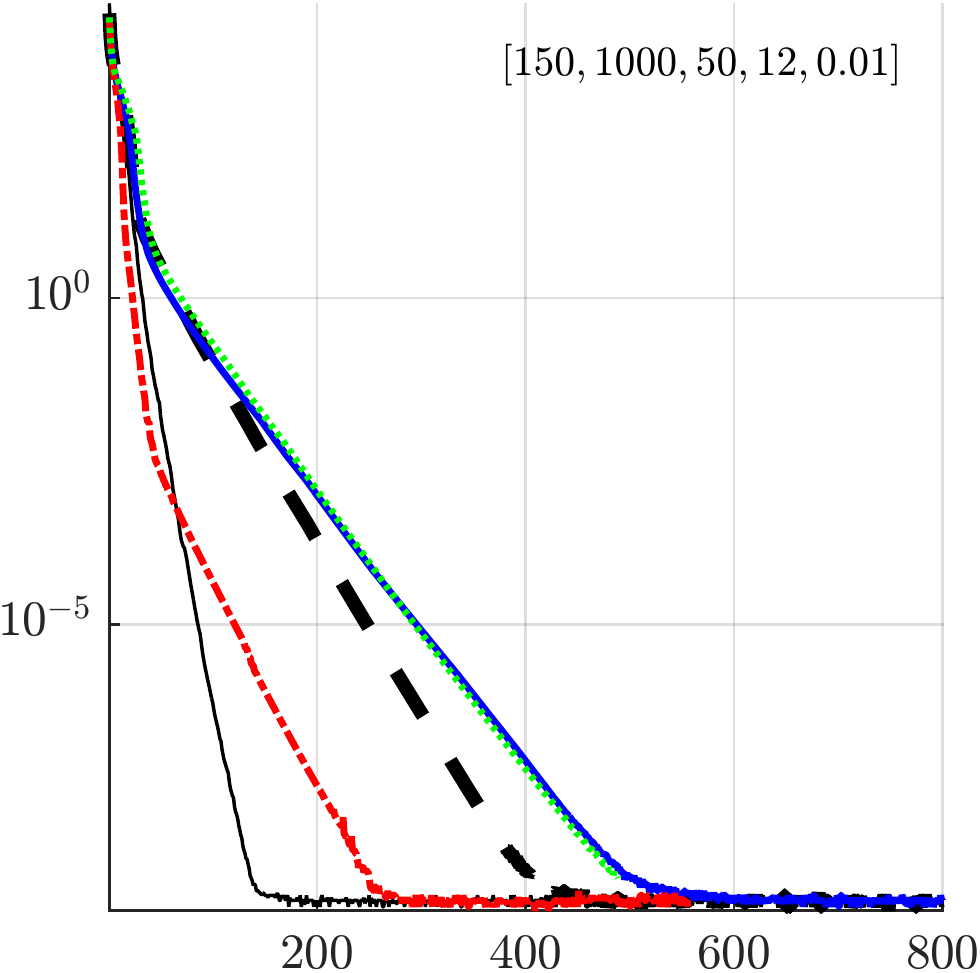}  
& \includegraphics[width=0.32\linewidth]{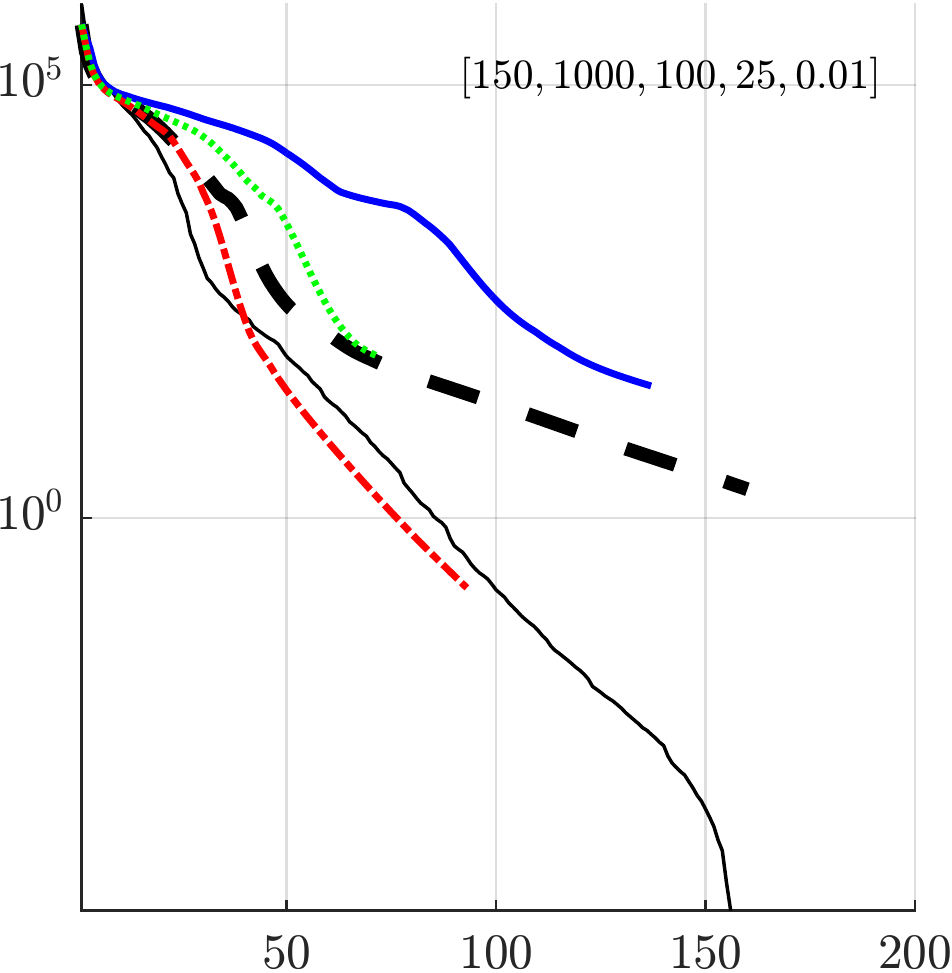}   
\\~\\
\includegraphics[width=0.32\linewidth]{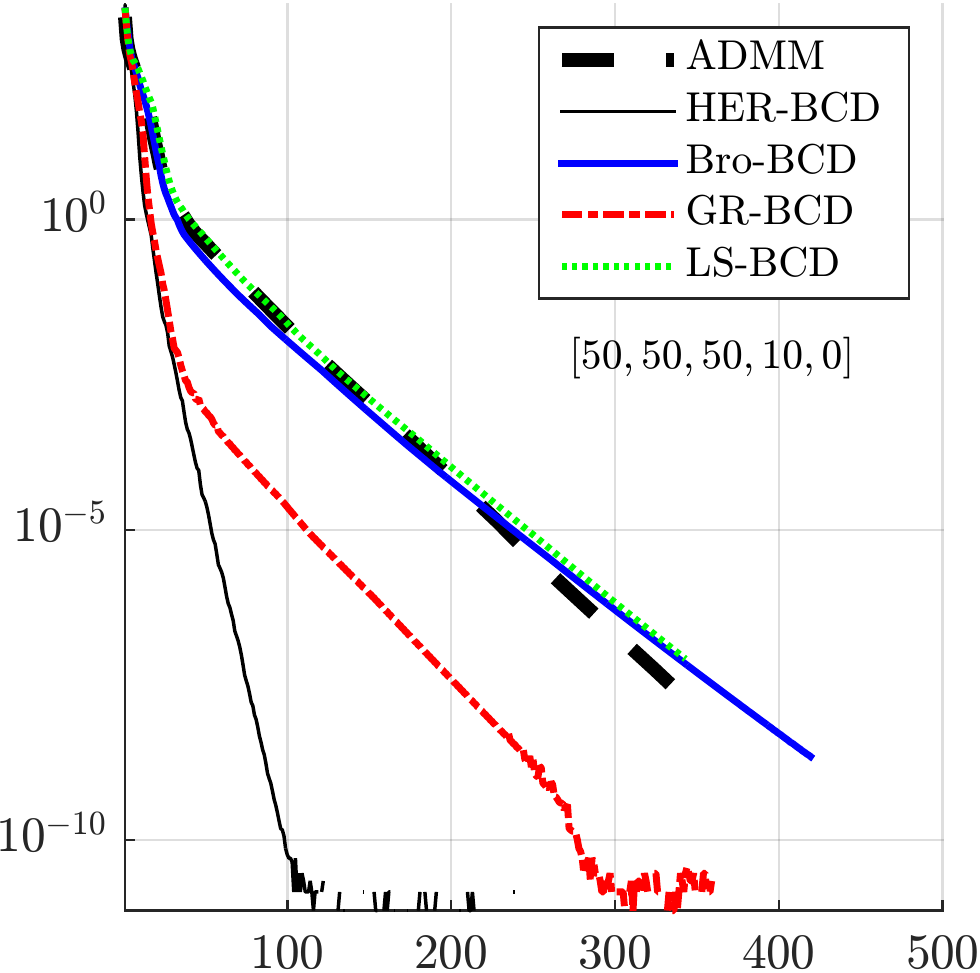}  
& \includegraphics[width=0.32\linewidth]{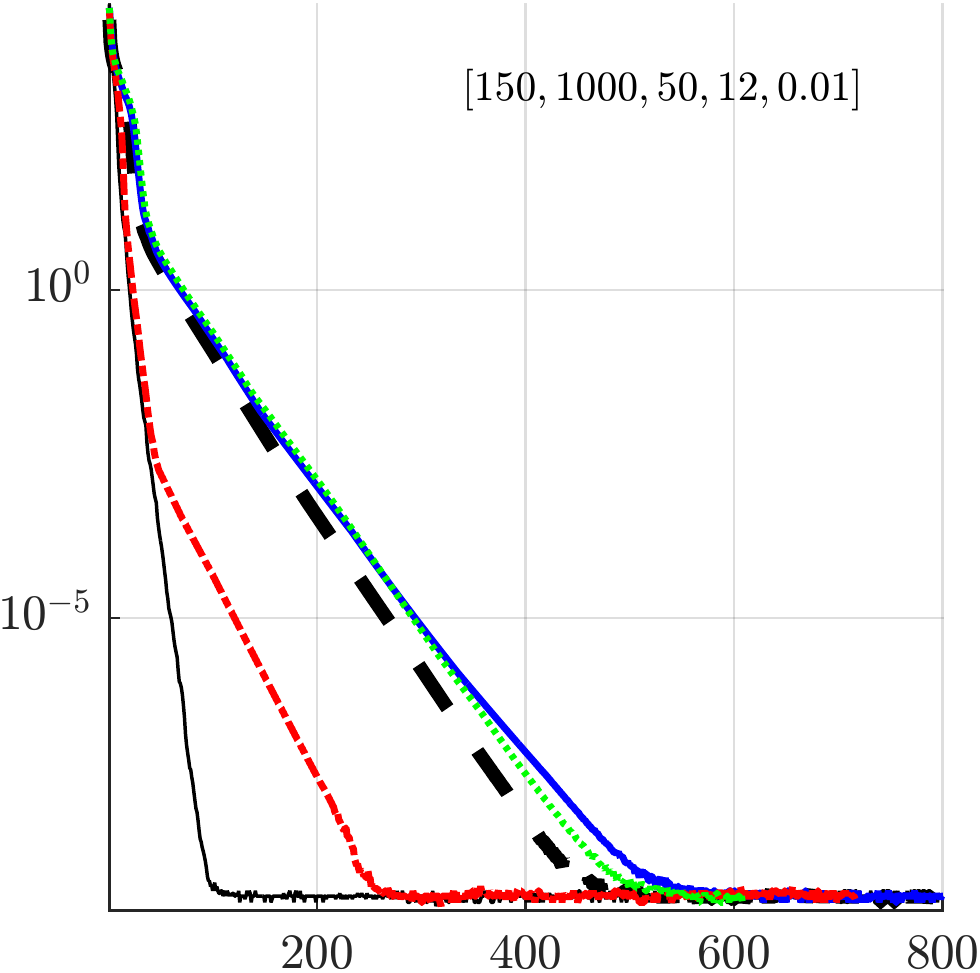}  
& \includegraphics[width=0.32\linewidth]{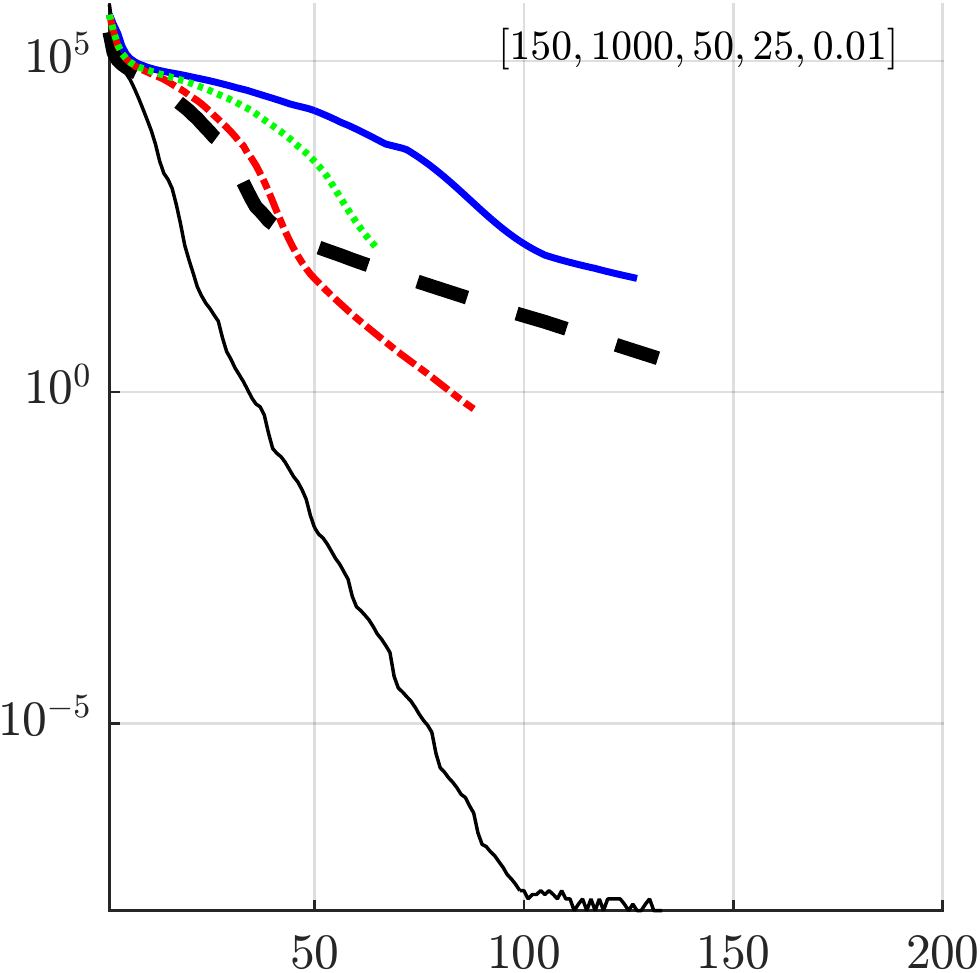}  
\\~\\
\includegraphics[width=0.32\linewidth]{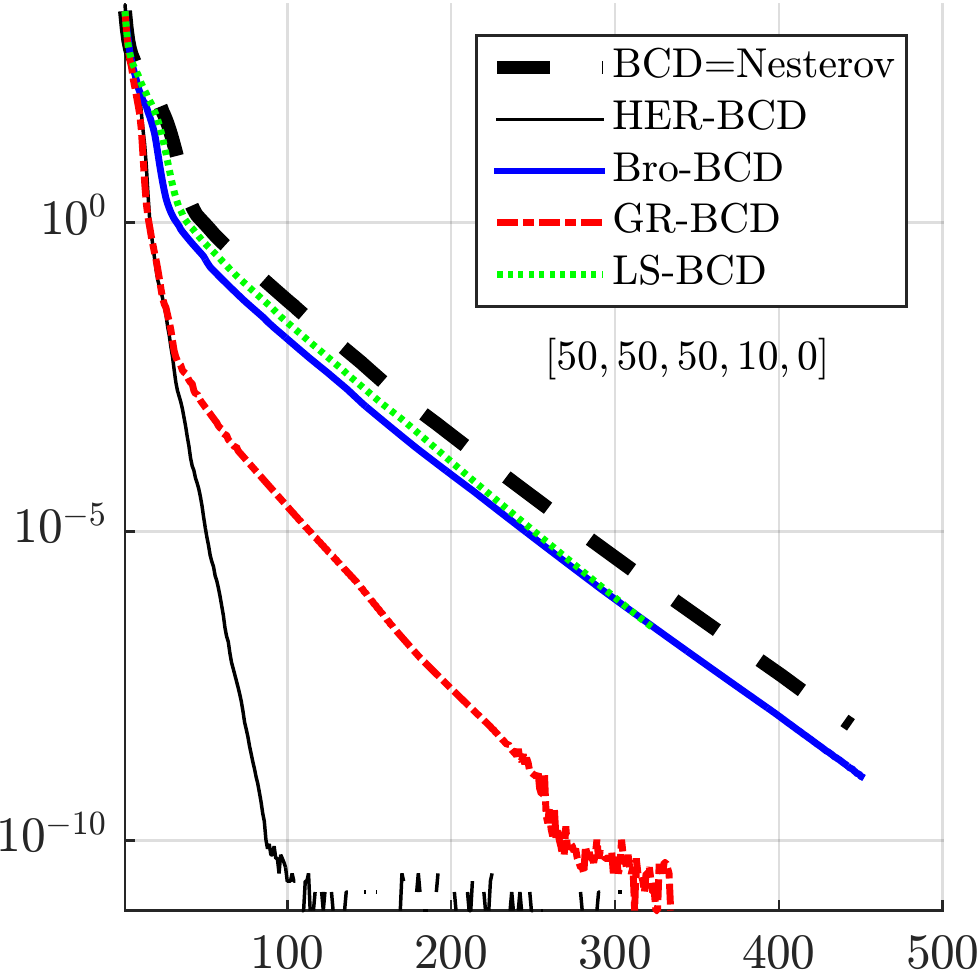}  
& \includegraphics[width=0.32\linewidth]{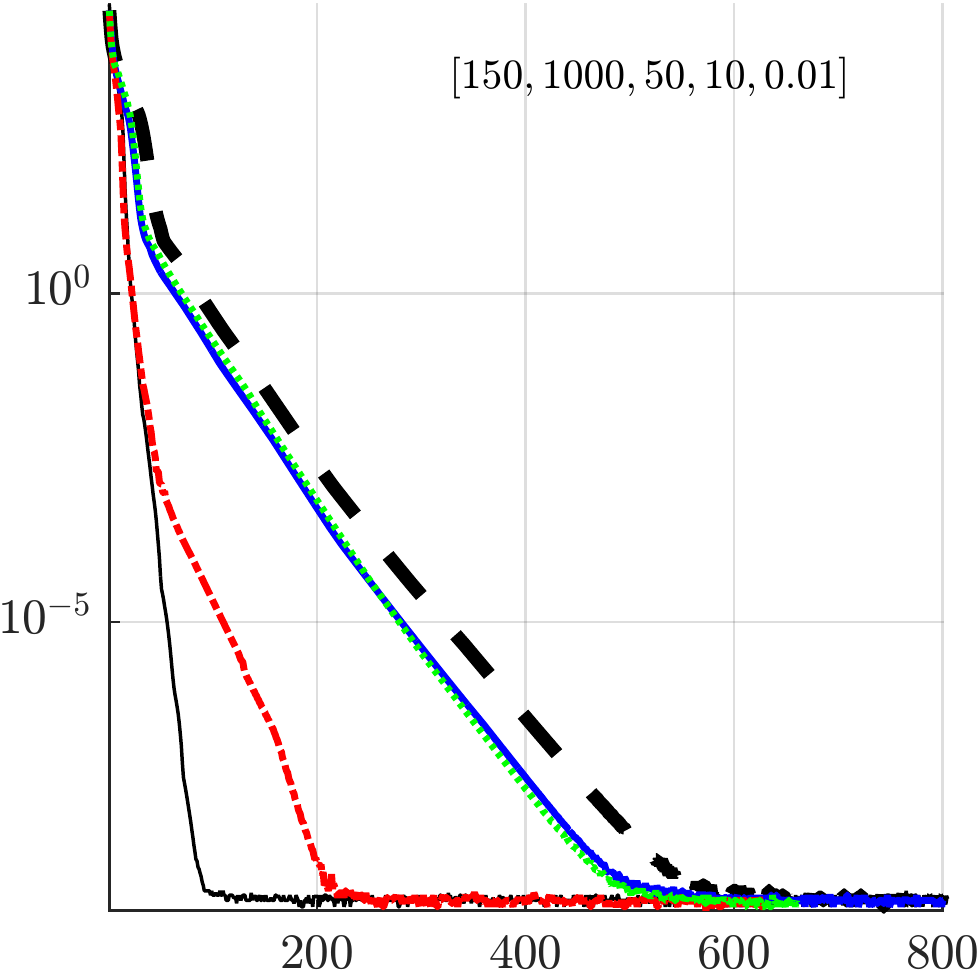}  
& \includegraphics[width=0.32\linewidth]{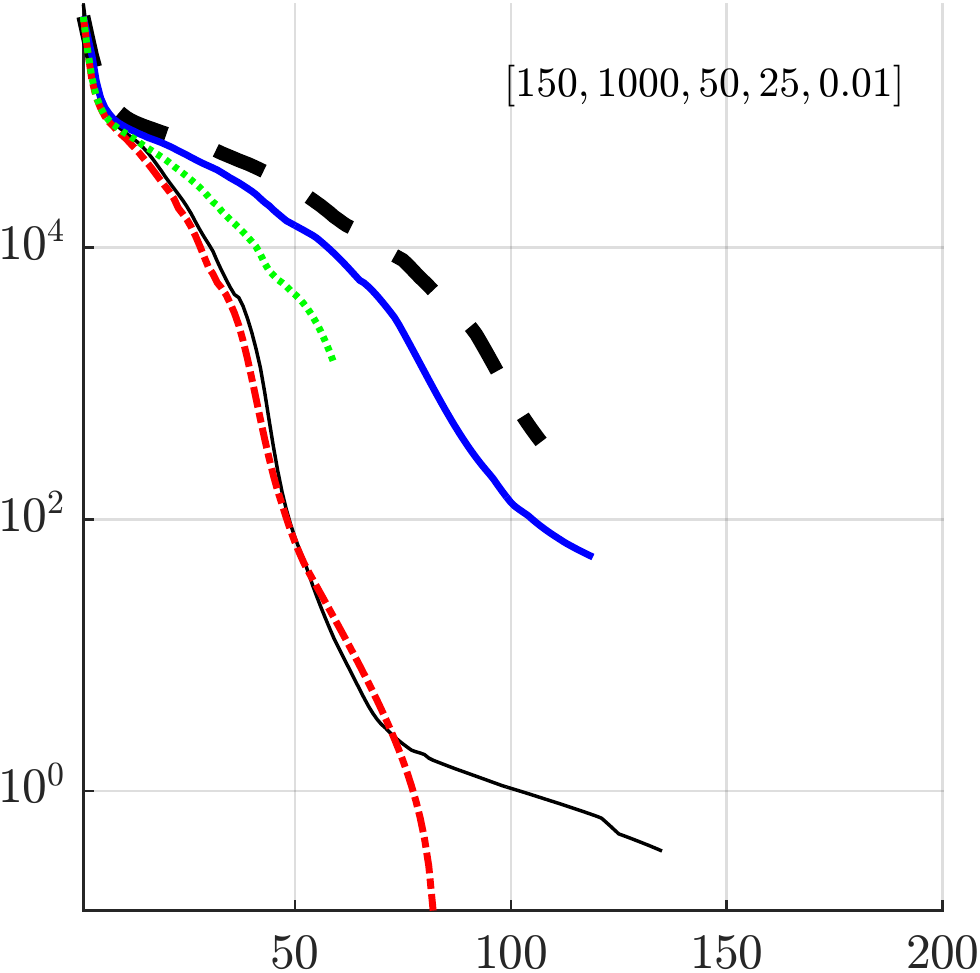}  
\end{tabular}
\caption{Summary of the results of Fig.\ref{fig:exp:ahals_BH50_50_50_10_n0} to Fig.\ref{fig:exp:nest_BH150_1000_50_25_n001}.}
\end{table}

\clearpage

As BCD and LS-BCD consistently perform slower than HER-BCD, we will not show the results on comparing HER-BCD with them in the subsequent plots.

% BH- AOADMM
% BH- AO Nest
%
% HSI : PaviaU(Left) IndianPines(Right)
%
\begin{center}\begin{figure}[ht!]
		\centering
		\begin{subfigure}[b]{0.49\textwidth}
			\centering
			\includegraphics[width=\linewidth]{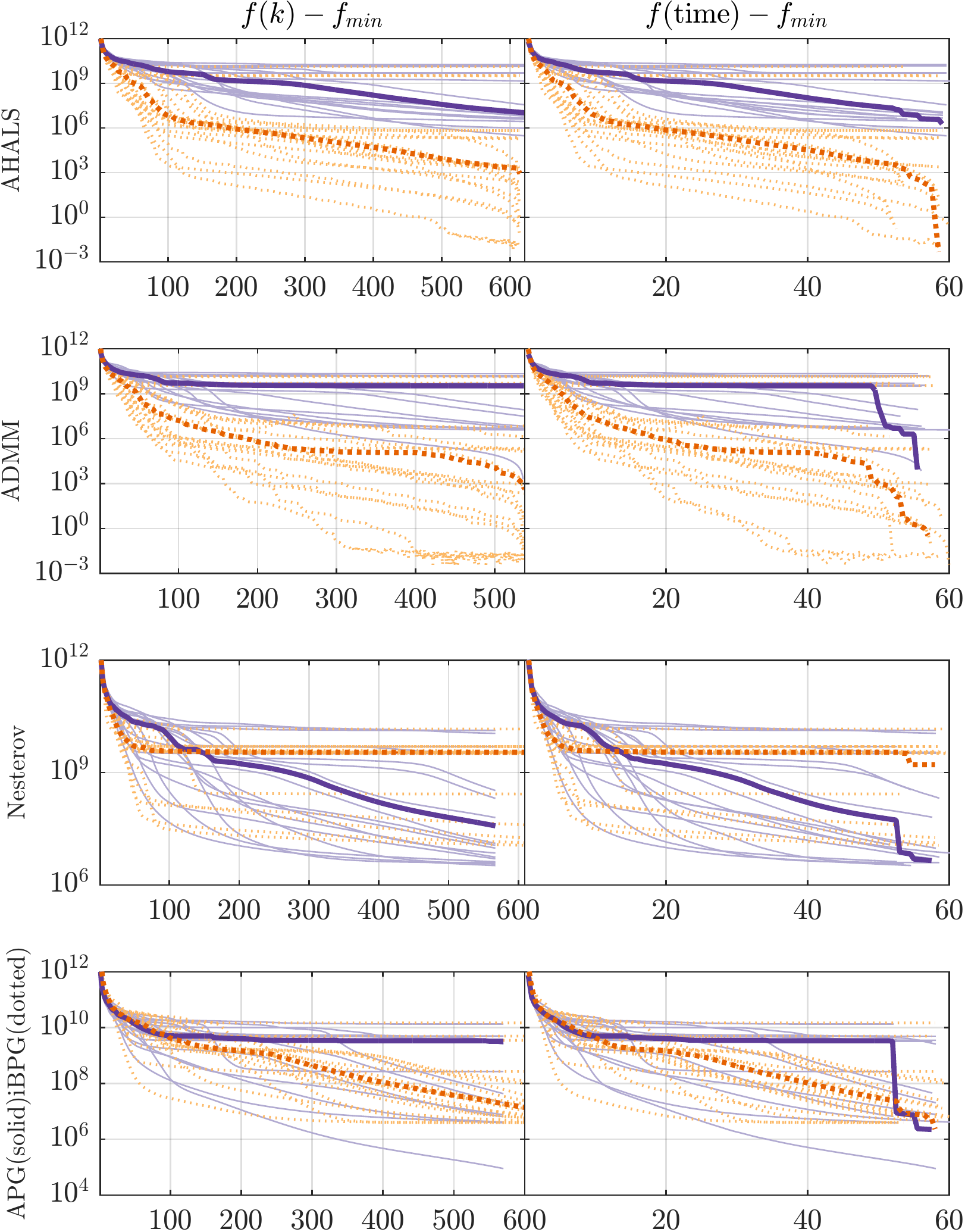}
			\caption{On PaviaU dataset.}
			\label{fig:exp:real_PaviaU10_supp}
		\end{subfigure}
		\hfill
		\begin{subfigure}[b]{0.5\textwidth}
			\centering
			\includegraphics[width=\linewidth]{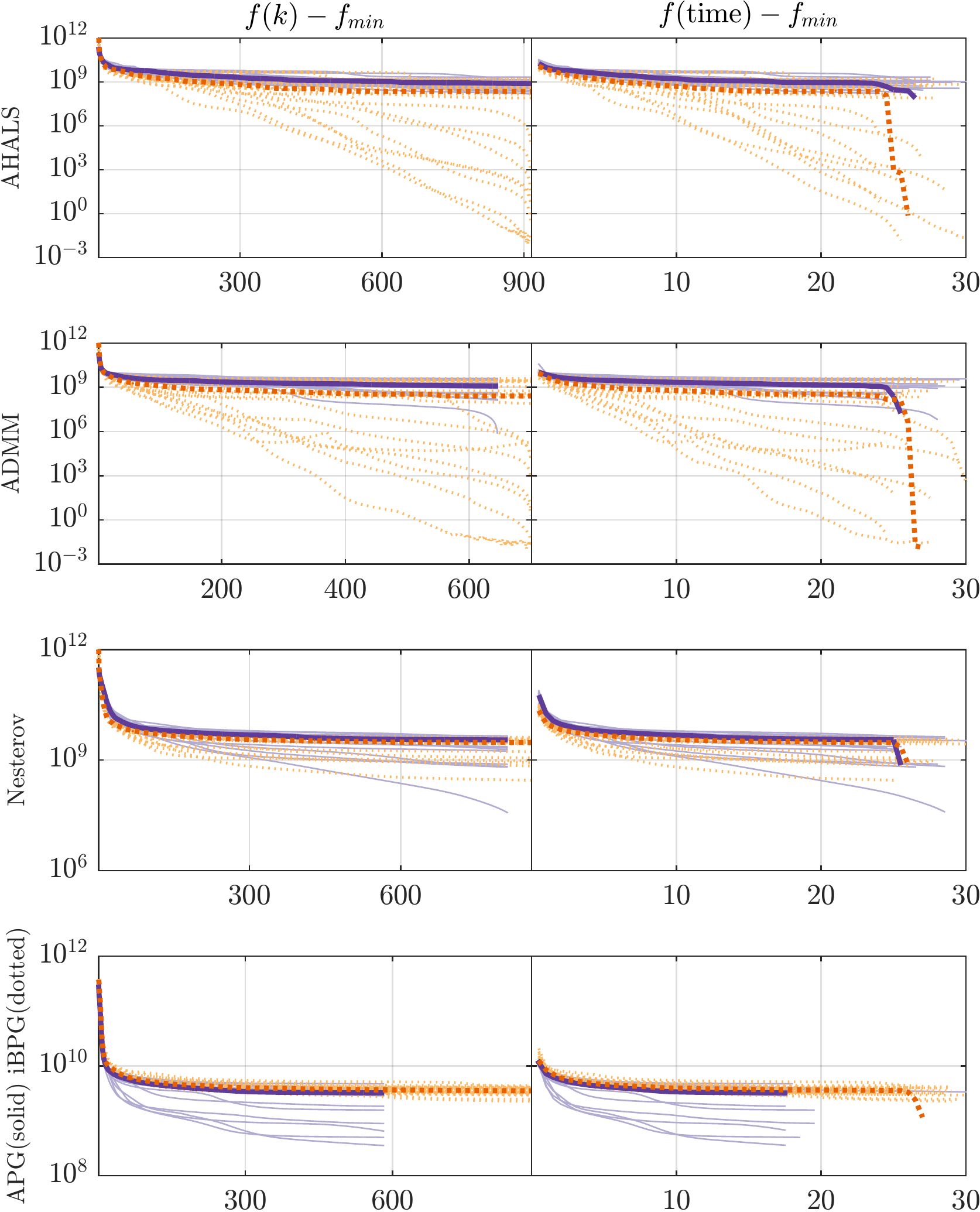}
			\caption{On Indian Pines dataset.}
			\label{fig:exp:real_IndianPines15_supp}
		\end{subfigure}
		\caption{Experiments with HSI data.
	    Gradient-based methods (Nesterov, HER-Nesterov, APG and iBPG)perform worse than AHALS and ADMM.
        In terms of $f$, the best run of HER-ADMM and HER-AHALS are about $10^7 - 10^9$ times better than AO-Nesterov, HER-Nesterov, APG and iBPG.
        In terms of $e$, the best run for HER-ADMM and HER-AHALS are about $10^{10} - 10^{12}$ times better than Nesterov, APG and iBPG.
        }
		\label{fig:exp:HSI_supp}
\end{figure}\end{center}

\end{document}